\documentclass[12pt,a4paper,leqno]{amsart}
\usepackage{amsmath,amsthm,amsfonts,amssymb}
\newtheorem{defin}{Definition}[section]
\newtheorem{lemma}[defin]{Lemma}
\newtheorem{propos}[defin]{Proposition}
\newtheorem{defprop}[defin]{Definition and Proposition}
\newtheorem{theor}[defin]{Theorem}
\newtheorem{corol}[defin]{Corollary}
\newtheorem{example}[defin]{Example}
\newenvironment{definition}{\begin{defin} \rm}{\hspace{2mm}\BOX \end{defin}}
\newenvironment{definition*}{\begin{defin} \rm}{\end{defin}}
\newcommand{\dln}[5]{\begin{picture}(#4,#5)(0,0)\multiput(0,0)(#1,#2){#3}{{\tiny .}}\end{picture}}

\newcommand{\bd}{\begin{definition}}
\newcommand{\ed}{\end{definition}}
\newcommand{\bdefn}{\begin{definition*}}
\newcommand{\edefn}{\end{definition*}}
\newcommand{\bl}{\begin{lemma}}
\newcommand{\el}{\end{lemma}}
\newcommand{\bp}{\begin{propos}}
\newcommand{\ep}{\end{propos}}
\newcommand{\bt}{\begin{theor}}
\newcommand{\et}{\end{theor}}
\newcommand{\bc}{\begin{corol}}
\newcommand{\ec}{\end{corol}}
\newcommand{\bds}{\begin{displaystyle}}
\newcommand{\eds}{\end{displaystyle}}
\newcommand{\bdm}{\begin{displaymath}}
\newcommand{\edm}{\end{displaymath}}
\renewcommand{\leq}{\leqslant}
\renewcommand{\geq}{\geqslant}
\newcommand{\BOX}{\rule{2mm}{2mm}}
\newcommand{\Zint}{\mathbb Z}     % Integer number field
\newcommand{\Rea}{\mathbb R}      % Real number field
\newcommand{\Cplx}{\mathbb C}     % Complex  number field
\newcommand{\Nat}{\mathbb N} 
\newcommand{\gf}{\hat{{\mathfrak g}}}
\newcommand{\gl}{\mathfrak {gl}}      % Lie algebra gl
\newcommand{\sll}{\mathfrak {sl}}      % Lie algebra sl
\newcommand{\agl}{\widehat{{\mathfrak {gl}}}}      % affine Lie algebra gl
\newcommand{\asll}{\widehat{{\mathfrak {sl}}}}     % affine Lie algebra sl
\newcommand{\daha}{{\mathcal H}}               % double degenerate affine Hecke algebra
\newcommand{\tf}{\mathfrak t}                % Cartan subalgebra of affine gl
\renewcommand{\ep}{\epsilon}
\newcommand{\vep}{\varepsilon}
\newcommand{\vro}{\varrho}
\newcommand{\ph}{\varphi}
\newcommand{\al}{\alpha}
\newcommand{\vr}{\mathrm v}
\newcommand{\row}{{\mathrm {row}}}
\newcommand{\col}{{\mathrm {col}}}
\newcommand{\reg}{{\mathrm {reg}}}
\newcommand{\Vaff}{V_{\mathrm {aff}}}
\newcommand{\Ve}{\mathfrak V}
\newcommand{\bef}{\mathfrak b}
\newcommand{\UU}{\operatorname{U}}
\newcommand{\YY}{\operatorname{Y}}
\newcommand{\Ch}{\operatorname{Ch}}

\newcommand{\hv}{h^{\vac}}
\newcommand{\vac}{\mathrm {vac}}
\newcommand{\vf}{\mathfrak v}
\newcommand{\uf}{\mathfrak u}
\newcommand{\End}{\mathrm {End}}
\newcommand{\Ker}{\mathrm {Ker}}
\newcommand{\Imm}{\mathrm {Im}}
\newcommand{\ov}{\overline}
\newcommand{\un}{\underline}
\newcommand{\F}{{\mathcal F}}

\newcommand{\bb}{b}

\newcommand{\fra}[1]{
\begin{picture}(8,8)(2,1.5)
\multiput(0,0)(0,8){2}{\line(1,0){8}}
\multiput(0,0)(8,0){2}{\line(0,1){8}}
\put(0,0){\makebox(8,8){{\tiny $ #1 $}}}
\end{picture}}

\newcommand{\semidir}{\begin{picture}(20,50)(0,-2.5)
\put(7,1){\line(0,1){20}}
\put(1,1){$\times$}
\end{picture}}

\begin{document}
%\thispagestyle{empty}
%%%%%%
\title[ ]{ Yangian actions on higher level irreducible integrable modules of $\agl_N$ }

\author[ ]{Denis Uglov}

\address{Research Institute for Mathematical Sciences, Kyoto University, 606 Kyoto, Japan.}

\email{duglov@kurims.kyoto-u.ac.jp}

\subjclass{Primary 17B67; Secondary 17B37, 17B81}
%
% Abstract
%
\begin{abstract}
An action of the Yangian of the general linear Lie algebra $\gl_N$ is defined on every irreducible integrable highest weight module of $\agl_N$ with level $\geq 2.$ This action is derived, by means of the Drinfeld duality and a subsequent semi-infinite limit, from a certain induced representation of the degenerate double affine Hecke algebra $\daha.$ 
Each vacuum module of $\agl_N$ is decomposed into irreducible Yangian representations by means of the intertwiners of $\daha.$ Components of this decomposition are parameterized by semi-infinite skew Young diagrams.
\end{abstract}

\maketitle

\tableofcontents
\section{Introduction}
It is  known that irreducible integrable level 1 modules of the affine Lie algebras $\asll_N$ and $\agl_N$ possess hidden Yangian symmetries \cite{HHTBP,Schoutens}. Apart from the standard actions of $\asll_N$ or $\agl_N$ such modules  admit also actions of the associative algebra $\YY(\gl_N)$ -- the Yangian of $\gl_N.$ 
These Yangian actions have a physical interpretation: the family of operators arising from the centre of $\YY(\gl_N)$ is interpreted as the algebra of conserved charges in a long-range interacting solvable model -- the Haldane-Shastry model in the case of $\asll_N$ and the Sutherland model with spin in the  case of $\agl_N.$ 
The irreducible integrable level 1 modules of the quantum affine algebras $\UU_q(\,\asll_N)$ and $\UU_q(\,\agl_N)$ possess hidden symmetries as well, the place  of the Yangian is now taken by a level 0 action of $\UU_q(\,\agl_N).$ In fact, every 
irreducible integrable level 1 module  $\UU_q(\,\agl_N)$ is a module of the quantum toroidal algebra \cite{GKV}. The two actions of $\UU_q(\,\agl_N)$ arise as actions of the two subalgebras -- each isomorphic to $\UU_q(\,\agl_N)$ --  of this large algebra.

Based on the combinatorial results of \cite{ANOT,BLS,KKN1,KKN2} one can expect, that higher level irreducible integrable modules of $\asll_N$ (at least for $N=2$) and  $\agl_N$ possess hidden Yangian symmetries as well. In the present paper we study the Yangian symmetries of the irreducible integrable modules of $\agl_N.$

\subsection{Statement of the results} \label{se:res}
Let $\gf$ be the affine Lie algebra $\asll_N,$ let $\Lambda$ be a dominant integral weight of $\gf$  of level $L \geq 2,$ and let $V(\Lambda)$ be the irreducible $\gf$-module with highest weight $\Lambda.$ Let $H$ be the infinite-dimensional Heisenberg Lie algebra generated by the elements $K,B(m)$ $(m \in \Zint_{\neq 0})$ modulo the relations 
\bdm
 [ B(m) , B(n) ] = K m \delta_{m+n,0}, \quad [ K , B(m)] = 0, 
\edm
and let $S_- = \Cplx[B(-1),B(-2),B(-3),\dots]$ be the standard bosonic  Fock module of $H$ where the central charge $K$ acts as the multiplication by $NL.$  The  object which we study in this paper is the irreducible module 
\begin{equation}
S_-\otimes V(\Lambda)  \label{eq:mainspace}
\end{equation}   
of the Lie algebra $\agl_N = H\oplus \asll_N.$

The Yangian $\YY(\gl_N)$ of the general linear Lie algebra $\gl_N$ is a canonical deformation of the universal enveloping algebra $\UU(\gl_N[u])$ in the class of Hopf algebras \cite{DrinfeldYang}. 
The centre of $\YY(\gl_N)$ is isomorphic to the algebra of polynomials in infinitely many variables $\Delta_1,\Delta_2,\Delta_3,\dots$ which it is customary to collect into a single formal series $\Delta(u) = 1 + u^{-1} \Delta_1 + u^{-2} \Delta_2 + u^{-3} \Delta_3 + \cdots $ called the {\em quantum determinant} of $\YY(\gl_N).$  The quotient of $\YY(\gl_N)$ by the ideal generated by the centre, denoted $\YY(\sll_N),$ is known as the Yangian of $\sll_N.$ The algebra  $\YY(\gl_N)$ is isomorphic to a tensor product of its centre and $\YY(\sll_N)$ (an extensive review of the Yangian can be found in \cite{MNO}).

As the main result of this article we define an action of the algebra $\YY(\gl_N)$ on each of the irreducible integrable $\agl_N$-modules of  the form (\ref{eq:mainspace}). 
With this done, we proceed to study the decomposition of $S_-\otimes V(\Lambda)$ as the Yangian module in the case where $S_-\otimes V(\Lambda)$  is the vacuum representation of $\agl_N,$ i.e. $\Lambda = L\Lambda_0.$ 
We find, that each component of this decomposition is an irreducible finite-dimensional $\YY(\gl_N)$-module. Yangian modules that appear as irreducible components of $S_-\otimes V(L\Lambda_0)$ form a set parameterized by semi-infinite skew Young diagrams of a special type which we will now describe. With any semi-infinite sequence of non-negative integers $\vec{h} = (h_i)_{i\in \Nat}$ one associates the semi-infinite skew diagram 
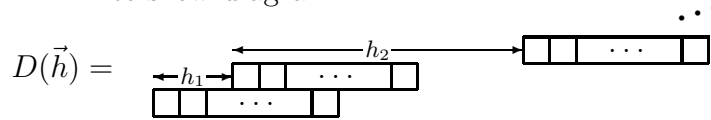
\begin{equation*}
D(\vec{h}) = \quad 
\begin{picture}(210,20)(0,15) 
{\thicklines
\put(0,0){\line(1,0){70}}\put(0,10){\line(1,0){100}} 
\put(30,20){\line(1,0){70}}
\multiput(140,20)(0,10){2}{\line(1,0){70}}
%%%%%%%%%%
\multiput(0,0)(10,0){3}{\line(0,1){10}} \multiput(60,0)(10,0){2}{\line(0,1){10}} 
\multiput(30,10)(10,0){3}{\line(0,1){10}} \multiput(90,10)(10,0){2}{\line(0,1){10}} 
\multiput(140,20)(10,0){3}{\line(0,1){10}} \multiput(200,20)(10,0){2}{\line(0,1){10.5}} 
}
%%%%%%%%%%
\put(9,15){\vector(-1,0){9}}\put(20,15){\vector(1,0){10}}
\put(79,25){\vector(-1,0){49}}\put(90,25){\vector(1,0){50}}
\put(10,10){{\makebox(10,10){{\scriptsize $h_1$}}}}\put(80,20){{\makebox(10,10){ {\scriptsize $h_2$} }}}
\put(20,0){\makebox(40,10){$\dots$}}\put(50,10){\makebox(40,10){$\dots$}}
\put(160,20){\makebox(40,10){$\dots$}}
\multiput(200,37)(6,2){3}{\circle*{2}}
\end{picture}
\end{equation*}
\hfill \\ 
where each row has $L$ squares. The diagram that corresponds to $\vec{h}^{\vac}= (\hv_i)_{i\in \Nat}$ where  $\hv_i = L$ if $i\equiv0 \bmod N$ and  $\hv_i = 0$ if otherwise, is called the {\em vacuum diagram.}    
\begin{defin} \label{def:siskewD}
A sequence $\vec{h} \in \Zint_{\geq 0}^{\Nat},$ or the corresponding $D(\vec{h}),$ is called a {\em skew diagram of type } ${\mathcal D}_L(N)$ if and only if:  {\em (i)} $ h_i = \hv_i$ for all but finite number of $i,$ {\em (ii)} each column of  $D(\vec{h})$ contains no more than $N$ squares. 
\end{defin}

\noindent An irreducible module $V$ of $\YY(\gl_N)$ is specified up to an isomorphism by a set of $N-1$ monic polynomials $P_1(u),\dots,P_{N-1}(u) \in \Cplx[u]$ called the {\em Drinfeld polynomials} of $V,$ and a formal series $f(u) \in 1 + u^{-1}\Cplx[[u^{-1}]]$ such that $\Delta(u)|_V = f(u)\cdot 1_{V}$  \cite{DrinfeldYang}.  
Accordingly, to each skew Young diagram $\vec{h}$ of type ${\mathcal D}_L(N)$ we attach $N-1$ polynomials $P_1^{\vec{h}}(u),\dots ,P_{N-1}^{\vec{h}}(u)$ and a formal series $f^{\vec{h}}(u).$ 
To define these, number the rows and columns of the plane where $D(\vec{h})$ is positioned by $\Zint$ in the matrix  order: rows -- from the top downward, columns -- from the left to the right so that the leftmost square in the first row of  $D(\vec{h})$ has the vertical coordinate $0$ and the horizontal coordinate $1+\sum_{i=1}^{\infty} \hv_i - h_i.$  To the square positioned at the intersection of $i$th row and $j$th column  one then associates the number $j-i$ called the {\em content} of this square. For all $k=1,\dots,N-1$ define\begin{equation*}
P_k^{\vec{h}}(u) = \prod_c ( u - k - c + 1) %\label{eq:DP}
\end{equation*}
where the product is taken over contents of bottom squares in columns of height $k$ in $D(\vec{h})$ (note that for any diagram of type ${\mathcal D}_L(N)$ there are only a finite number of such columns). With each skew diagram $D(\vec{h})$ of type ${\mathcal D}_L(N)$ associate the rational function   
\begin{equation}
f^{\vec{h}}(u) = \prod_{k=1}^{\infty}  \left(\frac{ u + 2 - r_1  - k - \sum_{i=1}^{k-1} h_i }{ u + 1 - k - \sum_{i=1}^{k-1} \hv_i  }\right)\left(\frac{ u + 1 - k - \sum_{i=1}^{k-1} \hv_i - L }{ u + 2 - r_1  - k - \sum_{i=1}^{k-1} h_i -L}\right)
\label{eq:QD1}
\end{equation}
where $r_1 = 1+\sum_{i=1}^{\infty} \hv_i - h_i.$ The following theorem gives the Yangian decomposition of the vacuum $\agl_N$-module. 
\begin{theor} \label{t:Yangdec}
Let $V_{{D}(\vec{h})}$ be the irreducible $\YY(\gl_N)$-module defined up to an isomorphism by the Drinfeld polynomials $P_1^{\vec{h}}(u),\dots,P_{N-1}^{\vec{h}}(u)$ and the formal series obtained by the expansion of  
$f^{\vec{h}}(u) 
$
in negative powers of $u.$

Then as the $\YY(\gl_N)$-module the vacuum representation $S_-\otimes V(L\Lambda_0)$  is isomorphic to the direct sum 
$
\bigoplus_{ D(\vec{h}) }V_{{D}(\vec{h})} $
 taken over all semi-infinite skew diagrams of type ${\mathcal D}_L(N).$
\end{theor}
\noindent The irreducible $\YY(\gl_N)$-modules appearing in the above decomposition belong to the class of {\em tame} modules. These are distinguished among all  $\YY(\gl_N)$-modules in that they admit semi-simple actions of a certain maximal commutative subalgebra of $\YY(\gl_N)$ \cite{C1,NT1,NT2}. The structure of a tame irreducible module is considerably simpler than that of a general  irreducible $\YY(\gl_N)$-module. In particular, the $\sll_N$-character of a tame irreducible module admits a simple explicit form given by a  product of skew Schur functions \cite{NT1,NT2,KKN1}. No explicit formula, on the other hand, seems to be at the present available for the character of a general irreducible  $\YY(\gl_N)$-module apart from the case $N=2.$  

For each skew diagram $D(\vec{h})$ of type ${\mathcal D}_L(N)$ the {\em degree} of $D(\vec{h})$ is defined as $$\deg D(\vec{h}) = \sum_{i=1}^{\infty} i(h_i - \hv_i).$$ Note that $ h_i = \hv_i $ if $ i > N \deg D(\vec{h}).$ A {\em finite part} $\ov{\!D}(\vec{h})$ is  the finite skew Young diagram obtained from   $D(\vec{h})$ by deleting all but the first $Nl$ rows, where $l$ is an arbitrary integer such that $l \geq \deg D(\vec{h}).$ 
The module $V_{{D}(\vec{h})}$ is isomorphic, up to a $\YY(\gl_N)$-automorphism, to one of the  so called {\em elementary} tame  $\YY(\gl_N)$-modules associated with an arbitrary finite part of $D(\vec{h})$ (see \cite{NT2} and Section \ref{sec:yangdecofvac}). The $\sll_N$-character of  $V_{{D}(\vec{h})}$ equals the skew Schur function $s_{\ov{\!D}(\vec{h})}(x_1,\dots,x_N)$ (cf. \cite{Mac} and Section \ref{sec:yangdecofvac}) labeled by an arbitrary finite part  $\ov{\!D}(\vec{h})$ (elementary modules labeled by different finite parts of the same semi-infinite diagram are isomorphic to each other). In particular, the dimension of   $V_{{D}(\vec{h})}$ is equal to the number of semi-standard tableaux on numbers $1,\dots,N$ of the shape $D(\vec{h}).$ 
The Yangian decomposition given in Theorem \ref{t:Yangdec} is reflected in the following expression for the character of the vacuum $\agl_N$-module: 
\bdm
 \Ch_{S_-\otimes V(L\Lambda_0)}(q,x) = \frac{\Ch_{V(L\Lambda_0)}(q,x)}{\prod_{k\geq 1}(1 - q^k)}= \sum_{D(\vec{h})} \: q^{\sum_{i=1}^{\infty} i(h_i - \hv_i)} \: s_{\ov{\!D}(\vec{h})}(x)
\edm 
where the summation is taken over all semi-infinite skew diagrams of type ${\mathcal D}_L(N)$ and $x = (x_1,\dots,x_N).$ Above,  $\ov{\!D}(\vec{h})$ stands for an arbitrary finite part of $D(\vec{h}).$

\subsection{Related results}
\subsubsection{Long-range solvable models and Yangian actions on level 1 irreducible integrable modules of $\asll_N$ and $\agl_N$} The subject of Yangian actions on infinite-dimensional representations of affine Lie algebras was initiated in the  work \cite{HHTBP}, where actions of $\YY(\gl_2)$ were defined on irreducible integrable highest weight modules of $\asll_2$ of level 1. 
The physical meaning of a $\YY(\gl_2)$-action of this type is as follows. 
The centre of $\YY(\gl_2)$ is interpreted as the algebra of commuting conserved charges in an appropriate field theory limit of the solvable Haldane-Shastry spin chain with long-range interaction, whereas  the subalgebra $\YY(\sll_2)$ is interpreted as the  non-abelian symmetry algebra which commutes with the conserved charges and, thus,  gives rise to degeneracies in their spectra. 
The results of  \cite{HHTBP} were subsequently extended and generalized in several directions. 
First of all, the  decompositions of irreducible integrable level 1 modules of $\asll_2$ with respect to the   $\YY(\gl_2)$-actions were obtained in \cite{BPS} and \cite{BLS}. These decompositions are given in terms of the {\em spinon} bases. In physical terms these bases are formed by elementary excitations of the Haldane-Shastry spin chain -- the spinons. 
\mbox{}From the viewpoint of the representation theory these bases provide a new description of the irreducible modules of the affine Lie algebra which, in particular, leads to new identities for  characters of these modules. The above  results  are known to generalize for the $\sll_N$ with $N>2$ (cf. \cite{Schoutens, BS}).       

 Irreducible integrable level 1 modules of the Lie algebra $\agl_N$ also admit actions of the Yangian $\YY(\gl_N)$ (\cite{U1,U2}). An irreducible  $\agl_N$-module is  simply a tensor product of an irreducible $\asll_N$-module and the bosonic Fock module of the Heisenberg algebra. In this case the physical meaning of the  $\YY(\gl_N)$-action changes in that the centre of $\YY(\gl_N)$ is now interpreted as the algebra of conserved charges in a suitable field theory limit of the Sutherland model with spin. 
The Sutherland model is the parent model of the Haldane-Shastry spin chain, and the latter is the limiting case of the former in which the dynamical degrees of freedom are frozen so that only the spin degrees of freedom remain relevant. In the language of the representation theory this freezing of dynamical degrees of freedom is expressed simply as the fact that an irreducible  $\asll_N$-module is the quotient of  an irreducible  $\agl_N$-module by the linear subspace generated by the creation operators of the Heisenberg algebra.    

The present work is partly an attempt to extend the above results to the higher level irreducible integrable modules of affine Lie algebras. What we obtain here is quite similar to the level 1 situation. The centre of the $\YY(\gl_N)$-action gives rise to a family of commuting Hamiltonians. 
The spectral decomposition of the vacuum level $L$ module of $\agl_N$ by these Hamiltonians is given by Theorem \ref{t:Yangdec}, the spectrum itself being encoded in the coefficients of the series obtained by the expansion of (\ref{eq:QD1}). In contrast with the level 1 case, however, we lack a physical interpretation  of this picture. Presumably, the above commuting Hamiltonians can be related to conserved charges of some suitable higher spin version of the Sutherland model with spin.
\nopagebreak
\subsubsection{Spectral decomposition of solvable lattice models}
The corner transfer matrix method in solvable lattice models leads to a combinatorial description of an irreducible integrable highest weight module of $\asll_N$ as the space of paths \cite{J}. Each path is an eigenvector of an infinite family of mutually commutative Hamiltonians defined by the local energy function associated with the $R$-matrix of the solvable lattice model. 
A series of recent papers \cite{ANOT,KKN1,KKN2} deals with the spectral decomposition of the path space with respect to this commutative family. Remarkably, it was found that in the case when level is 1 this spectral decomposition matches combinatorially the Yangian decomposition described in the previous section. More precisely, there is a one-to-one, $\sll_N$-character preserving, correspondence between eigenspaces of the commuting Hamiltonians in the path space, and irreducible components of the Yangian action on the $\asll_N$-module that is described by this path space. 
So far, Yangian actions have been defined only on level 1 irreducible integrable  $\asll_N$-modules (see, however, \cite{BMM} where Yangian actions of a different type from the one we consider here were defined on higher level modules), and it is not known, therefore, whether a similar correspondence, or an appropriate modification thereof,  exists in the case of  higher levels. Still, it is instructive to compare the Yangian decomposition of the irreducible $\agl_N$-module $S_-\otimes V(L\Lambda_0)$ given by Theorem \ref{t:Yangdec} with the spectral decomposition of the path space for the vacuum $\asll_N$-module $V(L\Lambda_0)$ given in \cite{KKN2}. The  components of these decompositions are labeled similarly, by semi-infinite skew Young diagrams. 
The one difference is that in the case of $V(L\Lambda_0)$ the sequence $(h_i )_{i\in \Nat},$ defining a skew Young diagram, contains only numbers not exceeding $L.$ The other difference is in the dimension of the component labeled by a skew Young diagram. For $S_-\otimes V(L\Lambda_0)$ this dimension equals the number of semi-standard tableaux (on the numbers $1,\dots,N$) of shape defined by the skew diagram. Whereas for  $V(L\Lambda_0)$ the dimension is generally smaller -- it equals the number of the so called non-movable semi-standard tableaux \cite{KKN2}.

\subsubsection{Representations of the Quantum Toroidal Algebra} 
The results of the present paper have natural counterparts for the quantum affine algebras. 
In this case instead of a Yangian action on the irreducible $\agl_N$-module $S_-\otimes V(\Lambda)$ we obtain a level 0 action of $\UU_q(\,\agl_N) = \UU(H)\otimes \UU_q(\,\asll_N)$ on the irreducible $\UU_q(\,\agl_N)$-module $S_-\otimes V_q(\Lambda),$ where $V_q(\Lambda)$ is the irreducible module of $\UU_q(\,\asll_N)$ with highest weight $\Lambda.$ The two actions of  $\UU_q(\,\agl_N)$ on $S_-\otimes V_q(\Lambda)$ extend, in fact, to a representation of a more general object -- the quantum toroidal algebra. 
The quantum toroidal algebra is a deformation of the enveloping algebra of the universal central extension of the Lie algebra of maps from the two-dimensional torus into $\sll_N.$ Some infinite-dimensional representations of the quantum toroidal algebra on level 1 modules of $\UU_q(\,\agl_N)$ and $\UU_q(\,\asll_N)$ were recently constructed and studied in \cite{Saito,STU,T,VV}. Our present work has mostly been motivated by a desire to generalize the latter results on the case of higher levels. We plan to report on this issue in a forthcoming article.

\subsection{An overview of the paper}
The paper is divided into three main sections. 
Section 2 is a summary of definitions and known results on the Degenerate Double Affine Hecke Algebra, representation theory of $\agl_N$ and $\asll_N$ (mainly we cover the Fock space module of $\agl_N$ and irreducible quotients thereof) and the Yangian. 
In Section 3 we define a Yangian action on an irreducible $\agl_N$-module of the form (\ref{eq:mainspace}). 
Section 4 contains the results on the Yangian decompositions of the vacuum $\agl_N$-modules, in particular --  a proof of Theorem \ref{t:Yangdec}. Below we highlight the main technical points of these sections.     

\subsubsection{Degenerate Double Affine Hecke Algebra}
The central role in the construction of the Yangian actions on the $\agl_N$-modules is played by the Degenerate Double Affine Hecke Algebra. We start from a certain family of parabolic induced  representations of this algebra. These representations were defined in the setting of Double Affine Hecke Algebra by I.Cherednik \cite{C3}, and were extensively studied in the recent work of Arakawa {\em et al.} \cite{AST}. 
The Drinfeld functor \cite{DrinfeldDual} applied to a representation of this type gives a Yangian action on the finite exterior product (wedge product) of vector spaces $\Vaff.$ Here $\Vaff = \Cplx[z^{\pm 1}] \otimes \Cplx^N\otimes \Cplx^L$ is the level 0 module of the Lie algebra $\asll_N\oplus \asll_L$ obtained by the affinization of the vector module $\Cplx^N\otimes \Cplx^L$ of $\sll_N\oplus\sll_L.$ 
The wedge product with infinite number of factors and appropriate asymptotic conditions is nothing but the well-known fermionic Fock space module of the Clifford algebra \cite{J}. 
The Yangian action on the finite wedge product gives rise to a Yangian action on the Fock space. Informally speaking, the latter action is an infinite limit of the former. 
A formal definition is contained in Section \ref{se:YonF}. This definition is made possible by a certain stability (cf. Section \ref{sec:inter-and-stab}) of the Yangian action on the finite wedge product when the number of factors in the latter grows by steps of $NL.$  

\subsubsection{From the Fock space to the irreducible modules of $\agl_N$}
The finite wedge product is equipped with the diagonal, level 0 action of the Lie algebra $\asll_N\oplus \asll_L$ (i.e. actions of both subalgebras $\asll_N$ and $\asll_L$ have level 0). 
This action survives in the Fock space limit. On the Fock space, however, the action of $\asll_N$ acquires the level $L$ and the  action of $\asll_L$ acquires the level $N.$ 
Moreover, on the Fock space, there is an action of the Heisenberg algebra $H$ that centralizes the action of $\asll_N\oplus \asll_L.$ The decomposition of the Fock space with respect to the action of $H\oplus \asll_N\oplus \asll_L$ was studied by I.Frenkel \cite{F}. 
His result is quoted in Theorem \ref{t:Fockdec}. 
\mbox{}Form this result it follows that every $\agl_N$-module of the form (\ref{eq:mainspace}) is realized as the quotient of a charge component of the Fock space (cf. Section \ref{sec:wedge}) by the linear subspace generated by a certain subalgebra of $\UU(\asll_L).$  With a specific choice of parameters, the $\YY(\gl_N)$-action leaves this subspace invariant, and therefore a $\YY(\gl_N)$-action is defined on (\ref{eq:mainspace}).

\subsubsection{Intertwiners of the Degenerate Double Affine Hecke Algebra and the Yangian decompositions of the vacuum modules of $\agl_N$}
The intertwiners of weight spaces  provide a powerful machinery which allows to study representations  of the Degenerate Double Affine Hecke Algebra, $\daha$. These intertwiners form a linearly independent family of elements of $\daha$ parameterized by a subset of the affine Weyl group. 
The induced representations of the Degenerate Double Affine Hecke Algebra analyzed in \cite{AST} depend on several parameters. If these parameters are generic, the corresponding representation is irreducible and admits a basis generated by the intertwiners from the distinguished cyclic vector called the highest weight vector of this representation (\cite{AST} Proposition 2.4.3 and Theorem 2.4.4). 
The parameters of the representations which are used in the present paper to define the Yangian actions on the irreducible modules of  $\agl_N$ are not generic. This means that there are  intertwiners which are not invertible on the highest weight vector, and/or that intertwiners applied to the highest weight vector no longer form a basis of the representation.  

In Section \ref{se:regel} we study intertwiners in the representations of $\daha$ that are used to define the Yangian actions on the vacuum modules of  $\agl_N.$ 
We introduce a subset of the affine Weyl group formed by {\em regular} elements. These elements are defined so  that the corresponding intertwiners generate from the highest weight vector a  linearly independent set of vectors and, moreover,  are invertible on this vector. 
We find that the set of regular elements is in one-to-one correspondence with the set of skew Young diagrams of a special type which we refer to as the type ${\mathcal D}_L^m$ (Section \ref{sec:regdiag}). 
To each skew Young diagram of this type, with the additional restriction that the number of squares in every column does not exceed $N,$ there corresponds a highest weight vector of the $\YY(\gl_N)$-action on the wedge product with the number of factors equal to $Lm$ (Section \ref{sec:Yhwinwedge}). 
When the number of factors becomes infinite, these highest weight vectors give rise to a family of highest weight vectors of the Yangian action on the Fock space. 
The members of this family are now labeled by the semi-infinite skew Young diagrams of type ${\mathcal D}_L(N)$ introduced in Definition \ref{def:siskewD}. In Section \ref{sec:yangdecofvac} we demonstrate that these highest weight vectors exhaust all Yangian highest weight vectors in the quotient of the suitable charge component of the Fock space that is isomorphic to the vacuum module of $\agl_N$ of level $L.$         
\subsection{Acknowledgments}
I am grateful to I. Cherednik, B. Feigin, T. Miwa, T. Suzuki and K. Takemura for numerous valuable discussions concerning the subject of this article.

\section{Preliminaries} 

\subsection{Representations of the Affine Lie Algebra $\asll_N$}
\subsubsection{Affine Lie Algebra $\asll_N$} \label{sec:sl}
Let $N$ be an integer $\geq 2$ and let ${\mathfrak h}_N$ be an $N$-dimensional vector space over $\Cplx$ with basis $\{ H_0,H_1,\dots,H_{N-1}\}.$ We let $\{ \Lambda_0,\Lambda_1,\dots,\Lambda_{N-1}\}$ be the corresponding dual basis of  ${\mathfrak h}_N^*,$ the dual space of ${\mathfrak h}_N.$ It is convenient to extend the index set so that $ \Lambda_s = \Lambda_{(s\bmod N)}$ for all $s \in \Zint.$ Then, for all $s\in \Zint,$ we set $\vep_s = \Lambda_{s} - \Lambda_{s-1}$ and $\tau_{s} = 2 \Lambda_s - \Lambda_{s-1}- \Lambda_{s+1}.$    

The $N\times N$ matrix $\| \tau_s(H_t) \|$ is called the generalized Cartan matrix of type $A_{N-1}^{(1)}.$ The associated affine Kac-Moody algebra is denoted $\asll_N.$ It is defined as the complex Lie algebra generated by elements $ E_s,F_s,H_s$ for $ 0\leq s < N,$ subject to the relations: 
\begin{alignat*}{4}
 & [H_s,H_t] = 0; & \quad & [E_s,F_t] = \delta_{st} H_s;\\ 
 & [H_s,E_t] = \tau_t(H_s) E_t; &  \quad & [H_s,F_t] = -\tau_t(H_s) F_t;  \\
 & ({\mathrm {ad}} E_s)^{1-\tau_t(H_s)} E_t = 0; & \quad & ({\mathrm {ad}} F_s)^{1-\tau_t(H_s)} F_t = 0 \quad (s \neq t),
\end{alignat*}
where $({\mathrm {ad}}a ) b = [a,b].$ 
 The abelian Lie algebra ${\mathfrak h}_N$ is known as the Cartan subalgebra of  $\asll_N,$  and the vectors  $\Lambda_s$ as the fundamental weights.    

Let $\sll_N$ be the finite-dimensional Lie subalgebra of $\asll_N$ generated by  $ E_s,F_s,H_s$ for $ 1\leq s < N.$ We denote by $\bar{{\mathfrak h}}_N$ the Cartan subalgebra of $\sll_N,$ and by $ \bar{\Lambda}_1,\dots,\bar{\Lambda}_{N-1} $ the set of fundamental weights of $\sll_N.$
Let $\Lambda = \sum_{s=0}^{N-1} a_s \Lambda_s $ $(a_s \in \Cplx)$ be a weight of  $\sll_N.$ We denote by $\bar{\Lambda}$ the finite part of $\Lambda,$ i.e.: $ \bar{\Lambda} = \sum_{s=0}^{N-1} a_s \bar{\Lambda}_s ,$ where $\bar{\Lambda}_0 := 0.$ We denote  by $\bar{Q}_N,$ $P_N^+$ and  $P_N^+(L)$  the root lattice of $\sll_N,$ the cone of dominant integral weights of  $\asll_N$ and  the set of dominant integral weights of level $L$ respectively:
\begin{align*}
\mbox{} & \bar{Q}_N = \oplus_{s=1}^{N-1} \Zint \bar{\tau}_s, \\
\mbox{} &P_N^+ = \{ \Lambda = \sum_{s=0}^{N-1} a_s \Lambda_s \: | \: a_s \in \Zint_{\geq 0} \},\\
\mbox{} &P_N^+(L) = \{ \Lambda \in P_N^+ \: | \: a_0 + a_1 + \cdots + a_{N-1} = L  \}.
\end{align*}

Let $L$ be an integer $\geq 2.$ In what follows we will often consider the pair of algebras $\asll_N$ and $\asll_L.$  In order to distinguish the two, we denote the generators of $\asll_L$ by $e_a,f_a,h_a $ $( 0\leq a < L),$ its Cartan subalgebra by ${\mathfrak h}_L,$ the corresponding fundamental weights by $\omega_0,\omega_1,\dots,\omega_{L-1},$ and set $\vartheta_a = \omega_a - \omega_{a-1}$ for all $a\in \Zint.$ As above, we set  $\omega_a = \omega_{(a\bmod L)},$ and denote by $\bar{\omega}$ the finite part of the weight $\omega \in   {\mathfrak h}_L^*.$

\subsubsection{The wedge product and the Fock space}  \label{sec:wedge} 
\hfill \\ 
\noindent {\bf  Vector modules and their affinizations.} 
Let $\Cplx^N = \oplus_{s=1}^N \Cplx \uf_s $ be the vector module of the finite-dimensional Lie algebra $\sll_N.$ Let $E_{st} \in \End(\Cplx^N)$ be the matrix units in the basis $\{\uf_1,\dots,\uf_N\}.$ The actions of the generators of $\sll_N$ on the vector module are given by  
\begin{align*}
E_s\cdot v & = E_{s,s+1}\cdot v, \\   F_s\cdot v & = E_{s+1,s}\cdot v, \\
H_s\cdot v & = (E_{ss} - E_{s+1,s+1})\cdot v\qquad \left( s=1,\dots,N-1; \; v \in \Cplx^N\right),
\end{align*}
so that the weight of the basis element $\uf_s$ is $\bar{\vep}_s.$ 
Likewise, let $\Cplx^L = \oplus_{a=1}^L \Cplx \vf_s $ be the vector module of the Lie algebra $\sll_L.$ Let $e_{ab} \in \End(\Cplx^L)$ be the matrix units in the basis $\{\vf_1,\dots,\vf_L\}.$ We define the actions of the generators of $\sll_L$ on $\Cplx^L$ by 
\begin{align*}
e_{L-a}\cdot v & = e_{a+1,a}\cdot v,\\  f_{L-a}\cdot v & = e_{a,a+1}\cdot v, \\
h_{L-a}\cdot v & = (e_{a+1,a+1} - e_{aa})\cdot v\qquad \left( a=1,\dots,L-1; \; v \in \Cplx^L\right),
\end{align*}
so that the weight of the basis element $\vf_a$ is $\bar{\vartheta}_{L+1-a}.$ 
The actions of $\sll_L$ and $\sll_N$ are naturally  extended to the tensor product  $\Cplx^L \otimes \Cplx^N,$ so that for $x\otimes y \in  \Cplx^L \otimes \Cplx^N$ we have 
\begin{displaymath}
 a\cdot (x\otimes y) = \begin{cases} (a\cdot x)\otimes y & \text{if $a \in \sll_L,$ }\\
                                   x \otimes  (a\cdot y) & \text{if $a \in \sll_N.$ }\end{cases}
\end{displaymath}
Clearly these two actions are mutually commutative.

Let $z$ be a formal variable, and let $\Vaff = \Cplx[z^{\pm 1}]\otimes\Cplx^L \otimes \Cplx^N$ be the affinization of the $\sll_L \oplus \sll_N$-module  $\Cplx^L \otimes \Cplx^N.$ We define level 0 actions of $\asll_L$ and  $\asll_N$ on $t\otimes x\otimes y \in \Vaff$ by setting $ a\cdot (t\otimes x\otimes y ) = t\otimes a\cdot (x\otimes y)$ if $a \in \sll_L$ or $a \in \sll_N,$ and   
\begin{alignat*}{4}
&  e_0\cdot (t\otimes x\otimes y ) = zt\otimes e_{1L}\cdot (x\otimes y ),& \quad& E_0\cdot (t\otimes x\otimes y ) = zt\otimes E_{N1}\cdot (x\otimes y ) , \\
& f_0\cdot (t\otimes x\otimes y ) = z^{-1}t\otimes e_{L1}\cdot (x\otimes y ),& \quad &  F_0\cdot (t\otimes x\otimes y ) = z^{-1}t\otimes E_{1N}\cdot (x\otimes y ).
\end{alignat*}
These actions of $\asll_L$ and  $\asll_N$ are obviously mutually commutative.

Let us now introduce a notation concerning a basis of $\Vaff.$ With each integer $k$ we associate the unique triple $ \ov{k},\dot{k},\un{k}$ such that $\ov{k} \in \{1,\dots,N\},$ $\dot{k} \in \{1,\dots,L\},$ $\un{k} \in \Zint,$ and   
\begin{displaymath}
k = \ov{k} - N (\dot{k} + L \un{k}).
\end{displaymath}
For each integer $k$ we define $u_k := z^{\un{k}}\otimes \vf_{\dot{k}} \otimes \uf_{\ov{k}}.$ Then the set $\{ u_k \: |\: k \in \Zint \}$ is a basis of $\Vaff.$

\noindent{\bf The wedge product.}
Consider the tensor product $\Vaff^{\otimes n}.$ As a linear space  $\Vaff^{\otimes n}$ is naturally isomorphic to $\Cplx[z_1^{\pm 1},\dots,z_n^{\pm 1}]\otimes (\Cplx^L)^{\otimes n}\otimes (\Cplx^N)^{\otimes n}.$ The algebras   $\asll_L$ and $\asll_N$ act on $\Vaff^{\otimes n}$ diagonally (i.e. by means of the comultiplication). 

Let $T_i$ $(i=1,\dots,n-1)$ be the permutation operator with respect to the factors $i$ and $i+1$ in the tensor product $\Vaff^{\otimes n}.$ We define the {\em wedge product } of  $\Vaff$ as the following quotient linear space:
\begin{equation}
\Vaff^{\wedge n} = \Vaff^{\otimes n} / \sum_{i=1}^{n-1} \Ker ( T_i - 1). \label{eq:wedgeprod}
\end{equation}
Denote by $\wedge$ the quotient map $\Vaff^{\otimes n} \rightarrow \Vaff^{\wedge n}.$ For integer $k_1,\dots,k_n$ define the {\em wedge vector } (or, simply, the {\em wedge}) as 
\begin{displaymath}
u_{k_1}\wedge u_{k_2}\wedge \cdots \wedge u_{k_n} = \wedge( u_{k_1}\otimes u_{k_2}\otimes \cdots \otimes u_{k_n}). 
\end{displaymath}
The factors of the wedge obey the standard fermionic exchange relations, i.e. we have
\begin{displaymath}
u_{k_i}\wedge u_{k_{i+1}} = - u_{k_{i+1}}\wedge u_{k_i} \qquad (i=1,\dots,n-1).
\end{displaymath}
Moreover, the set $\{ u_{k_1}\wedge u_{k_2}\wedge \cdots \wedge u_{k_n} \: | \: k_1 > k_2 > \dots > k_n \}$ is a basis of $\Vaff^{\wedge n}.$ We will call an element of this basis a {\em normally ordered wedge.}

Since the actions of $\asll_L$ and $\asll_N$ on $\Vaff^{\otimes n}$ commute  with $T_i,$ they factor through the quotient map $\wedge$ and give rise to mutually commutative  actions  of $\asll_L$ and $\asll_N$ on the wedge product. 

\noindent{\bf The Fock space.}
Let $\Vaff^{\wedge \frac{\infty}{2}}$ be the semi-infinite wedge product of the linear spaces $\Vaff.$    The {\em Fock space} ${\mathcal F}$ is defined as a subspace of $\Vaff^{\wedge \frac{\infty}{2}}$ spanned by semi-infinite wedges  $u_{k_1}\wedge u_{k_2} \wedge \cdots $ that satisfy the asymptotic condition: $ k_{i+1} = k_{i} - 1$ for all but finite number of $i\in\Nat.$ 
For each $M\in \Zint$ the component  ${\mathcal F}_M$ of {\em charge} $M$ is defined as the linear span of wedges $u_{k_1}\wedge u_{k_2} \wedge \cdots$ that satisfy the asymptotic condition: $ k_i = M-i + 1$ for all but finite number of $i\in\Nat.$ Clearly we have $\F = \oplus_{M \in \Zint} \F_M.$ For each $M\in \Zint$  the vector $|M \rangle = u_{M}\wedge u_{M-1}\wedge u_{M-2}\wedge \cdots \in \F_M$ is called the {\em vacuum of charge }$M.$ It is well-known that $\F$ admits a description as an irreducible Fock module of the Clifford algebra -- hence its name \cite{JM,KR,J}. 
For each integer $k$ and a normally ordered wedge $u_{k_1}\wedge u_{k_2} \wedge \cdots $ define the operators $\psi_k^*$ and $\psi_k$ as      
\begin{eqnarray*}
& \psi^*_k \cdot (u_{k_1}\wedge u_{k_2} \wedge \cdots ) & = u_{k}\wedge u_{k_1}\wedge u_{k_2} \wedge \cdots  , \\  
& \psi_k \cdot (u_{k_1}\wedge u_{k_2} \wedge \cdots ) & = \begin{cases} (-1)^{i-1} u_{k_1}\wedge u_{k_2} \wedge  \cdots \wedge u_{k_{i-1}}\wedge u_{k_{i+1}} \wedge \cdots  & \text{if $ k_i = k,$} \\  
   0  \quad  \text{if $ k_i \neq  k$ for all $i=1,2,\dots.$  }&  \end{cases}
\end{eqnarray*}
These operators satisfy the defining relations of the Clifford algebra:
\bdm
\psi_k^*\psi_l + \psi_l \psi^*_k = \delta_{kl} 1, \quad \psi_k^*\psi_l^* + \psi_l^* \psi^*_k = 0, \quad \psi_k\psi_l + \psi_l \psi_k = 0,
\edm 
and generate the entire Fock space from the vacuum vector $|0\rangle.$

The actions of the Lie algebras $\asll_L$ and $\asll_N$ on $\F$ are defined as in the finite situation -- by means of the  comultiplication. More precisely, let $a$ be an element of $\asll_L$ or $\asll_N.$ We define the action of $a$ on the semi-infinite wedge $u_{k_1}\wedge u_{k_2}\wedge \cdots \in \F$ as   
\begin{equation}
a\cdot (u_{k_1}\wedge u_{k_2}\wedge \cdots \:) = (a\cdot u_{k_1})\wedge u_{k_2}\wedge \cdots \; + u_{k_1}\wedge(a\cdot u_{k_2})\wedge \cdots  \; + \dots, \label{eq:slaction}
\end{equation}
with the additional normalization $a\cdot |0 \rangle  = 0$  if $a$ is an element of the Cartan subalgebra of $\sll_L$ or $\sll_N.$ 
\begin{propos}\mbox{} \label{p:slacts}\\
{\em (i)} The rule {\em (\ref{eq:slaction})} together with the above normalization give rise to  well-defined, mutually commutative actions of the algebras   $\asll_L$ and $\asll_N$ on $\F.$ \\
{\em (ii)} The level of the $\asll_L$ action on $\F$ equals $N,$ and the level of the $\asll_N$ action on $\F$ equals $L.$   
\end{propos}
\begin{proof}For the proof of this proposition it is convenient to express the actions of $\asll_L$ and $\asll_N$ in terms of the Clifford algebra generators (the fermions). For each $s=1,\dots,N;$ $a=1,\dots,L$ and $m \in \Zint$  define $\psi_s^a(m) = \psi_{k}$  where $k= s - N(a+Lm),$ and similarly for the operators $\psi_k^*.$  Let $:\;:$ denote the normal ordering with respect to the vacuum $|0 \rangle,$ and set 
\bdm
J_{st}^{ab}(m) = \sum_{ n \in \Zint } : \psi_s^a(m+n)^* \psi_t^b(n) : \quad (s,t=1,\dots,N; a,b=1,\dots,L; m \in \Zint). 
\edm
The operators $J_{st}^{ab}(m)$ generate a level 1  action of the Lie algebra $\agl_{NL}$ on $\F.$ We have:
\bdm
[J_{st}^{ab}(m),J_{pq}^{cd}(n)] = \delta_{tp}\delta_{bc}J_{sq}^{ad}(m+n) - \delta_{sq}\delta_{ad}J_{pt}^{cb}(m+n) + m \delta_{m+n,0} \delta_{tp} \delta_{sq} \delta_{bc}\delta_{ad} 1. 
\edm 
For all $s,t =1,\dots,N;$ $a,b=1,\dots,L;$ and $m \in \Zint$ define 
\bdm J_{st}(m) = \sum_{a=1}^L J_{st}^{aa}(m) \quad  \text{and} \quad J^{ab}(m) = \sum_{s=1}^N J_{ss}^{ab}(m).
\edm 
The operators $J_{st}(m)$ and $J^{ab}(m)$ generate an action of $\agl_N$ of level $L$ and an action of $\agl_L$ of level $N$ respectively. We have:
\begin{align*}
& [J_{st}(m),J_{pq}(n)] = \delta_{tp}J_{sq}(m+n) - \delta_{sq}J_{pt}(m+n) + L m \delta_{m+n,0} \delta_{tp} \delta_{sq} 1, \\  
&[ J^{ab}(m),J^{cd}(n)] = \delta_{bc}J^{ad}(m+n) - \delta_{ad}J^{cb}(m+n) + N m \delta_{m+n,0} \delta_{bc}\delta_{ad} 1.
\end{align*}
Now, the  generators of the $\asll_N$-action on $\F$ are expressed in terms of the fermions as:  
\begin{alignat*}{6}
& & H_s = J_{ss}(0) - J_{s+1,s+1}(0),  &\; \; E_s = J_{s,s+1}(0),  &\; \; F_s = J_{s+1,s}(0), & \quad  (s=1,\dots,N-1),  \\
& & &\; \;E_0 = J_{N1}(1), &\; \;F_0 = J_{1N}(-1). &     
\end{alignat*}
And  the  generators of the $\asll_L$-action on $\F$ are expressed as: 
\begin{alignat*}{6}
& & h_{L-a} =  J^{a+1,a+1}(0) - J^{aa}(0) ,& \; \;   e_{L-a} = J^{a+1,a}(0), &\; \;f_{L-a} = J^{a,a+1}(0),&  \; \; (a=1,\dots,L-1),  \\
& & & \; \;e_0 = J^{1L}(1), &\; \;f_0 = J^{L1}(-1). & 
\end{alignat*}
The proposition immediately follows from these expressions. \end{proof} 

\noindent For each $m\in \Zint_{\neq 0}$  define an operator $B(m)$ $\in$ $\End(\F)$ by setting for  $u_{k_1}\wedge u_{k_2}\wedge \cdots \in \F$: 
\begin{equation}
B(m)\cdot (u_{k_1}\wedge u_{k_2}\wedge \cdots \:) = (z^m\cdot u_{k_1})\wedge u_{k_2}\wedge \cdots \; + u_{k_1}\wedge(z^m \cdot u_{k_2})\wedge \cdots  \; + \dots \;. \label{eq:heisenbergaction}
\end{equation}
Expressing the $B(m)$ in terms of the Clifford algebra generators one proves the following proposition:
\begin{propos}\mbox{}\\
{\em (i)} The action {\em (\ref{eq:heisenbergaction})} is well-defined. The operators $\{B(m) \: | \: m\in \Zint_{\neq 0}\}$ satisfy the defining relations of the Heisenberg algebra $H$ with the central charge $LN,$ i.e.:  
\begin{displaymath}
[B(m),B(n)] = LN m \delta_{m+n,0} \cdot 1.
\end{displaymath}
{\em (ii)}
The above action of $H$ centralizes the actions of $\asll_L$ and $\asll_N$ on $\F.$
\end{propos}
\noindent Thus we have an  action of the Lie algebra $H\oplus \asll_N \oplus \asll_L$ on $\F.$ This action clearly leaves invariant the charge component $\F_M$ for each $M\in \Zint.$ 

\noindent{\bf Decomposition of the Fock space.} Let $S_- = \Cplx[B(-1),B(-2),B(-3),\dots \;]$ be the  Fock module of the Heisenberg algebra $H.$ The standard action of $H$ on $b \in S_{-}$ is defined as 
\begin{displaymath}    
B(m)\cdot b = B(m) b \quad(m < 0), \qquad B(m)\cdot b = [B(m),b] \quad (m > 0).
\end{displaymath}

With every pair $M,\Lambda$ where $M\in \Zint,$ and $\Lambda = \sum_{s=0}^{N-1} a_s \Lambda_s \in P_N^+(L)$ such that $\bar{\Lambda} \equiv \bar{\Lambda}_M \bmod \bar{Q}_N$ we associate the $\asll_L$ weight $\omega_{M,\Lambda}\in P_L^+(N) $ as follows:
\begin{displaymath}
\omega_{M,\Lambda}  = \omega_{l_1} + \omega_{l_2} + \dots + \omega_{l_N},
\end{displaymath}
where $l_1,l_2,\dots,l_N$ are integers defined by 
\begin{equation} \begin{aligned} 
  \mbox{}    & a_s = l_s - l_{s+1} \quad ( 1\leq s \leq N-1), \\ 
    \mbox{}   &    M = l_1 + l_2 + \cdots + l_N.
\end{aligned} \label{eq:lsss}\end{equation}
Note that $\omega_{M,\Lambda} \neq \omega_{M,\Lambda^{\prime}}$ if $\Lambda \neq \Lambda^{\prime}.$

Let $\Lambda \in P_N^+,$ and let $V(\Lambda)$ be the irreducible (integrable) highest weight module of the affine Lie algebra $\asll_N$ with the highest weight $\Lambda.$ Likewise, for  $\omega \in P_L^+,$ let $V(\omega)$ be the irreducible (integrable) highest weight module of the affine Lie algebra $\asll_L$ with the highest weight $\omega.$

The following theorem, due to I.Frenkel \cite{F}, gives the decomposition of the Fock space into irreducible  modules of the Lie algebra $H\oplus\asll_N\oplus\asll_L:$ 
\begin{theor} \label{t:Fockdec}
For each $M\in \Zint$ one has the following isomorphism of $H\oplus\asll_N\oplus\asll_L$-modules:
\begin{equation}
\F_M \simeq \bigoplus_{\{ \Lambda \in P_N^+(L) \:|\: \bar{\Lambda} \equiv \bar{\Lambda}_M\bmod \bar{Q}_N \}}  S_-\otimes V(\Lambda) \otimes V(\omega_{M,\Lambda}). \label{eq:Fockdec}
\end{equation}
\end{theor}
\noindent Note that each weight $\Lambda \in P_N(L)^+$ appears in the above decomposition when $M$ runs through the set $\{0,1,\dots,N-1\}.$ To describe the isomorphism (\ref{eq:Fockdec}) explicitly it is enough to point out the highest weight vectors of $H\oplus\asll_N\oplus\asll_L$ in $\F_M.$ With $l_1,\dots,l_N$ defined as in (\ref{eq:lsss}), these are given as follows. For each integer $l$ and each $s=1,\dots,N$ define  
\bdm
\Xi_s(l) = \begin{cases} \psi_{s + N(l-1)}^*\psi_{s + N(l-2)}^*\cdots \psi_{s}^* & \text{ if $ l > 0,$} \\
     1 & \text{ if $ l = 0,$} \\
\psi_{s - N(-l)}\psi_{s - N(-l-1)}\cdots \psi_{s - N} & \text{ if $ l < 0.$}
\end{cases} \edm
Then under the isomorphism (\ref{eq:Fockdec}) the vector  
\bdm 
\Xi_1(l_1)\Xi_2(l_2)\cdots \Xi_N(l_N) \cdot | 0 \rangle  \quad \in \F_M
\edm
corresponds to 
\bdm
1 \otimes | \Lambda \rangle \otimes |\: \omega_{M,\Lambda} \rangle
\edm 
where $| \Lambda \rangle $ is the highest weight vector in $V(\Lambda),$ $| \:\omega_{M,\Lambda}\rangle  $ is the highest weight vector in $V(\omega_{M,\Lambda})$ and $1$ is the highest weight vector in $S_-.$

\subsection{Degenerate Double Affine Hecke Algebra} \label{sec1:daha}
Here we collect the necessary facts about the Degenerate Double Affine Hecke Algebra. Exposition in this section closely follows the work \cite{AST}.
\subsubsection{Affine Root System}
Let $\bar{ \tf} = \oplus_{i=1}^n \Cplx \ep_i^{\vee}$ be the Cartan subalgebra of the Lie algebra $\gl_n(\Cplx)$  and let $\tf = \bar{\tf}\oplus\Cplx c \oplus \Cplx d$ be the  Cartan subalgebra of the affine Lie algebra $\agl_n(\Cplx).$ The non-degenerate  bilinear symmetric form $(\:,\:)$ on $\tf$ is defined by setting $(\ep_i^{\vee},\ep_j^{\vee})=\delta_{ij},$ $(c,d) = 1,$ $ (\ep_i^{\vee},c)=(\ep_i^{\vee},d)=(c,c)=(d,d) = 0.$ Let $\bar{ \tf}^* = \oplus_{i=1}^n \Cplx \ep_i$ be the dual space of $\bar{ \tf}$ and $\tf^* = \bar{\tf}^*\oplus\Cplx c^* \oplus \Cplx \delta $ be the dual space of $\tf,$ where $\ep_i, \delta $ and $c^*$ are the dual vectors of  $\ep_i^{\vee}, d $ and $c$ respectively. We identify $\tf^*$ with $\tf$ via the correspondences $\ep_i \mapsto \ep_i^{\vee}, \delta \mapsto c $ and $c^* \mapsto d.$ For a vector $\zeta \in \tf^*$ we denote by $\zeta^{\vee}$ the vector of $\tf$ obtained through this identification.   

Let $\bar{R},\bar{R}_+$ and $\bar{\Pi}$ be, respectively, the  root system, the set of positive roots and the set of simple roots of type $A_{n-1}:$ 
\begin{eqnarray*}
& & \bar{R} =\{ \alpha_{ij} = \ep_i - \ep_j \: | \: i\neq j \},\\
& & \bar{R}_+ =\{ \alpha_{ij} \:|\: i < j \},\\
& & \bar{\Pi} = \{ \alpha_1,\dots,\alpha_{n-1} \}\qquad ( \alpha_i := \alpha_{ii+1}).
\end{eqnarray*}
The affine root system $R,$ the set of positive roots $R_+$ and the set of simple roots $\Pi$ of type $A_{n-1}^{(1)}$ are defined by 
\begin{eqnarray*}
& & R =\{ \alpha + k\delta \: | \: \alpha \in  \bar{R}, k \in \Zint  \},\\
& & {R}_+ =\{ \alpha + k\delta \: | \: \alpha \in  \bar{R}_+ , k \geq 0 \}\sqcup \{ -\alpha + k\delta \: | \: \alpha \in  \bar{R}_+ , k > 0 \},\\
& & {\Pi} = \{ \alpha_0 := \delta - (\ep_1 - \ep_n)\} \sqcup \bar{\Pi}. 
\end{eqnarray*}

\subsubsection{Affine Weyl Group}

Let $\bar{W}$ be the Weyl group of the root system $\bar{R},$ it is isomorphic to the symmetric group ${\mathfrak S}_n.$ Let $\bar{P} = \oplus_{i=1}^n \Zint \ep_i$ be the weight lattice of  $\gl_n(\Cplx).$ The {\em affine Weyl group} $W$ is defined as the semidirect product  
\setlength{\unitlength}{0.1mm}\begin{displaymath}
 W = \bar{W} \semidir \;\bar{P},
\end{displaymath}
with the relations $w\cdot t_{\eta} \cdot w^{-1} = t_{w(\eta)},$ where $w$ and $t_{\eta}$ are the elements of $W$ that correspond to $w\in \bar{W}$ and $\eta \in \bar{P}$ respectively.  

Let $s_{\alpha} \in \bar{W}$ be the reflection that corresponds to the root $\alpha \in \bar{R}.$ The action of $W$ on an element $\xi \in \tf $ is given by 
\begin{eqnarray*}
& & s_{\alpha}(\xi) = \xi - \alpha(\xi) \alpha^{\vee}\qquad \left( \alpha \in \bar{R} \right),\\
& & t_{\eta}(\xi) = \xi + \delta(\xi)\eta^{\vee} - \left( \eta(\xi) + \frac{1}{2}(\eta,\eta)\delta(\xi)\right) c \qquad \left(\eta \in \bar{P}\right).
\end{eqnarray*}
This action leaves invariant the linear subspace $\tf^{\prime} = \bar{\tf} \oplus \Cplx c \subset \tf.$ 
The dual action of $W$ on a vector $\zeta \in \tf^*$ is given by 
\begin{eqnarray*}
& & s_{\alpha}(\zeta) = \zeta - (\alpha , \zeta) \alpha \qquad \left( \alpha \in \bar{R} \right),\\
& & t_{\eta}(\zeta) = \zeta  + (\delta, \zeta ) \eta - \left( (\eta, \zeta) + \frac{1}{2}(\eta,\eta)(\delta,\zeta)\right) \delta \qquad \left(\eta \in \bar{P}\right).
\end{eqnarray*}
This action leaves invariant the subspace $ \Cplx \delta \subset \tf^*.$ Therefore we have an action of $W$ on the linear  space $ (\tf^{\prime})^* = \tf^*/ \Cplx \delta $ which we identify with $\bar{\tf}^*\oplus \Cplx c^*:$ 
\begin{equation} 
\begin{aligned} 
 \mbox{} &   s_{\alpha}(\zeta)  = \zeta - (\alpha , \zeta) \alpha \qquad \left( \alpha \in \bar{R}, \zeta \in (\tf^{\prime})^* \right),\\
 \mbox{} &  t_{\eta}(\zeta)  = \zeta  +  \zeta (c) \eta   \qquad \left(\eta \in \bar{P}, \zeta \in (\tf^{\prime})^* \right). 
\end{aligned}  \label{eq:Wact}
\end{equation}
For an affine root $\alpha = \bar{\alpha} + k\delta $ $(\bar{\alpha}\in \bar{R}, k \in \Zint),$ define the corresponding affine reflection by $s_{\alpha} = t_{-k\bar{\alpha}}\cdot s_{\bar{\alpha}}.$ Set $s_i = s_{\alpha_i}$ for $i=0,\dots,n-1.$ We will identify the set $\{ 0,\dots,n-1\}$ with the abelian group $\Zint/n\Zint.$ Let $\pi = t_{\ep_1}\cdot s_1\cdots s_{n-1}.$ 
\begin{propos}
The group $W$ is isomorphic to the group defined by  generators \\$\pi,\pi^{-1},$ $  s_i \; (i\in \Zint/n\Zint)$ and relations 
\begin{eqnarray*}
 & & s_i^2 = 1, \quad \pi\cdot \pi^{-1} = 1,\quad  \pi\cdot s_i = s_{i+1} \cdot \pi,\\
  & & s_i\cdot s_j = s_j \cdot s_i \; (i-j \not\equiv \pm 1\bmod n),  \\ 
  & & s_i\cdot s_{i+1}\cdot s_i = s_{i+1}\cdot s_i \cdot s_{i+1}.
\end{eqnarray*}

\end{propos}
For $w\in W,$ let $S(w) = R_+\cap w^{-1}\left(R_-\right),$ where $R_- = R\setminus R_+$ is the set of negative roots. The length  $l(w)$ of $w$ is defined as the number of elements in the set $S(w).$ For $w\in W,$ an expression $w = \pi^k\cdot s_{i_1}\cdots s_{i_l}$ is called a reduced expression if $l = l(w).$ Let $W^a$ be the subgroup of $W$ generated by $s_0,\dots s_{n-1}.$ For $w,w^{\prime}\in W^a$ write $w \preceq w^{\prime}$ if $w$ can be obtained as a subexpression of a reduced expression for $w^{\prime}.$ The partial ordering $\preceq$ is extended on $W$ by $\pi^{k}w \preceq \pi^{k^{\prime}}w^{\prime}$ $\Leftrightarrow$ $k=k^{\prime}$ and $w \preceq w^{\prime}$ $(w, w^{\prime} \in W^a).$

\subsubsection{Degenerate Double Affine Hecke Algebra}

Let $\Cplx[W]$ denote the group algebra of $W$ and let $S[\tf^{\prime}]$ denote the symmetric algebra of $\tf^{\prime}.$ We have $\Cplx[W]=\Cplx[\bar{P}]\otimes \Cplx[\bar{W}] $ and $S[\tf^{\prime}] = S[\bar{\tf}]\otimes\Cplx[c].$ 

\begin{defin}[\cite{C4}]
The {\em degenerate double affine Hecke algebra } $\daha$ is the unital associative $\Cplx$-algebra defined by the following properties:  \\
{\em (i)} As a $\Cplx$-vector space 
\begin{gather*}
\daha = \Cplx[W]\otimes S[\tf^{\prime}],
\end{gather*}
{\em (ii)} The inclusions $\Cplx[W]\hookrightarrow \daha$ and $S[\tf^{\prime}]\hookrightarrow \daha$ are algebra homomorphisms.\\
{\em (iii)} The following relations hold in $\daha:$
\begin{align}
\mbox{} &  s_i\cdot \xi - s_i(\xi)\cdot s_i = - \alpha_i(\xi) \qquad (i=0,\dots,n-1, \; \xi \in \tf^{\prime}),\\
\mbox{} & \pi\cdot \xi = \pi(\xi)\cdot \pi \qquad  (\xi \in \tf^{\prime}).
\end{align}
\end{defin}
\begin{defin}
The {\em degenerate affine Hecke algebra }$\bar{\daha}$ is the following subalgebra of $\daha:$
\begin{displaymath}
\bar{\daha} = \langle w\in \bar{W}, \xi \in \bar{\tf}\rangle  \cong \Cplx[\bar{W}]\otimes S[\bar{\tf}].
\end{displaymath}
\end{defin}
Let us identify the group algebra $\Cplx[\bar{P}]$ with the ring of Laurent polynomials in variables $e^{\ep_1},\dots,e^{\ep_n}.$ The following proposition gives an alternative description of the algebra $\daha.$
\begin{propos}[\cite{AST}]
The algebra $\daha$ is the unital associative $\Cplx$-algebra such that as a vector space
\begin{displaymath}
 \daha = \Cplx[\bar{P}] \otimes \bar{\daha} \otimes \Cplx[c],
\end{displaymath}
the element $c$ is central,   
\begin{eqnarray*}
& & w\cdot f \cdot w^{-1} = w(f) \quad \left( w\in \bar{W}, f \in \Cplx[\bar{P}]\right),\\
& & [ \xi, f] = c \cdot \partial_{\xi}(f)  + \sum_{\alpha \in \bar{R}_+ } \alpha(\xi) \frac{(1-s_{\alpha})(f)}{1-e^{-\alpha}} \cdot s_{\alpha} \quad \left(\xi \in \bar{\tf}, f \in \Cplx[\bar{P}]\right),
\end{eqnarray*}
where $\partial_{\xi}(e^{\eta}) := \eta(\xi)e^{\eta}$ $(\xi \in \bar{\tf}, \eta \in \bar{P})$ and the inclusions  $\Cplx[\bar{P}] \hookrightarrow \daha,$ $ \bar{\daha} \hookrightarrow \daha $ are algebra homomorphisms. 
\end{propos}

\subsubsection{A representation of the Degenerate Double Affine Hecke Algebra} \label{sec:DAHArep}

Let $z_1,\dots,z_n$ be a set of formal variables, and consider the linear space  
\begin{equation}
\Ve = \Cplx[z_1^{\pm 1},\dots,z_n^{\pm 1}] \otimes \left(\Cplx^L\right)^{\otimes n} \cong  \left( \Cplx[z^{\pm 1}] \otimes \Cplx^L\right)^{\otimes n}. 
\end{equation}
Let $K_{ij}$ be the exchange operator for variables $z_i,z_j$ in $\Cplx[z_1^{\pm 1},\dots,z_n^{\pm 1}]$ and let $P_{ij}$ be the exchange operator of factors $i$ and $j$ in the tensor product  $\left(\Cplx^L\right)^{\otimes n}.$ For $  e \in \End\Cplx^L,$ we set 
\bdm
e_i :=  1\otimes \cdots \otimes \underset{i}{e} \otimes \cdots \otimes 1\; \in \End\left(\Cplx^L\right)^{\otimes n}. 
\edm
If $B$ is an element of $\End \Cplx[z_1^{\pm 1},\dots,z_n^{\pm 1}]$ or $\End\left(\Cplx^L\right)^{\otimes n}$ we will keep the same symbol $B$ to denote the natural extension of $B$ on $\Ve.$ 

Let $\kappa$ be a complex number, and let $\nu = \sum_{a=1}^L \nu(a) E_{aa} \in \End \Cplx^L$ be a diagonal matrix with arbitrary complex entries $\nu(1),\dots,\nu(L).$ For each $i=1,\dots,n$ introduce the {\em matrix Dunkl operator } \cite{C3,AST}:
\begin{multline}
d_i = \kappa z_i\frac{\partial}{\partial z_i} + \nu_i + \frac{n}{2L}-\frac{1}{2}+\sum_{j > i} \frac{z_j}{z_j-z_i}(K_{ij}-1)P_{ij} + r_{ij} - \label{eq:dunkl}\\ 
- \sum_{j < i} \frac{z_i}{z_i-z_j}(K_{ij}-1)P_{ij} + r_{ji}, 
\end{multline}
where $r_{ij}$ is the constant {\em classical $r$-matrix} associated with $\sll_L$: 
\begin{equation}
r_{ij} = \frac{1}{2}\sum_{a=1}^L (e_{aa})_i(e_{aa})_j + \sum_{1\leq a < b \leq L} (e_{ab})_i(e_{ba})_j.
\end{equation}

Introduce the map $\pi_{\kappa,\nu} : \daha \rightarrow \End\Ve$ by 
\begin{equation}\left.\begin{aligned}[8]
&t_{\ep_i} & \mapsto z_i^{-1} & \quad & \ep_i^{\vee} & \mapsto -d_i &  &(i=1,\dots,n), \\
&c & \mapsto \kappa 1, &   \quad & s_i & \mapsto K_{ii+1}P_{ii+1} & \quad  &(i=1,\dots,n-1).
\end{aligned}\right.\end{equation}

\begin{propos}[\cite{C3} Theorem 2.3]
The map $\pi_{\kappa,\nu}$ defines a left $\daha$-module structure on $\Ve.$
\end{propos}

The linear space $\Ve$ is also an $\asll_L$-module with the action of $\asll_L$ defined as in Section \ref{sec:wedge}. Let $\chi(1),\dots,\chi(L)$ be complex numbers, and let $\UU^{\prime}(\bef^{\chi})$ be the associative, non-unital subalgebra of $\UU(\asll_L)$ generated by the elements
\begin{align*}
\mbox{} &  h_a^{\chi} :=  h_a - \left( \chi(L+1-a) - \chi(L-a) \right)1 \quad ( a=1,\dots,L-1), \\
\mbox{} &  f_0,f_1,\dots,f_{L-1}.
\end{align*}

\begin{propos} \label{p:dahainv}
Let $\kappa = L,$ and let $\nu(a) = a + \frac{\chi(a)}{2}$ $\;(a=1,\dots,L).$ Then the double degenerate affine Hecke algebra action $\pi_{\kappa,\nu}(\daha)$ leaves invariant the subspace $\UU^{\prime}(\bef^{\chi})\Ve.$
\end{propos}
\noindent \begin{proof}
The actions of $\Cplx[\bar{P}]\hookrightarrow \daha$ and $\bar{W}\hookrightarrow \daha$ obviously commute with the action of $\UU(\asll_L).$ The action of $S[\bar{\tf}]=\langle \ep_1^{\vee},\dots,\ep_n^{\vee}\rangle \hookrightarrow \daha$ commutes with the action of the Cartan subalgebra $\bar{{\mathfrak h}}_L.$ Finally, for each $i=1,\dots,n$ we have 
\begin{align*}
\mbox{}&[d_i,f_{a}] = \frac{1}{2} h_{a}^{\chi} (e_{L-a,L-a+1})_i + \frac{1}{2} f_{a} \left( (e_{L-a,L-a})_i - (e_{L-a+1,L-a+1})_i\right) \quad (a=1,\dots,L-1), \\
\mbox{}& [d_i,f_0 ] = -\frac{1}{2}(h_1^{\chi}+\cdots + h_{L-1}^{\chi}) z_i^{-1}(e_{L1})_i + \frac{1}{2}f_0 \left((e_{LL})_i - (e_{11})_i\right).
\end{align*} 
\end{proof}

\subsection{The Yangian of the Lie algebra $\gl_N$} \label{se:yangian}
In this section we recall  basic definitions and few known facts about the Yangian of $\gl_N.$ Our basic references are \cite{NT1,NT2}. Some  background information on representations of the Yangian is summarized in Section \ref{sec:yangdecofvac}.
\subsubsection{The Yangian of  $\gl_N$} 
The {\em Yangian} of general Lie algebra $\gl_N$ is the associative unital algebra $\YY(\gl_N)$ over $\Cplx$ with the generators $T^{(m)}_{st}$ where $m=1,2,\dots$ and $s,t=1,\dots,N.$ Defining relations of  $\YY(\gl_N)$ can be written in terms of the generating series 
\bdm
T_{st}(u) = \delta_{st} + T_{st}^{(1)}u^{-1} + T_{st}^{(2)}u^{-2}+ \cdots
\edm
in a formal parameter $u$ as follows: for all indices $p,q,s,t=1,\dots,N$ we have
\bdm
(u-v)[T_{pq}(u),T_{st}(v)]= T_{sq}(v)T_{pt}(u)- T_{sq}(u)T_{pt}(v).
\edm
Here $v$ is another formal parameter. 

Let $E_{st} \in \End\Cplx^N$ be the standard matrix units. Combine all the series $T_{st}(u)$ into the single element 
\bdm
 T(u) = \sum_{s,t=1}^N E_{st}\otimes T_{st}(u)
\edm
of the algebra $ \End\Cplx^N \otimes \YY(\gl_N)[[u^{-1}]].$ The element $T(u)$ is invertible, let
\bdm
T(u)^{-1} = \widetilde{T}(u) =  \sum_{s,t=1}^N E_{st}\otimes \widetilde{T}_{st}(u)
\edm
\begin{propos}[\cite{NT1,NT2}]
The map $\sigma : {T}_{st}(u)\mapsto \widetilde{T}_{st}(-u)$ defines an involution of $\YY(\gl_N).$
\end{propos}
\begin{propos}[\cite{MNO}]
The coefficients at $u^{-1},u^{-2},\dots $ of the series 
\bdm
\Delta(u) = \sum_{\gamma \in {\mathfrak S}_N} {\mathrm {sign}}(\gamma) T_{1\gamma(1)}(u)T_{2\gamma(2)}(u-1)\cdots T_{N\gamma(N)}(u-N+1)
\edm
are free generators of the centre of the algebra $\YY(\gl_N).$
\end{propos}

\subsubsection{The Yangian of $\sll_N$} \label{sec:Ysl}
Let $\{ t_{\alpha}\}$ be a basis of $\sll_N$ in the vector representation, such that $ {\mathrm {tr}}_{\Cplx^N}(t_{\alpha}t_{\beta}) = \delta_{\alpha\beta}.$ The corresponding structure constants $f_{\alpha\beta\gamma}$ are totally antisymmetric and satsify the normalization condition
$$ f_{\alpha\beta\gamma}f_{\lambda\beta\gamma} = -2N\delta_{\alpha\lambda}. $$
The {\em Yangian } of the Lie algebra $\sll_N$ \cite{DrinfeldYang1} is the unital associative algebra $\YY(\sll_N)$ over $\Cplx$ with the generators $ X_{\alpha},  Y_{\alpha}$ where $\alpha =1,\dots,\dim\sll_N.$ The defining relations of $\YY(\sll_N)$ are  
\begin{gather*} 
[ \: X_{\alpha} \: , \: X_{\beta} \: ] \: = \: f_{\alpha\beta\gamma} \,  X_{\gamma} \tag*{(i)}  \\
[ \: X_{\alpha} \: , \: Y_{\beta} \: ] \: = \: f_{\alpha\beta\gamma}\,  Y_{\gamma} \tag*{(ii)}  \\
[\:   Y_{\alpha} \: ,\: [\: Y_{\beta} \: , \:   X_{\gamma} \: ] \: ] \: + \: \text{ cyclic in $\alpha,\beta,\gamma$} \: = \: A_{\alpha\beta\gamma}^{\lambda\mu\nu} \: \{ \: X_{\lambda}\:,\: X_{\mu}\: , \: X_{\nu} \: \} \tag*{(iii)}  \\ 
 [\: [\:  X_{\alpha} \: ,Y_{\beta}\:]\: ,\: [\:  X_{\gamma} \: , \: Y_{\delta}  \: ] \: ] \: + \: [\: [\:  Y_{\gamma} \: ,Y_{\delta}\:]\: ,\: [\:  X_{\alpha} \: , \: Y_{\beta}  \: ] \: ] \tag*{(iv)}\\
 \qquad \qquad  = \left( A_{\alpha\beta\rho}^{\lambda\mu\nu} f_{ \gamma\delta\rho } + A_{\gamma\delta\rho}^{\lambda\mu\nu} f_{\alpha\beta\rho}\right) \: \{ \: X_{\lambda}\:,\: X_{\mu}\: , \: Y_{\nu} \: \}
\end{gather*} 
where $$ 
 A_{\alpha\beta\gamma}^{\lambda\mu\nu} = \frac{1}{4} f_{\alpha\lambda\rho } f_{\beta\mu\sigma}  f_{\gamma\nu\tau}  f_{\rho\sigma\tau}
 \quad \text{and} \quad  \{ \: x_1 \: , \: x_2 \: , \: x_3 \: \} \: = \: \frac{1}{6} \sum_{ \pi \in S_3 } \: x_{\pi(1)} \:  x_{\pi(2)} \: x_{\pi(3)} $$ 
Let $(t_{\alpha})_{st}$ be the matrix of the $\sll_N$-generator $t_{\alpha}$ in the basis $\{\uf_s\}$ of $\Cplx^N:$  $t_{\alpha} \cdot \uf_t = \sum  \uf_s (t_{\alpha})_{st}.$ The algebra $\YY(\gl_N)$ contains $\YY(\sll_N)$ as a subalgebra. The assignments 
\begin{equation} X_{\alpha} \mapsto Q_{\alpha}^{(1)} \overset{{\mathrm {def}}}{=} \sum_{s,t = 1}^N  (t_{\alpha})_{st} T^{(1)}_{ts}, \quad  Y_{\alpha} \mapsto Q_{\alpha}^{(2)} \overset{{\mathrm {def}}}{=} -\sum_{s,t = 1}^N  (t_{\alpha})_{st}\left( T_{ts}^{(2)} - \frac{1}{2}\sum_r T^{(1)}_{tr}T^{(1)}_{rs}\right) \label{eq:embed}
\end{equation}
define the embedding. As algebra $\YY(\gl_N)$ is isomorphic to the tensor product of its centre and $\YY(\sll_N).$

\subsubsection{Yangian action on the wedge product} \label{sec:yangonwedge}
Recall the definitions of the tensor product $\Vaff^{\otimes n}$ and the wedge product $\Vaff^{\wedge n}$ given in Section \ref{sec:wedge}. Recall as well the definition of the $\daha$-module $\Ve$ given in Section \ref{sec:DAHArep}. There is a natural isomorphism of linear spaces
\bdm
\Vaff^{\otimes n} \cong  \Ve\otimes (\Cplx^N)^{\otimes n}.
\edm
In what follows we will identify these spaces by this isomorphism.

The action $\pi_{\kappa,\nu}$ of the degenerate double affine Hecke algebra is naturally extended from $\Ve$ to $\Vaff^{\otimes n}:$ an element $a$ of $\daha$ acts as $\pi_{\kappa,\nu}(a)\otimes {\mathrm {id}}.$ Likewise action of the matrix unit $(E_{st})_i\in \End (\Cplx^N)^{\otimes n}$ is extended to $\Vaff^{\otimes n}$ as ${\mathrm {id}}\otimes (E_{st})_i.$ For all $s,t=1,\dots,N$ we now define the following elements of $\End\Vaff^{\otimes n}[[u^{-1}]]:$
\begin{equation}\begin{aligned}
\mbox{} & L_{st}^i(u) = \delta_{st} + \frac{{\mathrm {id}}\otimes (E_{ts})_i}{u - d_i\otimes {\mathrm {id}}},\\
\mbox{} & \widehat{T}_{st}(u) = \sum_{s_1,s_2,\dots,s_{n-1}}  L_{s,s_1}^1(u)L_{s_1,s_2}^2(u)\cdots L_{s_{n-1},t}^n(u).
\end{aligned}\label{eq:yangdef1}\end{equation}
Here the denominators are to be expanded as  series in $u^{-1}.$

Define the map $\rho_{\kappa,\nu}$ from the system of generators of $\YY(\gl_N)$ into $\End\Vaff^{\wedge n}$ as follows: for each  $ v \in \Vaff^{\otimes n},$ set 
\begin{equation}
\rho_{\kappa,\nu}(T_{st}(u))\cdot \wedge(v) = \wedge \left(\widehat{T}_{st}(u) \cdot v \right).\label{eq:yangdef2}
\end{equation}
The content of the following proposition constitutes a part of the Drinfeld duality between modules of the degenerate affine Hecke algebra $\bar{\daha}$ and modules of the algebra $\YY(\gl_N)$ \cite{DrinfeldDual}.  
\begin{propos}
The map $\rho_{\kappa,\nu}$ defines a $\YY(\gl_N)$-module structure on  $\Vaff^{\wedge n}.$
\end{propos}

The actions of the Lie algebras $\asll_N$ and $\asll_L$ on $\Vaff^{\wedge n}$ were defined in Section \ref{sec:wedge}.

\begin{propos} \label{p:yinv}
Let $\kappa = L,$ and let $\nu(a) = a + \frac{\chi(a)}{2}$ $\;(a=1,\dots,L).$ Then the Yangian action $\rho_{\kappa,\nu}(\YY(\gl_N))$ leaves invariant the subspace $\UU^{\prime}(\bef^{\chi})\Vaff^{\wedge n}.$
\end{propos}
\noindent \begin{proof}
This is a straightforward corollary to  Proposition \ref{p:dahainv}.
\end{proof}

\section{Yangian actions on irreducible integrable modules of $\agl_N$}

\subsection{Yangian action on the Fock space} \label{se:YonF}
In this section we define a Yangian action on each charge component $\F_M$ of the Fock space $\F.$ This action is, informally speaking, a limit of the action on the wedge product $\Vaff^{\wedge n}$ when $n$ is sent to infinity. The formal definition is made possible by stability (cf. Proposition \ref{p:stab}) of the (renormalized) Yangian action on  $\Vaff^{\wedge n}$ when $n$ is incremented by multiples of $NL.$ The procedure we follow here to define a $\YY(\gl_N)$-action on the Fock space is practically the same as the one used in \cite{TU,STU}. The only, inessential, difference is that now we are dealing with a more general representation of $\daha$ than the polynomial representation that underlies the construction of \cite{TU,STU}.

\subsubsection{Intertwining relations and stability of the Yangian action on the wedge product} \label{sec:inter-and-stab}
Let $\F_M$ be the component of charge $M\in \Zint$ of the Fock space $\F$ (cf. Section \ref{sec:wedge}). The linear space $\F_M$ has a basis formed by normally ordered semi-infinite wedges
\begin{equation}
u_{(k_i)}:= u_{k_1}\wedge u_{k_2}\wedge u_{k_3}\wedge \cdots \qquad (k_1 > k_2 > k_3 > \dots\;),\label{eq:siwedge}
\end{equation}
such that the sequence of {\em momenta} $(k_i)_{i\in \Nat}$ satisfies the asymptotic condition: $k_i = M-i+1$ for all but finite number of $i\in \Nat.$ 
Let $|M\rangle $ be the vacuum vector of $\F_M,$ let $(o_i)_{i\in \Nat}$ be the sequence of momenta labeling $|M\rangle:$ $ o_i = M-i+1.$

The linear space $\F_M$ is $\Zint_{\geq 0}$-graded: for any semi-infinite wedge (\ref{eq:siwedge}) the degree is defined as 
\begin{equation}
 \deg u_{(k_i)} = \sum_{i \in \Nat} \un{o_i} - \un{k_i}.
\end{equation}
Let $\F_M^d$ be the homogeneous component of $\F_M$ of degree $d.$ We will define a Yangian action on each of the spaces $\F_M^d$ $(d=0,1,2,\dots\:).$

Let $s \in \{0,1,\dots,NL - 1\}$ be such that $M\equiv s \bmod NL.$ Let $l$ be a non-negative integer, and define a linear subspace of the wedge product (\ref{eq:wedgeprod}) with $s+lNL$ factors as  
\begin{equation}
V_{s+lNL}  = \bigoplus_{ \un{k_{s+lNL}} \leq \un{o_{s+lNL}} } \Cplx u_{k_1}\wedge u_{k_2} \wedge \cdots \wedge u_{k_{s+lNL}}. \label{eq:Vs}
\end{equation}
If $s=l=0$ we set $V_{s+lNL} =\Cplx.$ Above, and in what follows we use the notation $u_{k_1}\wedge u_{k_2} \wedge \cdots \wedge u_{k_{s+lNL}}$ to mean exclusively  a normally ordered wedge.

The vector space (\ref{eq:Vs}) has a grading similar to that one of the vector space $\F_M.$ Now the the degree is defined as 
\begin{equation}
\deg u_{k_1}\wedge u_{k_2} \wedge \cdots \wedge u_{k_{s+lNL}} = \sum_{i =1}^{s+lNL} \un{o_i} - \un{k_i}.
\end{equation}
Notice that this degree is a non-negative integer since $k_1>k_2>\cdots>k_{s+lNL}$ and $\un{k_{s+lNL}} \leq \un{o_{s+lNL}}$ imply  $\un{k_{i}} \leq \un{o_{i}}$ for all $i=1,\dots,s+lNL.$ Let $V_{s+lNL}^d$ $(d=0,1,2,\dots\:)$ denote the homogeneous component of $V_{s+lNL}$ of degree $d.$

Recall that we have the Yangian action $\rho_{\kappa,\nu}$ (cf. Section \ref{sec:yangonwedge}) on the wedge product $\Vaff^{\wedge n}.$ In this section we will denote this action by $\rho_{\kappa,\nu}^{(n)}$ to indicate explicitly the number of factors in the wedge product. For $n=0$ we define  $\rho_{\kappa,\nu}^{(0)}(T_{pq}(u)) = \delta_{pq}.$ The following proposition is a straightforward consequence of the definition of the matrix Dunkl operators (cf. \ref{eq:dunkl}). 

\begin{propos}
For each $d=0,1,2,\dots$ the homogeneous component $V_{s+lNL}^d$ is invariant with respect to the Yangian action $\rho_{\kappa,\nu}^{(s+lNL)}.$  
\end{propos}

\begin{defin}
For each $d=0,1,2\dots$ define a linear map $\vro_l^d : V_{s+lNL}^d \rightarrow \F_M$ by 
\bdm
\vro_l^d(w) = w \wedge |M-s-lNL\rangle \qquad (w \in V_{s+lNL}^d ).
\edm
\end{defin}
\noindent Clearly we have $\deg\vro_l^d(w) = \deg w$ and hence  
\bdm
\vro_l^d : V_{s+lNL}^d \rightarrow \F_M^d \qquad (d=0,1,2,\dots ).
\edm
\begin{propos} \label{p:VFockisom}
Suppose $l \geq d.$ Then $\vro_l^d$ is an isomorphism of vector spaces.
\end{propos}
\noindent \begin{proof} Suppose $w = u_{k_1}\wedge u_{k_2}\wedge \cdots$ is a normally ordered wedge from $\F_M^d.$
Then $k_i = o_i$ for all $i > s+dNL.$ For if otherwise, the degree of $w$ must be greater or equal to $d+1.$ Thus for $l\geq d$ we have $w = \bar{w}\wedge |M-s-lNL\rangle$ where $\bar{w} \in V_{s+lNL}^d.$ Since normally ordered wedges form a basis of $\F_M^d,$ the surjectivity of the map $\vro_l^d$ follows.

If $w,w^{\prime}\in V_{s+lNL}^d$ are two distinct normally ordered wedges, then  $w\wedge |M-s-lNL\rangle$ and $w^{\prime}\wedge |M-s-lNL\rangle$ are distinct and, as implied by the definition (\ref{eq:Vs}), normally ordered wedges in $\F_M.$ Thus $\vro_l^d$ is injective. \end{proof}

\begin{corol} \label{cor:vectisom}For each triple of non-negative integers $d,l,m$ such that $d\leq l < m$ the linear map
\bdm
\vro_{l,m}^d : V_{s+lNL}^d \rightarrow  V_{s+mNL}^d : \vro_{l,m}^d(w) = w\wedge u_{M-s-lNL}\wedge u_{M-s-lNL-1}\wedge \cdots \wedge u_{M-s-mNL+1}
\edm
is an isomorphism of vector spaces.
\end{corol}
\noindent\begin{proof}\hspace{5cm} $ \vro_{l,m}^d = (\vro_{m}^d)^{-1}\circ \vro_{l}^d.$ \end{proof}

\noindent{\bf Intertwining relations.} 
Let us recall (cf. Section \ref{sec:wedge}) that the Lie algebra $\sll_L$ acts on the wedge product $\Vaff^{\wedge n}.$ In particular, we have the action of Cartan subalgebra  ${\bar{\mathfrak h}}_L = \langle h_1,\dots,h_{L-1} \rangle $ which commutes with the Yangian action $\rho_{\kappa,\nu}^{(n)}.$ For each $a=1,\dots,L$ define $e_a \in {\bar{\mathfrak h}}_L$ as $h_a = e_a - e_{a+1},$ and let $ e_a^{(n)}$ be the operator that gives the action of $e_a$ on $\Vaff^{\wedge n}.$  

Let $x_1,\dots,x_L$ be a set of indeterminates, and for each $l=0,1,2,\dots $ define the following element of $\Cplx[x_1,\dots,x_L][[u^{-1}]]:$
\begin{align*}
\mbox{} & \Phi(u | x_1,\dots,x_L)_{s,l} = \prod_{a=1}^L \frac{u - \gamma_a(x) + N}{u - \gamma_a(x) },  \\
\mbox{} & \gamma_a(x) = \frac{s}{L} + (\kappa + N) l + \kappa (\un{o_1} + \delta(s >0)) + \nu(a) + \frac{1}{2} x_a,
\end{align*}
where the denominator is to be expanded as a series in $u^{-1}.$

\begin{propos} Let $0\leq d \leq l,$ and let $w \in V^d_{n:=s+lNL}.$ For all $p,q=1,\dots,N$ we have the following equality {\em ( intertwining relation )} in $V^d_{n + NL}[[u^{-1}]]:$
\begin{equation}
\rho_{\kappa,\nu}^{(n+NL)}\left(T_{pq}(u)\right)\cdot \vro_{l,l+1}^d(w)= \vro_{l,l+1}^d\left( \Phi(u|e_1^{(n)},\dots,e_L^{(n)})_{s,l}\cdot \rho_{\kappa,\nu}^{(n)}\left(T_{pq}(u)\right)\cdot w \right). \label{eq:intertwining}
\end{equation}
\end{propos}
\noindent \begin{proof} Set $n^{\prime}=n+NL.$ Let $w$ be a normally ordered wedge from $V^d_{n},$ and let  $\bar{w} =\vro_{l,l+1}^d(w).$ The vector   $\bar{w}$ is a normally ordered wedge from $V^d_{n^{\prime}},$ we have
\begin{alignat}{3}
& \bar{w} = u_{k_1}\wedge u_{k_2}\wedge \cdots \wedge u_{k_{n^{\prime}}}  = \wedge( f \otimes x \otimes y),  &\\
\intertext{where}
& f = (z_1^{\un{k_1}}z_2^{\un{k_2}}\cdots z_n^{\un{k_n}})(z_{n+1}\cdots z_{n^{\prime}})^m  \in C_z:=\Cplx[z_1^{\pm 1},\dots,z_{n+LN}^{\pm 1}] & \label{eq:fvec}\\ 
&x = (\vf_{\dot{k}_1}\vf_{\dot{k}_2} \cdots \vf_{\dot{k}_n})(\overbrace{\vf_1\vf_1\cdots\vf_1}^{\text{$N$ times}})\ldots (\overbrace{\vf_L\vf_L\cdots\vf_L}^{\text{$N$ times}})  \in C_L:=\left(\Cplx^L\right)^{\otimes n^{\prime}}& \label{eq:xvec}\\
&y = (\uf_{\ov{k_1}}\uf_{\ov{k_2}} \cdots \uf_{\ov{k_n}})(\uf_{N}\uf_{N-1} \overbrace{\cdots \uf_{1})\ldots (\uf_{N}\uf_{N-1}}^{\text{$L$ copies}} \cdots \uf_{1})  \in C_N:=\left(\Cplx^N\right)^{\otimes n^{\prime}}.& 
\end{alignat}
Here we dropped signs of tensor products and set $m:= \un{k_{n+1}}=\cdots = \un{k_{n^{\prime}}}.$ Observe that $m > \un{k_1},\dots,\un{k_n}.$

Define linear subspaces ${\mathcal L}_1,{\mathcal L}_2 \subset C_z\otimes C_L$ as follows:
\begin{alignat*}{3}
 &{\mathcal L}_1 =  \Cplx\{ z_1^{m_1}\cdots z_{n^{\prime}}^{m_{n^{\prime}}}\otimes x^{\prime} \: | \: x^{\prime}\in C_L, m_1,\dots,m_{n^{\prime}} \leq m, \#\{ m_i | m_i = m\} < NL \}  & \mbox{}\\
 &{\mathcal L}_2 = \Cplx\{(z_1^{m_1}\cdots z_n^{m_n})(z_{n+1}\cdots z_{n^{\prime}})^m \otimes (\vf_{a_1}\cdots \vf_{a_n})(\vf_{a_{n+1}}\cdots \vf_{a_{n^{\prime}}})\:|\: & \mbox{}\\
 & \qquad  \qquad  m_1,\dots,m_{n} < m, \exists a \in \{1,\dots,L\} \:\text{s.t.} \: \#\{a_i |  n < i \leq n^{\prime}, a_i = a\} > N  \} & \mbox{}
\end{alignat*}
\begin{lemma} Let  $y^{\prime} \in C_N,$ and let $v_1\in {\mathcal L}_1,$ $v_2\in {\mathcal L}_2.$ 
Then    
\begin{align}
\mbox{} & \wedge( v_1 \otimes y^{\prime} ) \in \oplus_{d^{\prime} > l} V_{n^{\prime}}^{d^{\prime}}, \tag{i}\\ 
\mbox{} & \wedge( v_2 \otimes y^{\prime} ) = 0.\tag{ii}
\end{align}
\end{lemma}
\noindent\begin{proof} \mbox{} \\
(i) The vector $\wedge( v_1 \otimes y^{\prime} )$ is a linear combination of normally ordered wedges 
\bdm
u_{(k_i)} = u_{k_1}\wedge u_{k_2} \wedge \cdots \wedge u_{k_n}\wedge u_{k_{n+1}}\wedge \cdots \wedge u_{k_{n^{\prime}}}
\edm
such that  $\un{k_{n+1}} < \un{o_{n+1}}.$ This inequality implies that $\deg u_{(k_i)} \geq l+1.$ \\ 
(ii) The vector $\wedge( v_2 \otimes y^{\prime} )$ is a linear combination of normally ordered wedges 
\bdm
u_{(k_i)} = u_{k_1}\wedge u_{k_2} \wedge \cdots \wedge u_{k_n}\wedge u_{k_{n+1}}\wedge \cdots \wedge u_{k_{n^{\prime}}}
\edm
such that there are at least two equal numbers among $k_{n+1},\dots,k_{n^{\prime}}.$ Hence $u_{(k_i)}=0.$ \end{proof}

\noindent
Let $d_i^{(n^{\prime})}$ be the matrix Dunkl operator acting on $C_z\otimes C_L.$ The operator $d_i^{(n)}$ naturally acts on $C_z\otimes C_L$ leaving the last $NL$ variables in $C_z$ and last $NL$ factors in $C_L$ untouched. 

For each $a=1,\dots,L,$ and each $i \in \{ n+(a-1)N+1,n+(a-1)N+2,\dots,n+aN\}$ we set
\bdm
\delta_i = \kappa m + \nu(a) + \frac{1}{2}e_a^{(n)} + \frac{n}{L} + n + aN - i.   
\edm
\begin{lemma}
Let $P(x_1,\dots,x_{n^{\prime}})$ be a polynomial in indeterminates $x_1,\dots,x_{n^{\prime}}.$ Let $f$ and $x$ be the vectors defined in {\em (\ref{eq:fvec})} and  {\em (\ref{eq:xvec})}. Then we have
\bdm
P(d_1^{(n^{\prime})},\dots,d_{n^{\prime}}^{(n^{\prime})})\cdot ( f \otimes x) \equiv P(d_1^{(n)},\dots,d_{n}^{(n)}, \delta_{n+1}, \cdots \delta_{n^{\prime}})\cdot ( f \otimes x)\bmod \left({\mathcal L}_1 + {\mathcal L}_2 \right).
\edm
\end{lemma}
\noindent \begin{proof}
We have the following relations 
\begin{align*}
\mbox{} & d_i^{(n^{\prime})} \cdot ( f \otimes x) \equiv \left(d_i^{(n)}  - \sum_{j=n+1}^{n^{\prime}} \sum_{a < b} (e_{ba})_i (e_{ab})_j \right)\cdot (f \otimes x)\bmod {\mathcal L}_1 \qquad ( i=1,\dots,n), \\
\mbox{} & d_i^{(n^{\prime})} \cdot ( f \otimes x) \equiv \left(\delta_i  + \sum_{j=1}^{n} \sum_{a < b} (e_{ab})_i (e_{ba})_j \right)\cdot (f \otimes x)\bmod {\mathcal L}_1 \qquad ( i=n+1,\dots,n^{\prime}).
\end{align*}
Observing that $d_i^{(n^{\prime})}\cdot {\mathcal L}_1 \subset {\mathcal L}_1$ we obtain the required statement. \end{proof}

\noindent Continuing with the proof of the proposition we find from the two preceding  lemmas and the definitions (\ref{eq:yangdef1}) and  (\ref{eq:yangdef2}) that (\ref{eq:intertwining}) holds modulo $\oplus_{d^{\prime} > l} V_{n^{\prime}}^{d^{\prime}}[[u^{-1}]].$ Both sides of  (\ref{eq:intertwining}), however, belong to  $ V_{n^{\prime}}^d[[u^{-1}]]$ since the Yangian actions preserve the degree. Hence $(\ref{eq:intertwining})$ holds exactly. \end{proof}

For each $l=0,1,2,\dots$ introduce the following element of $\End \Vaff^{\wedge n} [[u^{-1}]]:$
\begin{equation}
\varphi(u)_{s+lNL} = \prod_{m=0}^{l-1} \Phi(u|e_1^{(s+lNL)},\dots,e_L^{(s+lNL)})_{s,m} \qquad \left( \varphi(u)_s = 1\right), \label{eq:phi}
\end{equation}
and define $\bar{\rho}_{\kappa,\nu}^{(s+lNL)}\left(\YY(\gl_N)\right)$ as a renormalized Yangian action ${\rho}_{\kappa,\nu}^{(s+lNL)}\left(\YY(\gl_N)\right):$ 
\begin{equation}
\bar{\rho}_{\kappa,\nu}^{(s+lNL)}(T_{pq}(u)) = \frac{1}{\varphi(u)_{s+lNL}} \cdot {\rho}_{\kappa,\nu}^{(s+lNL)}(T_{pq}(u)).\label{eq:rhobar}
\end{equation}
Since $\sll_L$ leaves $V_{s+lNL}^d$ $(d\in \Zint_{\geq 0})$ invariant, so does the new action $\bar{\rho}_{\kappa,\nu}^{(s+lNL)}\left(\YY(\gl_N)\right).$ Continuing from  Corollary \ref{cor:vectisom} we have 
\begin{propos}\label{p:stab}
For each triple of non-negative integers $d,l,m$ such that $d\leq l < m$ the map $\vro_{l,m}^d$ is an isomorphism of $\YY(\gl_N)$-modules with actions $\bar{\rho}_{\kappa,\nu}^{(s+lNL)}$ and $\bar{\rho}_{\kappa,\nu}^{(s+mNL)}.$ 
\end{propos}
Recall the definition of the algebra $\UU^{\prime}(\bef^{\chi}) \hookrightarrow \UU(\asll_L)$ (cf. Section \ref{sec:DAHArep}).
\begin{propos}  \label{p:yinv1}
The statement of Proposition {\em \ref{p:yinv}} remains in force if $\rho_{\kappa,\nu}^{(n)}$ is replaced by $\bar{\rho}_{\kappa,\nu}^{(n)}.$
\end{propos}
\noindent\begin{proof} 
The action of the Cartan subalgebra $\bar{{\mathfrak h}}_L \subset \sll_L$ leaves $\UU^{\prime}(\bef^{\chi}) \Vaff^{\wedge n}$ invariant and commutes with  $\rho_{\kappa,\nu}(\YY(\gl_N)).$ The proposition now follows from the definition of  $\bar{\rho}_{\kappa,\nu}$ (cf. (\ref{eq:phi}) and (\ref{eq:rhobar})). \end{proof}

\subsubsection{Yangian action on the Fock space} \label{sec:YangonF}
\begin{defprop} \label{dp:YFock} Let $0\leq d \leq l,$ and define a Yangian action $\bar{\rho}_{\kappa,\nu}: \YY(\gl_N) \rightarrow \End\F_M^d$ as  
\bdm
 \bar{\rho}_{\kappa,\nu}\left(\YY(\gl_N)\right) = \vro_{l}^d \circ  \bar{\rho}_{\kappa,\nu}^{(s+lNL)}\left(\YY(\gl_N)\right) \circ \left(\vro_{l}^d\right)^{-1}.
\edm
This definition does not depend on the choice of $l$ in the right-hand side as long as $l \geq d.$
\end{defprop}
\noindent Thus a Yangian action is defined on each homogeneous component $\F_M^d,$ and hence on the entire linear space $\F_M.$

The following proposition is a semi-infinite counterpart of Proposition \ref{p:yinv}:

\begin{propos} \label{p:yinvFock}
Let $\kappa = L,$ and let $\nu(a) = a + \frac{\chi(a)}{2}$ $\;(a=1,\dots,L).$ Then the Yangian action $\bar{\rho}_{\kappa,\nu}(\YY(\gl_N))$ leaves invariant the subspace $\UU^{\prime}(\bef^{\chi})\F_M.$
\end{propos}
\noindent \begin{proof}
Let $v \in \F_M^d$ and let $l\geq d.$ From  Proposition \ref{p:VFockisom} it follows that there is a unique $v^{\prime} \in V^d_{s+lNL}$ where $s=M\pmod {NL}$ such that 
\bdm
v = v^{\prime}\wedge |M-s-lNL\rangle .
\edm
Let $g$ be one of the generators of $\UU^{\prime}(\bef^{\chi})=\langle h_1^{\chi},\dots,h_L^{\chi},f_0,\dots,f_{L-1}\rangle.$ For all large enough $l$ we have 
\bdm
g\cdot v = (g^{\prime}\cdot v^{\prime})\wedge |M-s-lNL\rangle
\edm
where $g^{\prime}$ stands for the action of $g$ on $\Vaff^{\wedge(s+lNL)}.$ If $g = f_0$ then $g^{\prime}\cdot v^{\prime} \in  V^{d+1}_{s+lNL},$ otherwise: $g^{\prime}\cdot v^{\prime} \in  V^{d}_{s+lNL}.$  

Let $a$ be an element of $\YY(\gl_N).$ Provided $l$ is large enough   Definition and Proposition \ref{dp:YFock} yields
\bdm 
\bar{\rho}_{\kappa,\nu}(a)\cdot g \cdot v =  \left(\bar{\rho}^{(s+lNL)}_{\kappa,\nu}(a)\cdot g^{\prime}\cdot v^{\prime}\right)\wedge |M-s-lNL\rangle.
\edm 
By Proposition $\ref{p:yinv1}$ the right-hand side above is a linear combination of vectors  
\bdm 
  (h^{\prime}\cdot w^{\prime})\wedge |M-s-lNL\rangle
\edm 
where $h$ is a generator of $\UU^{\prime}(\bef^{\chi})$ (we keep the meaning of the prime the same as above) and $w^{\prime}$ belongs to either $V^{d+1}_{s+lNL}$ or  $V^{d}_{s+lNL}.$ Taking again $l$ sufficiently large the last vector is seen to equal
\bdm
h \cdot (w^{\prime}\wedge |M-s-lNL\rangle)
\edm
which proves the proposition. \end{proof} 

\subsection{Yangian actions on irreducible integrable modules of $\agl_N$}
\subsubsection{From the Yangian action on the Fock space to Yangian actions on irreducible modules of  $\agl_N$}  \label{sec:Yonirr}
For each $\asll_N$ weight  $\Lambda \in P_N(L)^+$ we define in this section a Yangian action on the irreducible module (cf. Section \ref{sec:wedge}) 
\begin{equation}
  S_- \otimes V(\Lambda)
\end{equation}
of the Lie algebra $\agl_N = H\oplus \asll_N.$  

Let us go back to Theorem \ref{t:Fockdec}. There exists a unique $M \in \{0,1,\dots,N-1\}$ such that $ \Lambda \equiv \ov{\Lambda}_M \bmod \ov{Q}_N.$ 
Let $\omega_{M,\Lambda}$ be the corresponding dual weight of $\asll_L.$ Recall that  $\omega_{M,\Lambda} \neq \omega_{M,\Lambda^{\prime}}$ if $\Lambda \neq \Lambda^{\prime}.$ Define numbers $\chi_{M,\Lambda}(1),\dots,\chi_{M,\Lambda}(L)$ from (cf. Section \ref{sec:sl}):
\bdm
 \bar{\omega}_{M,\Lambda} = \sum_{a=1}^L \chi_{M,\Lambda}(L+1 - a) \bar{\vartheta}_a,
\edm 
and let $\UU^{\prime}({\mathfrak b}^{\chi_{M,\Lambda}})$ be the subalgebra of $\UU(\asll_L)$ associated to  $\chi_{M,\Lambda}(1),\dots,\chi_{M,\Lambda}(L)$ (cf. Section \ref{sec:DAHArep}).

Theorem  \ref{t:Fockdec} implies that we have the following isomorphism of  $\agl_N$-modules  
\bdm
 S_- \otimes V(\Lambda) \simeq \F_M / \UU^{\prime}({\mathfrak b}^{\chi_{M,\Lambda}})\F_M. 
\edm 
By Proposition \ref{p:yinvFock} the Yangian action $\bar{\rho}_{\kappa,\nu}(\YY(\gl_N))$ where $\kappa = L$ and $\nu(a) = a + \frac{\chi_{M,\Lambda}(a) }{2} $ $(a=1.\dots,L)$ factors through  the quotient map
\bdm
\F_M  \rightarrow \F_M / \UU^{\prime}({\mathfrak b}^{\chi_{M,\Lambda}})\F_M
\edm
and therefore defines a Yangian action on $ S_- \otimes V(\Lambda).$ We will denote this Yangian action by $\rho(\YY(\gl_N)).$

\subsubsection{Explicit expressions for the actions of the generators of $\YY(\sll_N)$}

With  notations of Section \ref{sec:wedge} (cf. proof of Proposition \ref{p:slacts}) for $n \in \Zint$ introduce the following operator on the Fock space:
\bdm
B(n) = \sum_{a=1}^L J^{aa}(n) = \sum_{s=1}^N J_{ss}(n).
\edm 
The operators $\{ B(n) | n \neq 0\}$ are generators of the Heisenberg algebra action on the Fock space $\F.$  The operator $B(0)$ commutes with all generators of $\agl_{NL}$ and acts on each charge component $\F_M$ $(M\in\Zint)$ of $\F$ as the multiplication by the  constant $M$. For $\alpha =1,\dots,\dim \sll_N$ and  $n \in \Zint$ set 
$$
 J_{\alpha}(n) = \sum_{s,t=1}^N (t_{\alpha})_{st} J_{st}(n), \quad J_{\alpha}^{ab}(n) = \sum_{s,t=1}^N (t_{\alpha})_{st} J_{st}^{ab}(n) \quad (a,b = 1,\dots,L).
$$ 
Here $(t_{\alpha})_{st}$ are the matrix elements introduced in Section \ref{sec:Ysl}. The operators $$\{ J_{\alpha}(n) \: | \: \alpha \in \{1,\dots,\dim \sll_N\}, n \in \Zint \}$$ are generators of the level $L$ action of $\asll_N$ on $\F.$ They satisfy the defining relations:
$$ [ J_{\alpha}(m) ,J_{\beta}(n) ] = f_{\alpha \beta \gamma} J_{\gamma}(m+n) + L m \delta_{\alpha \beta} \delta_{m+n,0} 1.
$$ 
Here and in what follows we assume summation over repeating Greek indices.
For $a=1,\dots,L$ and  $n \in \Zint$ set $I^{aa}(n) = {J}^{aa}(n) - \frac{1}{L}B(n).$ Then $ \sum_{a=1}^L I^{aa}(n) = 0$ and the operators $\{ J^{ab}(n),I^{aa}(n) | a,b=1,\dots,L; a\neq b; n\in \Zint \}$ are generators of the level $N$ action of $\asll_L$ on $\F.$   

\hfill 

\noindent Each charge component $\F_M,$ $M\in \Zint$ is an irreducible level 1  module of the Lie algebra $\agl_{NL}$ (cf. proof of Proposition \ref{p:slacts}). Therefore it is possible, in principle, to express the generators of the Yangian  action $\bar{\rho}_{\kappa,\nu}(\YY(\gl_N))$ defined in Section \ref{sec:YangonF} in terms of the generators of the $\agl_{NL}$-action. In practice the computation needed to obtain such  expressions is rather involved, and we have not carried it out for all generators of $\YY(\gl_N),$ but only for  generators of the subalgebra $\YY(\sll_N)$ of $\YY(\gl_N).$ To state the result introduce the following operators on $\F$: 
\begin{eqnarray}
&T_{\alpha} = \sum_{m \in \Zint}\sum_{a=1}^L \sum_{s,t=1}^N m \psi_s^a(m)^*  \psi_t^a(m) (t_{\alpha})_{st} \qquad & (\alpha = 1,\dots,\dim\sll_N),\label{eq:TD}\\
&W_{\alpha} = \frac{1}{2} \sum_{n \in \Zint} d_{\alpha\beta\gamma} \raisebox{-1pt}{\vdots} J_{\beta}(-n) J_{\gamma}(n) \raisebox{-1pt}{\vdots} &(\alpha = 1,\dots,\dim\sll_N).
\end{eqnarray}
Here the triple dots stand for the normal ordering (do not confuse with the fermion normal ordering $:\; :$ introduced in Section \ref{sec:wedge}):
$$ 
\raisebox{-1pt}{\vdots} J_{\beta}(-n) J_{\gamma}(n) \raisebox{-1pt}{\vdots} = \begin{cases} J_{\beta}(-n) J_{\gamma}(n) & \text{ if $ n \geq 0,$} \\ J_{\gamma}(n) J_{\beta}(-n)  & \text{ if $ n < 0.$} \end{cases}
$$ 
The $d_{\alpha\beta\gamma}$ is the $3$-symbol of $\sll_N$ defined from the following relation in $\End(\Cplx^N):$
$$ t_{\alpha} \cdot t_{\beta} = \frac{\delta_{\alpha\beta}}{N} \cdot 1 + \frac{1}{2} f_{\alpha\beta\gamma} t_{\gamma} + \frac{1}{2} d_{\alpha\beta\gamma} t_{\gamma}
$$  
The $d_{\alpha\beta\gamma}$ is a completely symmetric, traceless tensor, it satisfies the normalization condition 
$$  d_{\alpha\beta\gamma}d_{\lambda \beta\gamma} = \frac{ 2(N^2 - 4)}{N} \delta_{\alpha\lambda}. $$ 
Now for the generators $ Q_{\alpha}^{(1)}, Q_{\alpha}^{(2)}$ of the action of $\YY(\sll_N) \subset \YY(\gl_N)$ (cf. Section \ref{sec:Ysl},(\ref{eq:embed})) on a charge component $\F_M$ of the Fock space $\F$ we have the following  expressions:
\begin{gather}
\begin{aligned} \mbox{}
&\bar{\rho}_{\kappa,\nu} \left( Q_{\alpha}^{(1)} \right)  = J_{\alpha}(0),  \\ 
\mbox{}&\bar{\rho}_{\kappa,\nu} \left( Q_{\alpha}^{(2)} \right) = \left(\frac{M}{N} - \frac{N}{4} \right) J_{\alpha}(0) + \label{eq:Q2Fock}
\end{aligned} \\  
+ (L - \kappa) T_{\alpha}  + \nonumber \\ 
+ \sum_{a=1}^L \left( \frac{I^{aa}(0)}{2} + a - \nu(a)\right) J^{aa}_{\alpha}(0) + \sum_{1 \leq a < b \leq L} J^{ab}(0)J^{ba}_{\alpha}(0) + \nonumber \\
+\sum_{n > 0} \left( \sum_{1\leq a \neq b \leq L} J^{ab}(-n) J_{\alpha}^{ba}(n) +  \sum_{1\leq a \leq L} I^{aa}(-n) J_{\alpha}^{aa}(n) \right) + \nonumber \\
+ \left( \frac{1}{N} + \frac{1}{L}\right) \sum_{n > 0} B(-n)J_{\alpha}(n)  + \frac{1}{N}  \sum_{n > 0} J_{\alpha}(-n) B(n) + \nonumber \\
+ \frac{1}{2} f_{\alpha\beta\gamma} \sum_{n > 0} J_{\beta}(-n) J_{\gamma}(n) +  \frac{1}{2} W_{\alpha}.\nonumber 
\end{gather}
In Section \ref{sec:DAHArep} we associated with every sequence of complex numbers $\chi(a) $ $(a=1,\dots,L)$ the subalgebra $\UU^{\prime}( {\mathfrak b}^{\chi} )$ of $\UU(\asll_N).$ The action of $\UU^{\prime}( {\mathfrak b}^{\chi} )$ on the Fock space is generated by the operators 
\begin{align*} & \{ I^{aa}(0) - \chi(a) 1 , J^{cd}(0) \: | \:  1\leq a \leq L; 1\leq c < d \leq L \}, \\ 
&  \{  I^{aa}(n), J^{cd}(n) \: | \:  1\leq a \leq L; 1\leq c \neq d \leq L ; n < 0 \}.  
\end{align*}
The expressions (\ref{eq:Q2Fock}) show that with the parameters of the Yangian action $\bar{\rho}_{\kappa,\nu}$ fixed as  $\kappa = L$ and $ \nu(a) = a + \frac{\chi(a)}{2},$ the operators $\bar{\rho}_{\kappa,\nu}(Q_{\alpha}^{(1)})$ and $\bar{\rho}_{\kappa,\nu}(Q_{\alpha}^{(2)})$ leave the subspace $\UU^{\prime}( {\mathfrak b}^{\chi} ) \F_M$ invariant. This gives an alternative proof of Proposition \ref{p:yinvFock}  for the action of the subalgebra $\YY(\sll_N).$ Following the preceding section we obtain the $\YY(\gl_N)$-action $\rho$ on each irreducible $\agl_N$-module $S_-\otimes V(\Lambda).$ From (\ref{eq:Q2Fock}) it follows that the generators of the subalgebra $\YY(\sll_N)$ act on  $S_-\otimes V(\Lambda)$ as:
\begin{gather}
\begin{aligned}
\mbox{}&{\rho}\left( Q_{\alpha}^{(1)} \right)  = J_{\alpha}(0), \\
\mbox{}&{\rho}\left( Q_{\alpha}^{(2)} \right) = \left(\frac{M}{N} - \frac{N}{4} \right) J_{\alpha}(0) + \label{eq:Q2ir}
\end{aligned} \\ 
+ \left( \frac{1}{N} + \frac{1}{L}\right) \sum_{n > 0} B(-n)J_{\alpha}(n)  + \frac{1}{N}  \sum_{n > 0} J_{\alpha}(-n) B(n) + \nonumber \\
+ \frac{1}{2} f_{\alpha\beta\gamma} \sum_{n > 0} J_{\beta}(-n) J_{\gamma}(n) +  \frac{1}{2} W_{\alpha}.\nonumber 
\end{gather}
Here we understand that the symbols $B(n), J_{\alpha}(n)$ stand for the operators giving the canonical action of $\agl_N$ on  $S_-\otimes V(\Lambda).$ Notice that the expression for ${\rho}\left( Q_{\alpha}^{(2)} \right)$ contains $M$ which is an arbitrary integer such that $ \ov{\Lambda} \equiv \ov{\Lambda}_{M\bmod N} \bmod \ov{Q}_N.$ The dependence of the Yangian action $\rho$ on $M$ is trivial: the term $(M/N - N/4) J_{\alpha}(0)$ in (\ref{eq:Q2ir})  can be eliminated by performing the automorphism
$$ Q_{\alpha}^{(1)} \mapsto  Q_{\alpha}^{(1)}, \quad  Q_{\alpha}^{(2)} \mapsto  Q_{\alpha}^{(2)} - (M/N - N/4) Q_{\alpha}^{(1)}  $$
of the algebra $\YY(\sll_N).$ 

\noindent{\bf Yangian actions on irreducible integrable Level 1 modules of $\agl_N$ and $\asll_N$.} The results obtained thus far apply with obvious modifications  for the case when the level $L$ equals 1. In this case each charge component of the Fock space is by itself an irreducible module of $\agl_N.$ More precisely we  have the isomorphism of  $\agl_N$-modules: 
$$ \F_M \simeq  S_-\otimes V(\Lambda_{M\bmod N}). $$ 
We would now like to point out the relation of our construction of the Yangian actions on the irreducible integrable modules of $\agl_N$ to the $\YY(\asll_N)$-actions on irreducible integrable level 1 modules of $\asll_N$ obtained by \cite{HHTBP} and  \cite{Schoutens}. 

\hfill 

\noindent When $L=1$ the operator $T_{\alpha}$ (cf. (\ref{eq:TD})) can be written as follows:
$$ T_{\alpha} = -\frac{1}{2} J_{\alpha}(0) - \frac{1}{N} \sum_{n\in \Zint} B(-n) J_{\alpha}(n)  - \frac{1}{N+2} W_{\alpha}. $$ 
Therefore  the expressions (\ref{eq:Q2Fock}) assume the form:
\begin{gather}
\begin{aligned} \mbox{} &\bar{\rho}_{\kappa,\nu} \left( Q_{\alpha}^{(1)} \right)  = J_{\alpha}(0),  \\
\mbox{} & \bar{\rho}_{\kappa,\nu} \left( Q_{\alpha}^{(2)} \right) = c_{\kappa,\nu}(M,N) \ \: J_{\alpha}(0) +  \label{eq:Q2FockL1} \end{aligned}\\
+ \left( \frac{\kappa}{N} + 1 \right) \sum_{n > 0} B(-n)J_{\alpha}(n)  + \frac{\kappa}{N}  \sum_{n > 0} J_{\alpha}(-n) B(n) + \nonumber  \\
+ \frac{1}{2} f_{\alpha\beta\gamma} \sum_{n > 0} J_{\beta}(-n) J_{\gamma}(n) +  \frac{2\kappa + N}{2(N+2)} W_{\alpha}\nonumber \\
\text{where}\qquad  c_{\kappa,\nu}(M,N) = \frac{\kappa +1}{2} + \frac{\kappa M}{N} - \frac{N}{4} - \nu(1) .\nonumber 
\end{gather}
Clearly the operators 
\begin{align*}
&\bar{\rho}_{\kappa} \left( Q_{\alpha}^{(1)} \right) \stackrel{{\mathrm {def}}}{=} \bar{\rho}_{\kappa,\nu} \left( Q_{\alpha}^{(1)} \right), \\
& \bar{\rho}_{\kappa} \left( Q_{\alpha}^{(2)} \right) \stackrel{{\mathrm {def}}}{=}  \bar{\rho}_{\kappa,\nu} \left( Q_{\alpha}^{(2)} \right) - c_{\kappa,\nu}(M,N) \ \: J_{\alpha}(0)
\end{align*}
still  satisfy the defining relations of $\YY(\sll_N).$  

\hfill 

\noindent The irreducible level 1 $\asll_N$-module $V(\Lambda_{M})$ $(M=0,\dots,N-1)$ is isomorphic to the quotient of $S_-\otimes V(\Lambda_{M})$ by the image of the algebra $H_- = \langle B(n) \: | \: n < 0 \rangle .$ From the expressions (\ref{eq:Q2FockL1}) it follows that with $\kappa = 0$ the operators $\bar{\rho}_{\kappa} \left( Q_{\alpha}^{(1)} \right)$ and $ \bar{\rho}_{\kappa} \left( Q_{\alpha}^{(2)} \right)$ leave the image of $H_-$ invariant. Therefore an action of $\YY(\sll_N)$ is induced on  $V(\Lambda_{M}).$ Denote this action by $\rho_-.$ Then 
\begin{equation}
\begin{aligned}
 \mbox{}&  \rho_- \left(Q_{\alpha}^{(1)}\right) =   J_{\alpha}(0),  \\
 \mbox{}&  \rho_- \left(Q_{\alpha}^{(2)}\right) = 
 \frac{1}{2} f_{\alpha\beta\gamma} \sum_{n > 0} J_{\beta}(-n) J_{\gamma}(n) +  \frac{N}{2(N+2)} W_{\alpha}.\label{eq:Q2irL1}
\end{aligned}
\end{equation}
On the other hand $V(\Lambda_{M})$ $(M=0,\dots,N-1)$ is isomorphic to the kernel of the algebra $H_+ = \langle B(n) \: | \: n > 0 \rangle $ in $S_-\otimes V(\Lambda_{M}).$ 
From the expressions (\ref{eq:Q2FockL1}) it follows that if  $\kappa = -N,$ the operators $\bar{\rho}_{\kappa} \left( Q_{\alpha}^{(1)} \right)$ and $ \bar{\rho}_{\kappa} \left( Q_{\alpha}^{(2)} \right)$ leave the kernel  of $H_+$ invariant. Therefore an action of $\YY(\sll_N)$ is induced on  $V(\Lambda_{M}).$ Denote this action by $\rho_+.$ Then we have 
\begin{equation}
\begin{aligned}
 \mbox{}&  \rho_+ \left(Q_{\alpha}^{(1)}\right) =   J_{\alpha}(0),  \\
 \mbox{}&  \rho_+ \left(Q_{\alpha}^{(2)}\right) = 
 \frac{1}{2} f_{\alpha\beta\gamma} \sum_{n > 0} J_{\beta}(-n) J_{\gamma}(n) -  \frac{N}{2(N+2)} W_{\alpha}.\label{eq:Q2irL12}
\end{aligned}
\end{equation}
These  expressions for the Yangian generators coincide with those of \cite{Schoutens}.

\section{Yangian decompositions of vacuum modules of $\agl_N$} 

\subsection{Regular elements of the affine Weyl group} \label{se:regel}

\subsubsection{Intertwiners}
In this section we  summarize several facts concerning intertwiners of weight spaces in modules of the degenerate double affine Hecke algebra $\daha$ \cite{AST}. For an $\daha$-module $V$ and $\zeta \in (\tf^{\prime})^*$ (cf. Section \ref{sec1:daha}) define the weight space $V_{\zeta}$ with respect to the action of $S[\tf^{\prime}] \hookrightarrow \daha$ as     
\bdm
V_{\zeta} = \{ v \in V\: | \: \xi\cdot v = \zeta(\xi) v \: \text{ for any $ \xi \in \tf^{\prime}$}\}.
\edm
Define the following elements of $\daha:$
\begin{align*}
 & \ph_i  = 1 + s_i \alpha_i^{\vee} \qquad (i \in \{0,\dots,n-1\} \cong \Zint/n\Zint ), \\
& \ph_{\pi}  = \pi. 
\end{align*}
For any $\xi \in  \tf^{\prime}$ and $i \in \Zint/n\Zint$ one has 
\begin{equation}
\left.\begin{aligned} 
\mbox{}&\ph_i\cdot \xi = s_i(\xi) \cdot \ph_i, \\
\mbox{}&\ph_{\pi}\cdot \xi = \pi(\xi) \cdot \ph_{\pi}.
\end{aligned}   \right. \label{eq:phixi}
\end{equation}
\begin{propos}
The elements $\ph_{\pi}$ and $\ph_i$ $(i\in \Zint/n\Zint)$ satisfy the following relations: 
\begin{align*}
&\ph_i \cdot \ph_{i+1} \cdot \ph_i =  \ph_{i+1} \cdot \ph_i \cdot \ph_{i+1}, \\
& \ph_i \cdot \ph_j = \ph_j\cdot \ph_i \quad (i-j \not\equiv \pm 1 \bmod n ),\\
& \ph_{\pi} \cdot \ph_i = \ph_{i+1}\cdot \ph_{\pi},\\
& \ph_i^2 = 1 - {\alpha_i^{\vee}}^2.
\end{align*}
\end{propos}
\noindent For $w = \pi^k \cdot s_{i_1}\cdots s_{i_l} \in W,$ where the right-hand side is a reduced expression, define the {\em intertwiner } associated with $w$ as 
\bdm 
\ph_{w} = \ph_{\pi}^k \cdot \ph_{i_1}\cdots \ph_{i_l} \in \daha.
\edm
By the above proposition $\ph_{w}$ does not depend on the choice of a reduced expression, and by (\ref{eq:phixi}) we have
\bdm
\ph_{w}\cdot \xi = w(\xi)\cdot \ph_{w} \quad (\xi \in \tf^{\prime}).
\edm
Which leads to (cf. (\ref{eq:Wact}) for  the action of $W$ on $(\tf^{\prime})^*$):
\begin{propos} \label{p:weightofph}
Let $V$ be an $\daha$-module. Let $\zeta \in (\tf^{\prime})^*$ and $w\in W.$ Then $ \ph_{w} \cdot v \in V_{w(\zeta)}$ for any $v \in V_{\zeta}.$
\end{propos}
\begin{propos}[\cite{AST} Proposition 1.4.3] \label{p:sqint}
For $w \in W,$  
\begin{align}
&\ph_{w^{-1}}\cdot \ph_w = \prod_{\alpha \in S(w)} (1 - {\alpha^{\vee}}^2), \tag{i} \label{eq:phi2}\\ 
&\ph_w = w \cdot \prod_{\alpha \in S(w)} \alpha^{\vee} + \sum_{u \prec w } u \cdot f_u, \tag{ii}  \label{eq:hiterm} 
\end{align}
where $f_u \in S[\tf^{\prime}].$ 
\end{propos}

\subsubsection{Genericity and irreducibility}
In Section \ref{sec:DAHArep} a representation $\pi_{\kappa,\nu}$ of the degenerate double affine Hecke algebra $\daha$ was defined on the linear space 
\begin{equation*}
\Ve = \Cplx[z_1^{\pm 1},\dots,z_n^{\pm 1}] \otimes \left(\Cplx^L\right)^{\otimes n} \cong  \left( \Cplx[z^{\pm 1}] \otimes \Cplx^L\right)^{\otimes n}. 
\end{equation*}
Recall that each weight component with respect to the natural action of $\sll_L$ on $\Ve$ is left invariant by $\pi_{\kappa,\nu}(\daha).$ Let $n$ be divisible by $L:$ $n = Lm,$ where $m$ is a positive integer. \mbox{} From now on we will always assume that $n$ is of this form. Let $\Ve_0$ be the linear subspace of $\Ve$ of $\sll_L$ weight zero. The following vector        
\begin{equation}
v_0 = \underbrace{\vf_1\otimes \cdots \otimes \vf_1}_m \otimes \underbrace{\vf_2\otimes \cdots \otimes \vf_2}_m \otimes \cdots \otimes \underbrace{\vf_L\otimes \cdots \otimes \vf_L}_m
\end{equation}
is easily seen to be a {\em cyclic vector} of $\Ve_0$ with respect to the action of $\daha:$ $\Ve_0 = \pi_{\kappa,\nu}(\daha)\cdot v_0.$ Moreover, $v_0$ is a weight vector of $S[\tf^{\prime}] \hookrightarrow \daha$ of the weight  
\begin{equation}
\zeta_0 = \kappa c^* + \sum_{a=1}^L \sum_{i= (a-1)m + 1}^{a m} \left( i - \nu(a) - am \right) \ep_i  \in (\tf^{\prime})^*.  \label{eq:hweight}
\end{equation}
Let $\beta = (\beta_1,\dots,\beta_L)$ be the rectangular partition $(m,m,\dots,m)$ of $n,$ and define 
\begin{align*}
&\bar{R}_{\beta} = \{ \alpha_{ij} \in \bar{R}\:|\: (a-1)m < i,j \leq a m \;\text{for some $a=1,\dots,L$}\},\\ 
&\bar{R}_{\beta +} = \bar{R}_+ \cap  \bar{R}_{\beta},   \\
& W^{\beta} = \{ w \in W \: | \: S(w) \subset R_+\setminus \bar{R}_{\beta +}\},\\
& \bar{W}^{\beta} = \{ w \in \bar{W} \: | \: S(w) \subset \bar{R}_+\setminus \bar{R}_{\beta +}\}.
\end{align*}
\begin{defin}
A weight $\zeta \in (\tf^{\prime})^*$ of the form {\em (\ref{eq:hweight})} is said to be {\em generic} if and only if $(\zeta,\alpha) \not\in \{-1,0,1\}$ for every $\alpha \in \cup_{w \in W^{\beta}} S(w).$
\end{defin}
The following fact has been proved in \cite{AST}: 
\begin{propos}
If $\zeta_0 \in (\tf^{\prime})^*$ is generic, then $\Ve_0$ is an irreducible $\daha$-module with a basis $\{ \ph_{w} \cdot v_0 \:| \: w \in W^{\beta}\}.$  
\end{propos}
\mbox{} From now on we  fix  values of the parameters $\kappa$ and $\nu(1),\dots,\nu(L)$ as $\kappa = L$ and $\nu(a) = a$ $(a=1,\dots,L).$ The weight (\ref{eq:hweight}) that corresponds to these parameters is clearly not generic since (cf. Proposition \ref{p:dahainv}) $\Ve_0$ contains the proper invariant subspace $\Ve_0 \cap \UU^{\prime}({\mathfrak b}^0)\Ve.$  

\noindent{\bf The set $\bar{W}^{\beta}.$} In what follows we will find it convenient to label elements of the set $\bar{W}^{\beta}$ by tableaux of the diagram associated with the partition $\beta.$ For all $i=1,\dots,n$ and $w\in \bar{W}$ define  $w(i)$ by $w(\ep_i) = \ep_{w(i)}.$ To each $w \in \bar{W}^{\beta}$ we associate the tableau  
\setlength{\unitlength}{1.5pt}
\begin{equation*} {\mathrm {tab}}_w = \quad
\begin{picture}(250,28)(0,20)
\put(0,50){\line(1,0){110}}\put(160,50){\line(1,0){60}}
\put(0,40){\line(1,0){110}}\put(110,40){\makebox(50,10){$\cdots$}}\put(160,40){\line(1,0){60}}
\put(0,30){\line(1,0){110}}\put(110,30){\makebox(50,10){$\cdots$}}\put(160,30){\line(1,0){60}}
\put(0,10){\line(1,0){110}}\put(160,10){\line(1,0){60}}
\put(0,0){\line(1,0){110}}\put(110,0){\makebox(50,10){$\cdots$}}\put(160,0){\line(1,0){60}}
%%%%%%%%%%%%%%%%%%%%%%%%%%%%%%%
\put(0,0){\line(0,1){12}}\put(50,0){\line(0,1){12}}\put(100,0){\line(0,1){12}}\put(170,0){\line(0,1){12}}\put(220,0){\line(0,1){12}}
\put(0,50){\line(0,-1){22}}\put(50,50){\line(0,-1){22}}\put(100,50){\line(0,-1){22}}\put(170,50){\line(0,-1){22}}\put(220,50){\line(0,-1){22}}
\put(0,10){\makebox(50,20){$\cdots$}}\put(170,10){\makebox(50,20){$\cdots$}}
%%%%%%%%% entries %%%%%%%%%%%%%
\put(0,40){\makebox(50,10){\tiny $ w(1)$}}\put(50,40){\makebox(50,10){\tiny $ w(2)$}}\put(170,40){\makebox(50,10){\tiny $ w(m)$}}
\put(0,30){\makebox(50,10){\tiny $ w(m+1)$}}\put(50,30){\makebox(50,10){\tiny $ w(m+2)$}}\put(170,30){\makebox(50,10){\tiny $ w(2m)$}}
\put(0,0){\makebox(50,10){\tiny $ w((L-1)m+1)$}}\put(50,0){\makebox(50,10){\tiny $ w((L-1)m+2)$}}\put(170,0){\makebox(50,10){\tiny $ w(Lm)$}}
\end{picture}
\end{equation*} \\
\setlength{\unitlength}{1.2pt}
\mbox{}%\vspace{15pt}
\\
By definition of the set $\bar{W}^{\beta}$ the numbers inscribed into ${\mathrm {tab}}_w$ increase along the rows from left to right. The assignment $ w \mapsto {\mathrm {tab}}_w$ defines a one-to-one correspondence between $\bar{W}^{\beta}$ and tableaux of $\beta$ with this property. Another parameterization of  $\bar{W}^{\beta}$ is occasionally useful. 
We will say that  $\bb = (b_1,\dots,b_n) \in \{1,\dots,L\}^n$ is a {\em weight 0  spin configuration}  if the multiplicity of each $a=1,\dots,L$ in $\bb$ is $m.$ To each $w\in \bar{W}^{\beta}$ we associate a weight 0 spin configuration $\bb$ by setting $b_i = a$ if $i$ is inscribed into the row of depth $a$ in ${\mathrm {tab}}_w.$ 
The correspondence between $\bar{W}^{\beta}$ and the set of  weight 0 spin configurations defined by this rule is clearly one-to-one, we denote by $w_{\bb}$ the image of $\bb$ under the inverse correspondence.     

\subsubsection{Regular elements of the affine Weyl group}
Let $r=(r_1,\dots,r_n)$ be a sequence of integers such that $r_1 \leq r_2 \leq \cdots \leq r_n.$ With each such sequence  we associate unique  $a = (a_i) \in \{1,\dots,L\}^n$ and $\lambda = (\lambda_i) \in \Zint^n$ by $r_i = a_i - L \lambda_i.$ As a shorthand we write $r = a - L\lambda.$ Observe that $\lambda_1 \geq \lambda_2 \geq \cdots \geq \lambda_n$ and $a_{i} \leq  a_{i+1}$ if $\lambda_i = \lambda_{i+1}.$ 
Let ${\mathcal R}_0$ be the set of all non-decreasing sequences $r=a - L\lambda$ such that $a$ is a weight 0 spin configuration.

For each $\eta \in \bar{P}$ define $\gamma_{\eta}$ as the element of $\bar{W}$ with shortest possible length such that $\gamma_{\eta}(\eta) \in \bar{P}_- $ $=$ $\{ \eta \in \bar{P}\:|\: (\eta,\alpha) \leq 0 \; \text{for $\alpha \in \bar{R}_+$}\}.$ 
For $a = (a_i) \in \{1,\dots,L\}^n$ and $w\in \bar{W}$ we define $w(a)$ as $(a_{w^{-1}(i)}).$  
To each  $r = a - L\lambda \in {\mathcal R}_0$ we associate an element of the affine Weyl group  as follows: 
\begin{equation}
x_r = t_{\lambda}\cdot \gamma_{\lambda}^{-1} \cdot w_{\gamma_{\lambda}(a)}, \label{eq:xr}
\end{equation}
where $\lambda = \sum_{i=1}^n \lambda_i \ep_i \in \bar{P}_+.$ Observe that $x_r \in W^{\beta}$ (cf. Lemma 2.2.4 in \cite{AST}).
\begin{defin}
Let $r \in {\mathcal R}_0.$ The sequence $r$ and the element $x_r$ are called {\em regular} if and only if $(\zeta_0,\alpha) \not\in \{-1,0,1\}$ for each $\alpha \in S(x_r).$ 
\end{defin}
\noindent Proposition \ref{p:sqint} (\ref{eq:phi2}) shows that if $x_r$ is regular, then $\pi_{\kappa,\nu}(\ph_{x_r}) \cdot v_0$ does not belong to the invariant subspace $\UU^{\prime}({\mathfrak b}^0)\Ve$ because $v_0 \not\in  \UU^{\prime}({\mathfrak b}^0)\Ve.$ 

Let ${\mathcal R}_0^{\reg}$ be the subset of all regular sequences in ${\mathcal R}_0.$ In the next paragraph we give a combinatorial description of ${\mathcal R}_0^{\reg}.$ 
 
\subsubsection{Parameterization of regular elements by skew Young diagrams} \label{sec:regdiag}
Recall that a tableau of a skew Young diagram is called {\em standard} if numbers inscribed in this tableau increase along each row from left to right and along each column from top to bottom.
\begin{lemma} \label{l:reg1}
Let $r = a -L \lambda \in {\mathcal R}_0,$ and let $\eta = \gamma_{\lambda}(\lambda),$ $w=w_{\gamma_{\lambda}(a)}$ {\em (cf.(\ref{eq:xr})).} Then $r$ is regular if and only if the following two conditions are satisfied:  
\begin{align}
& {\mathrm {tab}}_w \;\text{ { is standard.}} \\
& \eta_{w((L-1)m+i)} - \eta_{w(i)} \in \{0,1\} \; \text{{ for all} $i=1,\dots,m.$} \label{eq:condii}
\end{align}
\end{lemma}
%%%%%%%%%%%%%%%%%%%%%%%%%%%%%%%%%%%%%%%%%%%%%%%%%%%%%%%%%%%%%%%%%%%%%%%%%%%%%%%%%%
\noindent\begin{proof}
A computation gives $S(x_r) = S_1\sqcup S_2,$ where 
\begin{align}
& S_1 = \bigsqcup_{\{ \alpha \in \bar{R}_+\setminus \bar{R}_{\beta +}\:|\: w(\alpha) \in  \bar{R}_-\}} 
\{ \alpha + k\delta \:|\: 0 \leq k \leq (\eta,w(\alpha))-1\}, \\
& S_2 = \quad \bigsqcup_{\{ \alpha \in \bar{R}_+\:|\: w(\alpha) \in  \bar{R}_+\}} 
\{ -\alpha + k\delta \:|\: 1 \leq k \leq -(\eta,w(\alpha))-1\}. 
\end{align}
For $i=1,2$ let  $\zeta_0(S_i) = \cup_{\alpha \in S_i}\{ (\zeta_0,\alpha )\}.$ Explicitly:
\bdm
\zeta_0(S_1)=\bigcup^{1\leq c < d \leq L} \!\!\!\! \bigcup^{1\leq i,j\leq m }_{ w((c-1)m+i)>w((d-1)m+j)}
\!\!\!\!\!\!\!\!\!\!\!\!\!\!\!\!\!\!\!\!\!\! \{ d - c + i-j + Lk \:|\: 0 \leq k < \eta_{w((c-1)m+i)} - \eta_{w((d-1)m+j)}\}.
\edm
Note that in the last formula we have $\eta_{w((c-1)m+i)} - \eta_{w((d-1)m+j)} > 0$ because $r$ is a non-decreasing sequence. Suppose ${\mathrm {tab}}_w$ is standard. Then for $1\leq i,j \leq m,$ $c<d$ the inequality $w((c-1)m+i)>w((d-1)m+j)$ implies $i>j.$ Hence all elements of $\zeta_0(S_1)$ are greater or equal 2. Suppose ${\mathrm {tab}}_w$ is not standard, i.e. there are $c$ and $i$ such that $w((c-1)m+i)>w(cm+i).$ Taking $j=i,$ $d=c+1,$ $k=0$ in the formula for $\zeta_0(S_1)$ we see that $\zeta_0(S_1)$ contains 1. Thus $-1,0,1 \not\in \zeta_0(S_1)$ if and only if  ${\mathrm {tab}}_w$ is standard.

Consider now  $\zeta_0(S_2).$ We have $\zeta_0(S_2)=\zeta_0(S_2)^{\prime}\cup\zeta_0(S_2)^{\prime\prime}$ where the set 
\bdm
\zeta_0(S_2)^{\prime} = \bigcup_{c=1}^L \bigcup_{(c-1)m<i<j\leq cm}\{j-i+Lk\:|\: 1\leq k < \eta_{w(j)}-\eta_{w(i)}\}   
\edm 
is either empty or else contains only integers greater than $L,$ and where 
\begin{align}
& \zeta_0(S_2)^{\prime\prime}=\bigcup^{1\leq c < d \leq L} S_{c,d}^+ \cup S_{c,d}^-, \nonumber \\ 
& S_{c,d}^+ =  \bigcup^{1\leq i <j\leq m }_{ w((c-1)m+i)<w((d-1)m+j)}
\!\!\!\!\!\!\!\!\!\!\!\!\!\!\!\!\!\!\!\!\!\!\! \{ c - d + j-i + Lk \:|\: 1 \leq k < \eta_{w((d-1)m+j)} -\eta_{w((c-1)m+i)} \},\nonumber \\
& S_{c,d}^- =  \bigcup^{m \geq i \geq j \geq 1 }_{ w((c-1)m+i)<w((d-1)m+j)}
\!\!\!\!\!\!\!\!\!\!\!\!\!\!\!\!\!\!\!\!\!\!\! \{ c - d + j-i + Lk \:|\: 1 \leq k < \eta_{w((d-1)m+j)} -\eta_{w((c-1)m+i)} \}. \label{eq:Smin}
\end{align}
For all $1\leq c < d \leq L$ the set  $S_{c,d}^+$ is either empty or else contains only integers that are greater or equal $2.$ Consider $S_{c,d}^- .$
Assume until the end of the proof that ${\mathrm {tab}}_w$ is standard. Then for any $i=1,\dots,m$ and $1\leq c < d\leq L$ we have $w(i) \leq w((c-1)m+i)<w((d-1)m+i)\leq w((L-1)m+i),$ which implies 
\begin{equation} \eta_{w(i)} \leq \eta_{w((c-1)m+i)} \leq  \eta_{w((d-1)m+i)}\leq \eta_{w((L-1)m+i)} \label{eq:ineqqqq}
\end{equation}
 because $\eta \in \bar{P}_-.$ For all $1\leq c < d \leq L$ and $1\leq j \leq i \leq m$ that appear in (\ref{eq:Smin}) we have $w((c-1)m+i)<w((d-1)m+j)\leq w((d-1)m+i).$  Suppose the condition (\ref{eq:condii}) is satisfied. Then the last chain of inequalities and (\ref{eq:ineqqqq}) give  $\eta_{w((d-1)m + j)} - \eta_{w((c-1)m + i)} \in \{0,1\}$ which implies that $S_{c,d}^-$ is empty.

Now suppose that $(\ref{eq:condii})$ is not satisfied: there is $i$ such that $\eta_{w((L-1)m+i)} - \eta_{w(i)} > 1.$ Then $S^{-}_{1,L}$ contains 1.  \end{proof}

With every spin configuration  $b=(b_i) \in \{1,\dots,L\}^n$ we associate a sequence $p(b) = (p_1,\dots,p_n)$ of $\sll_L$ weights called the {\em path} of $b,$ setting (cf. Section \ref{sec:sl})
\bdm    
p_i = \bar{\vartheta}_{b_1}+\bar{\vartheta}_{b_2}+\cdots +\bar{\vartheta}_{b_i}.
\edm 
A path $p(b)$ is called {\em dominant} if each $p_i$ is a dominant integral weight.
\begin{corol} \label{cor:cor2}
Let $r \in {\mathcal R}_0.$ With notations of {\em Lemma \ref{l:reg1}} let $b = \gamma_{\lambda}(a).$ Then  $r$ is regular if and only if the following two conditions are satisfied: 
\begin{align}
& \text{The path $p(b)$ is dominant.} \label{eq:coni}\\
& \text{For each pair $1\leq i< j \leq n$ such that $b_i=1,b_j=L$ }\label{eq:conii}\\ 
&\text{and $\#\{k\leq i\:|\:b_k = 1\} =\#\{k\leq j\:|\:b_k = L\}$ one has  $\eta_j - \eta_i \in \{0,1\}.$} \nonumber 
\end{align}
\end{corol}
\noindent We will call any subset $D$ of $\Zint^2$ a {\em diagram} and employ the usual graphic representation of $D:$ a point $(i,j)\in \Zint^2$ is represented by the unit square in the plane $\Rea^2$ with the centre $(i,j),$ the coordinates $i$ and $j$ on  $\Rea^2$ increasing from top to bottom and from left to right respectively.

With any element $r=(r_i)=a-L\lambda$ of ${\mathcal R}_0$ we associate the diagram $D_r$ as follows: for all $i=1,\dots,n$ set $n_i = \#\{j\leq i \:|\: a_i=a_j, \lambda_i=\lambda_j\} +\#\{j > i \:|\: a_i=a_j, \lambda_i<\lambda_j\},$ and define  
\bdm
D_r = \bigsqcup_{i=1}^n \{ (n_i,r_i)\}.
\edm
Observe that: (i) $(n_i,r_i) \neq (n_j,r_j)$ if $i\neq j;$ (ii) intersection of $D_r$ with the $i$th row $\{ (i,j)\:|\: j \in \Zint \}$ is empty if $i<1$ or $i>m,$ and contains exactly $L$ squares if $1\leq i\leq m;$ (iii)  intersection of $D_r$ with every column is either empty or else a vertical bar, i.e. a vertically connected set of squares positioned in the same column. (iv) The map $r \mapsto D_r$ is injective.

\begin{example} \label{ex:ex1}{\em (i)} Let $L=3$ and $m=6.$ The diagram $D_r$ that corresponds to the sequence  $r=((-2)^3,(-1)^2,(0)^2,(1)^3,(2)^4,(3)^4)$ is 
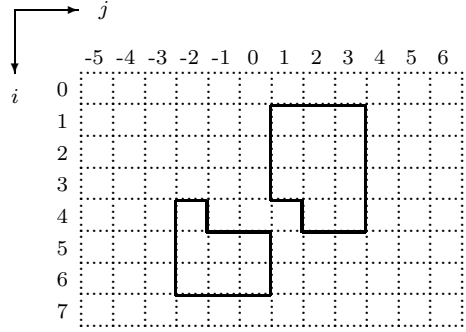
\begin{em} \begin{tiny}
\newsavebox{\hln} \newsavebox{\vln} 
\begin{center}
%%%%%%%%%%%%%%%%%%%%%%%%% picture 1
\begin{picture}(100,100)(0,0)
\put(-20,100){\vector(1,0){20}} \put(3,95){\makebox(10,10){$j$}}
\put(-20,100){\vector(0,-1){20}}  \put(-25,68){\makebox(10,10){$i$}}
\put(120,80){.}
\savebox{\hln}(120,0){\dln{2}{0}{60}{120}{0}}
\savebox{\vln}(0,80){\dln{0}{2}{40}{0}{80}}
\multiput(0,0)(0,10){9}{\usebox{\hln}}
\multiput(0,0)(10,0){13}{\usebox{\vln}}

\put(-10,70){\makebox(10,10){0}}\put(-10,60){\makebox(10,10){1}}\put(-10,50){\makebox(10,10){2}}
\put(-10,40){\makebox(10,10){3}}\put(-10,30){\makebox(10,10){4}}\put(-10,20){\makebox(10,10){5}}
\put(-10,10){\makebox(10,10){6}}\put(-10,0){\makebox(10,10){7}}

\put(0,80){\makebox(10,10){-5}}\put(10,80){\makebox(10,10){-4}}\put(20,80){\makebox(10,10){-3}}\put(30,80){\makebox(10,10){-2}}\put(40,80){\makebox(10,10){-1}}\put(50,80){\makebox(10,10){0}}\put(60,80){\makebox(10,10){1}}\put(70,80){\makebox(10,10){2}}\put(80,80){\makebox(10,10){3}}\put(90,80){\makebox(10,10){4}}\put(100,80){\makebox(10,10){5}}\put(110,80){\makebox(10,10){6}}

\put(0,0){\begin{picture}(100,100)(-0.5,0)
{\thicklines
\put(30,10){\line(0,1){30}}\put(40,30){\line(0,1){10}}\put(60,10){\line(0,1){20}}\put(60,40){\line(0,1){30}}\put(70,30){\line(0,1){10}}\put(90,30){\line(0,1){40}}

\put(30,10){\line(1,0){30}}\put(30,40){\line(1,0){10}}\put(40,30){\line(1,0){20}}
\put(60,40){\line(1,0){10}}\put(70,30){\line(1,0){20}}\put(60,70){\line(1,0){30}}
}
\end{picture}}
\sbox{\hln}{}\sbox{\vln}{}
\end{picture}
\end{center}
\end{tiny}\end{em}
{\em (ii)} Let $L=3$ and $m=5.$ The diagram $D_r$ that corresponds  to the sequence \\ $r=(-2,-1,(0)^3,(1)^4,(2)^4,(3)^2)$ is 
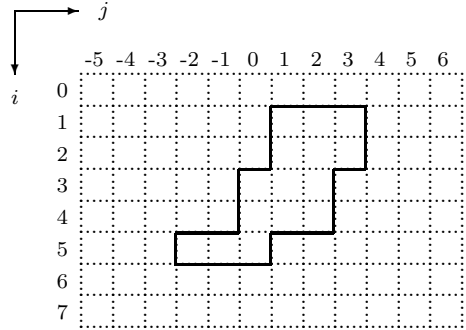
\begin{em} \begin{tiny}
\begin{center}
%%%%%%%%%%%%%%%%%%%%%%%%%%%%%%%%%%% picture 2
\begin{picture}(100,100)(0,0)
\put(-20,100){\vector(1,0){20}} \put(3,95){\makebox(10,10){$j$}}
\put(-20,100){\vector(0,-1){20}}  \put(-25,68){\makebox(10,10){$i$}}
\put(120,80){.}
\savebox{\hln}(120,0){\dln{2}{0}{60}{120}{0}}
\savebox{\vln}(0,80){\dln{0}{2}{40}{0}{80}}
\multiput(0,0)(0,10){9}{\usebox{\hln}}
\multiput(0,0)(10,0){13}{\usebox{\vln}}

\put(-10,70){\makebox(10,10){0}}\put(-10,60){\makebox(10,10){1}}\put(-10,50){\makebox(10,10){2}}
\put(-10,40){\makebox(10,10){3}}\put(-10,30){\makebox(10,10){4}}\put(-10,20){\makebox(10,10){5}}
\put(-10,10){\makebox(10,10){6}}\put(-10,0){\makebox(10,10){7}}

\put(0,80){\makebox(10,10){-5}}\put(10,80){\makebox(10,10){-4}}\put(20,80){\makebox(10,10){-3}}\put(30,80){\makebox(10,10){-2}}\put(40,80){\makebox(10,10){-1}}\put(50,80){\makebox(10,10){0}}\put(60,80){\makebox(10,10){1}}\put(70,80){\makebox(10,10){2}}\put(80,80){\makebox(10,10){3}}\put(90,80){\makebox(10,10){4}}\put(100,80){\makebox(10,10){5}}\put(110,80){\makebox(10,10){6}}

\put(0,0){\begin{picture}(100,100)(-0.5,0)
{\thicklines
\put(30,20){\line(0,1){10}}\put(50,30){\line(0,1){20}}\put(60,20){\line(0,1){10}}\put(60,50){\line(0,1){20}}\put(80,30){\line(0,1){20}}\put(90,50){\line(0,1){20}}

\put(30,20){\line(1,0){30}}\put(30,30){\line(1,0){20}}\put(50,50){\line(1,0){10}}
\put(60,30){\line(1,0){20}}\put(60,70){\line(1,0){30}}\put(80,50){\line(1,0){10}}
}
\end{picture}}
\sbox{\hln}{}\sbox{\vln}{}
\end{picture}
\end{center}
\end{tiny}\end{em}
\end{example} 
\mbox{} \hfill \\
\noindent A diagram $D \subset \Zint^2$ is called a {\em skew Young diagram} (or, simply, a skew diagram) if there exist two finite non-increasing sequences of integers: $\lambda = (\lambda_1 \geq \lambda_2 \geq \cdots \geq \lambda_l)$ and $\mu = (\mu_1 \geq \mu_2 \geq \cdots \geq \mu_l)$ such that $\lambda_i \geq \mu_i$ $(i=1,\dots,l)$ and    
\bdm 
 D = \{ (i,j) \in \Zint^2 \: | \:  1 \leq i \leq l, \: \mu_i < j \leq \lambda_i\}.
\edm
The skew Young diagram associated with a pair of sequences $\lambda$ and $\mu$ will be denoted, as usual, by $\lambda/\mu.$

\begin{defin}  \label{def:skewdiag}
We will say that $D\subset \Zint^2$ is {\em a skew diagram of type} ${\mathcal D}_L^m$ if and only if $D$ satisfies the following conditions:\\ 
{\em (i)} $D$ is a skew Young diagram with the total number of rows equal to $m.$\\
{\em (ii)} Each row of $D$ contains $L$ squares.\\
{\em (iii)} Rows of $D$ occupy the strip $\{ (i,j)\:|\: j\in \Zint, 1\leq i \leq m\}.$
\end{defin}
\noindent Below are shown two examples of skew diagrams, the one on the left is  of type ${\mathcal D}_2^6,$ the one on the right -- of type  ${\mathcal D}_4^5.$  \\
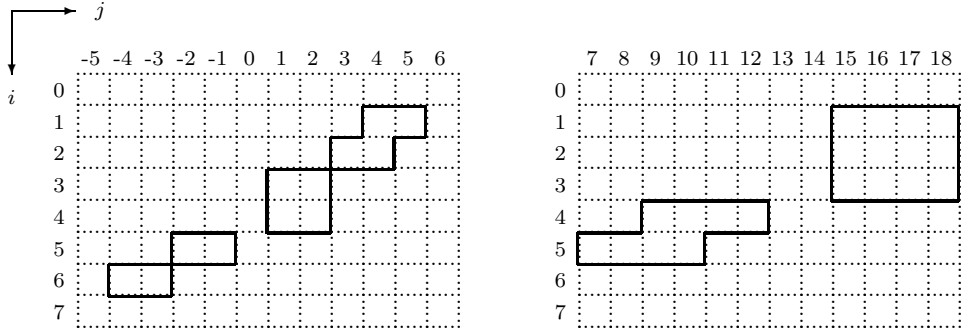
\begin{tiny}
%%%%%%%%%%%%%%%%%%%%%%%%%%%%%%%%%%%%%%% picture 3 
\begin{picture}(150,110)(-30,0)
\put(-20,100){\vector(1,0){20}} \put(3,95){\makebox(10,10){$j$}}
\put(-20,100){\vector(0,-1){20}}  \put(-25,68){\makebox(10,10){$i$}}
\put(120,80){.}
\savebox{\hln}(120,0){\dln{2}{0}{60}{120}{0}}
\savebox{\vln}(0,80){\dln{0}{2}{40}{0}{80}}
\multiput(0,0)(0,10){9}{\usebox{\hln}}
\multiput(0,0)(10,0){13}{\usebox{\vln}}

\put(-10,70){\makebox(10,10){0}}\put(-10,60){\makebox(10,10){1}}\put(-10,50){\makebox(10,10){2}}
\put(-10,40){\makebox(10,10){3}}\put(-10,30){\makebox(10,10){4}}\put(-10,20){\makebox(10,10){5}}
\put(-10,10){\makebox(10,10){6}}\put(-10,0){\makebox(10,10){7}}

\put(0,80){\makebox(10,10){-5}}\put(10,80){\makebox(10,10){-4}}\put(20,80){\makebox(10,10){-3}}\put(30,80){\makebox(10,10){-2}}\put(40,80){\makebox(10,10){-1}}\put(50,80){\makebox(10,10){0}}\put(60,80){\makebox(10,10){1}}\put(70,80){\makebox(10,10){2}}\put(80,80){\makebox(10,10){3}}\put(90,80){\makebox(10,10){4}}\put(100,80){\makebox(10,10){5}}\put(110,80){\makebox(10,10){6}}

\put(0,0){\begin{picture}(150,100)(-0.5,0)
{\thicklines
\put(30,20){\line(0,1){10}}\put(50,20){\line(0,1){10}}\put(60,30){\line(0,1){20}}\put(80,30){\line(0,1){30}}\put(90,60){\line(0,1){10}}\put(100,50){\line(0,1){10}}\put(110,60){\line(0,1){10}}

\put(30,20){\line(1,0){20}}\put(30,30){\line(1,0){20}}\put(60,30){\line(1,0){20}}
\put(60,50){\line(1,0){40}}\put(80,60){\line(1,0){10}}\put(100,60){\line(1,0){10}}\put(90,70){\line(1,0){20}}

\put(10,10){\line(0,1){10}}\put(30,10){\line(0,1){10}}
\put(10,10){\line(1,0){20}}\put(10,20){\line(1,0){20}}
}
\end{picture}}
\sbox{\hln}{}\sbox{\vln}{}
\end{picture}
\end{tiny}
%%%%%%%%%%%%%%%%%%%%%%%%%%%%%%%  picture 4
\begin{tiny}
\begin{picture}(150,110)(-30,0)
%\put(-20,100){\vector(1,0){20}} \put(3,95){\makebox(10,10){$j$}}
%\put(-20,100){\vector(0,-1){20}}  \put(-25,68){\makebox(10,10){$i$}}
\put(120,80){.}
\savebox{\hln}(120,0){\dln{2}{0}{60}{120}{0}}
\savebox{\vln}(0,80){\dln{0}{2}{40}{0}{80}}
\multiput(0,0)(0,10){9}{\usebox{\hln}}
\multiput(0,0)(10,0){13}{\usebox{\vln}}

\put(-10,70){\makebox(10,10){0}}\put(-10,60){\makebox(10,10){1}}\put(-10,50){\makebox(10,10){2}}
\put(-10,40){\makebox(10,10){3}}\put(-10,30){\makebox(10,10){4}}\put(-10,20){\makebox(10,10){5}}
\put(-10,10){\makebox(10,10){6}}\put(-10,0){\makebox(10,10){7}}

\put(0,80){\makebox(10,10){7}}\put(10,80){\makebox(10,10){8}}\put(20,80){\makebox(10,10){9}}\put(30,80){\makebox(10,10){10}}\put(40,80){\makebox(10,10){11}}\put(50,80){\makebox(10,10){12}}\put(60,80){\makebox(10,10){13}}\put(70,80){\makebox(10,10){14}}\put(80,80){\makebox(10,10){15}}\put(90,80){\makebox(10,10){16}}\put(100,80){\makebox(10,10){17}}\put(110,80){\makebox(10,10){18}}

\put(0,0){\begin{picture}(150,100)(-0.5,0)
{\thicklines
\put(0,20){\line(0,1){10}}\put(20,30){\line(0,1){10}}\put(40,20){\line(0,1){10}}\put(60,30){\line(0,1){10}}\put(80,40){\line(0,1){30}}\put(120,40){\line(0,1){30}}

\put(0,20){\line(1,0){40}}\put(0,30){\line(1,0){20}}\put(20,40){\line(1,0){40}}
\put(40,30){\line(1,0){20}}\put(80,40){\line(1,0){40}}\put(80,70){\line(1,0){40}}
}
\end{picture}}
\sbox{\hln}{}\sbox{\vln}{}
\end{picture}
\end{tiny} \\ \mbox{} \\
The diagram of Example \ref{ex:ex1}(ii) is of type ${\mathcal D}_3^5.$
\begin{propos} \label{p:regdiag}
%A sequence $r \in {\mathcal R}_0$ is regular if an only if $D_r$ is a skew diagram of type ${\mathcal D}_L^m.$ 
The map $r \mapsto D_r$ defines a one-to-one correspondence between ${\mathcal R}_0^{\reg}$ and the set of skew diagrams of type  ${\mathcal D}_L^m.$ If $r \in {\mathcal R}_0$ is not regular, then $D_r$ is not a skew Young diagram.
\end{propos}
\noindent\begin{proof} The proof is based on Corollary \ref{cor:cor2}. Let $i \in \{1,\dots,n\}$ be such that $a_i=1.$ The square of $D_r$ which corresponds to this $i$ is $(n_i,-\lambda_i L + 1).$ There exists a unique $j\in \{1,\dots,n\}$ such that $a_j = L$ and $n_j = n_i.$ \mbox{} From the last equality it follows that the square  $(n_j,-\lambda_j L + L)$ which corresponds to $j$ is located in the same row as  $(n_i,-\lambda_i L + 1).$ 

Suppose $r \in {\mathcal R}_0$ is regular, so that conditions (\ref{eq:coni}) and (\ref{eq:conii}) are satisfied. 
Condition (\ref{eq:conii}) implies that $ \lambda_j \in \{ \lambda_i -1,\lambda_i , \lambda_i +1 \}.$ Conditions (\ref{eq:coni}) and (\ref{eq:conii}) restrict possible values of $\lambda_j $ to $\lambda_i$ and $\lambda_i+1.$ For any $k\in \{1,\dots,n\}$ such that $ 1 < a_k < L$ and $n_k = n_i = n_j$ we have $|\lambda_k-\lambda_i|,|\lambda_k-\lambda_j|\leq |\lambda_i-\lambda_j|.$ Let $\lambda_j = \lambda_i,$ then $\lambda_k = \lambda_i$ and the square which corresponds to $k$ is located in the column $-\lambda_i L + a_k.$ 
Thus the $L$ squares of $D_r$ that intersect the $n_i$th row form the horizontal bar $\{ (n_i,j)\:|\: -\lambda_i L + 1 \leq j \leq -\lambda_i L + L \}.$ Let $\lambda_j = \lambda_i + 1,$ so that the square which corresponds to $j$ is located in the column $-(\lambda_i+1)L+L.$ 
In this case the square which corresponds to $k$ such that $a_k = L-1$ and $n_k=n_i=n_j$ can be located either in the column $-\lambda_i L+L-1$ or the column  $-(\lambda_i+1) L+L-1.$ If it is located in the column  $-\lambda_i L+L-1,$ then intersection of $D_r$ with the $n_i$th row is the horizontal bar $\{ (n_i,j)\:|\: -(\lambda_i +1)L + L \leq j \leq -\lambda_i L + L-1 \}.$ If it is located in the column  $-(\lambda_i+1)L+L-1,$ then we repeat the above arguments for $k$ such that $a_k = L-2$ and $n_k = n_i = n_j.$ 
Continuing by induction we find that intersection of $D_r$ with every row located in the strip $\{ (i,j)\:|\: j\in \Zint, 1\leq i \leq m\}$ is a horizontal bar of  $L$ squares. 
Suppose that $D_r$ has two rows located in two adjacent rows of $\Zint^2$ so that the leftmost square of the top row is on the left of the leftmost square of the bottom row:    
\begin{equation}\begin{picture}(100,15)(0,0) {\thicklines 
\put(0,10){\begin{picture}(50,10)(0,0) \put(0,0){\line(0,1){10}}\put(50,0){\line(0,1){10}}\put(0,10){\line(1,0){50}}\put(0,0){\line(1,0){50}}\end{picture}}
\put(40,0){\begin{picture}(50,10)(0,0) \put(0,0){\line(0,1){10}}\put(50,0){\line(0,1){10}}\put(0,10){\line(1,0){50}}\put(0,0){\line(1,0){50}}\end{picture}}
}
\put(15,5){\vector(-1,0){15}}\put(25,5){\vector(1,0){15}} \put(15,0){\makebox(10,10){{\scriptsize $ h$}}}
\end{picture} \label{eq:bconf} \quad \text{ where $h > 0.$}\end{equation}
Then there exist $i,j \in \{1,\dots,n\}$ such that $a_i=a_j,$ $\lambda_i > \lambda_j$ and $n_i = n_j - 1.$ However, for any $r \in {\mathcal R}_0$ the definition of $n_i$ implies that if $a_i=a_j,$ $\lambda_i > \lambda_j$ then $n_i > n_j.$ Thus we have a contradiction, and hence a configuration of the type (\ref{eq:bconf}) is impossible, which shows that  $D_r$ is a skew Young diagram.

Suppose now that $r \in {\mathcal R}_0$ is not regular. Let $p_i = \sum_{a=1}^L x_a^{(i)} \bar{\vartheta}_a $ be the $i$th element of the path $p(b)$ (cf. Corollary \ref{cor:cor2}). Suppose the condition (\ref{eq:coni}) is violated: there are $a\in \{2,\dots,L\}$ and $i\in \{1,\dots,n\}$ such that $x_a^{(i)} > x_{a-1}^{(i)}.$ 
Then one can always find $j$ $(j \leq i)$ such that $b_j = a$ and $x_a^{(j)} > x_{a-1}^{(j)}.$ Let $k$ be such that $j= \gamma_{\lambda}^{-1}(k)$ so that $b_j = a_k =a.$ The square of $D_r$ which  corresponds to this $k$ is located at the intersection of the column $-\lambda_k L + a$ with the row $x_a^{(j)}$ of $\Zint^2.$ 

Suppose that intersection of the column  $-\lambda_k L + a-1$ with  $D_r$ is not empty. Then it is necessarily a  vertical bar. The vertical coordinate of the bottom square of this bar is $x_{a-1}^{(j)}.$ Suppose that this bottom square is lower or on the level with the square   $(x_a^{(j)},-\lambda_k L + a),$ i.e. $x_{a-1}^{(j)} \geq   x_a^{(j)}.$ 
The last inequality contradicts $x_{a-1}^{(j)} <   x_a^{(j)},$ and therefore the bottom square of the vertical bar of $D_r$ located in the column $-\lambda_k L + a-1$ is above the bottom square of the vertical bar of $D_r$ located in the column $-\lambda_k L + a,$ i.e. $D_r$ is not a skew Young diagram. 

Suppose that intersection of the column  $-\lambda_k L + a-1$ with  $D_r$ is  empty. Then intersection of the row  $x_a^{(j)}$  with  $D_r$ is not a horizontal bar and is not empty, i.e. $D_r$ is not a skew Young diagram. 

Now suppose that the condition (\ref{eq:conii}) is violated: there exist $i,j$ such that $x_1^{(i)} = x_L^{(j)}$ and $|\lambda_i - \lambda_j| > 1.$ Then the $x_1^{(i)}$th row of $\Zint^2$ contains two squares of $D_r:$ $(x_1^{(i)},-\lambda_i L + 1 )$ and $(x_L^{(j)},-\lambda_j L + L )$ separated by horizontal distance $|-\lambda_j L + L +\lambda_i L - 1| > L.$ 
Since each row with numbers $1,\dots,m$ contains $L$ squares of $D_r,$ intersection of $D_r$ with the $x_1^{(i)}$th row is not a horizontal bar, i.e. $D_r$ is not a skew Young diagram. 

Thus it has been proven that $r$ is regular if and only if $D_r$ is a skew diagram of type ${\mathcal D}_L^m.$ 
Let $D$ be any diagram of $n$ squares. 
To this $D$ we associate the non-decreasing sequence $r = (r_1 \leq r_2 \leq \cdots \leq r_n) := r(D)$ by setting for each $s \in \Zint:$ $\#\{ r_i \:|\: r_i = s\} $ $=$  the number of squares of $D$ that belong to the $s$th column of $\Zint^2.$ 
Note that the map $D \mapsto r(D)$ is not injective on the set of all diagrams with $n$ squares.  Let $D$ be a skew diagram of type ${\mathcal D}_L^m.$ Then we  have $r(D) \in {\mathcal  R}_0.$ 
It is easy to see that   $D_{r(D)} = D.$ Hence $r(D)\in  {\mathcal  R}_0^{\reg},$ and $r \mapsto D_r$ is a surjection  between ${\mathcal  R}_0^{\reg}$ and the set of skew diagrams of type  ${\mathcal D}_L^m.$ Injectivity of this map is obvious. Thus it is a bijection  with the inverse map given by $ D \mapsto r(D).$ \end{proof}

\noindent Thus the sequence $r$ from Example \ref{ex:ex1}(i) is not regular because $D_r$ is not a skew Young diagram. 

 Let us number squares of a skew diagram in the natural order, that is column by column from left to right, and from top to bottom within every column. For example, this is the numbered skew diagram of type ${\mathcal D}_{3}^5$ from Example \ref{ex:ex1}(ii):  
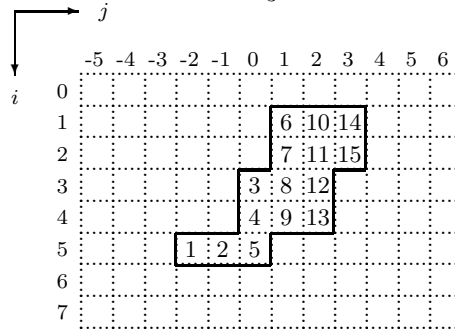
\begin{tiny}
\begin{center}
%%%%%%%%%%%%%%%%%%%%%%%%%%%%%%%%% picture 5 
\begin{picture}(100,100)(0,0)
\put(-20,100){\vector(1,0){20}} \put(3,95){\makebox(10,10){$j$}}
\put(-20,100){\vector(0,-1){20}}  \put(-25,68){\makebox(10,10){$i$}}
\put(120,80){.}
\savebox{\hln}(120,0){\dln{2}{0}{60}{120}{0}}
\savebox{\vln}(0,80){\dln{0}{2}{40}{0}{80}}
\multiput(0,0)(0,10){9}{\usebox{\hln}}
\multiput(0,0)(10,0){13}{\usebox{\vln}}

\put(-10,70){\makebox(10,10){0}}\put(-10,60){\makebox(10,10){1}}\put(-10,50){\makebox(10,10){2}}
\put(-10,40){\makebox(10,10){3}}\put(-10,30){\makebox(10,10){4}}\put(-10,20){\makebox(10,10){5}}
\put(-10,10){\makebox(10,10){6}}\put(-10,0){\makebox(10,10){7}}

\put(0,80){\makebox(10,10){-5}}\put(10,80){\makebox(10,10){-4}}\put(20,80){\makebox(10,10){-3}}\put(30,80){\makebox(10,10){-2}}\put(40,80){\makebox(10,10){-1}}\put(50,80){\makebox(10,10){0}}\put(60,80){\makebox(10,10){1}}\put(70,80){\makebox(10,10){2}}\put(80,80){\makebox(10,10){3}}\put(90,80){\makebox(10,10){4}}\put(100,80){\makebox(10,10){5}}\put(110,80){\makebox(10,10){6}}

\put(0,0){\begin{picture}(150,100)(-0.5,0)
{\thicklines
\put(30,20){\line(0,1){10}}\put(50,30){\line(0,1){20}}\put(60,20){\line(0,1){10}}\put(60,50){\line(0,1){20}}\put(80,30){\line(0,1){20}}\put(90,50){\line(0,1){20}}

\put(30,20){\line(1,0){30}}\put(30,30){\line(1,0){20}}\put(50,50){\line(1,0){10}}
\put(60,30){\line(1,0){20}}\put(60,70){\line(1,0){30}}\put(80,50){\line(1,0){10}}
}
{\scriptsize
\put(30,20){\makebox(10,10){1}}\put(40,20){\makebox(10,10){2}}
\put(50,40){\makebox(10,10){3}}\put(50,30){\makebox(10,10){4}}\put(50,20){\makebox(10,10){5}}
\put(60,60){\makebox(10,10){6}}\put(60,50){\makebox(10,10){7}}\put(60,40){\makebox(10,10){8}}\put(60,30){\makebox(10,10){9}}
\put(70,60){\makebox(10,10){10}}\put(70,50){\makebox(10,10){11}}\put(70,40){\makebox(10,10){12}}\put(70,30){\makebox(10,10){13}}
\put(80,60){\makebox(10,10){14}}\put(80,50){\makebox(10,10){15}}
}
\end{picture}}
\sbox{\hln}{}\sbox{\vln}{}
\end{picture}
\end{center}
\end{tiny}

\noindent The {\em content} of a square $(i,j) \in \Zint^2$ is defined as the difference $j-i.$ In view of Propositions \ref{p:weightofph} and \ref{p:sqint}(\ref{eq:hiterm}) the vector $\pi_{\kappa,\nu}(\ph_{x_r})\cdot v_0$ such that $r \in {\mathcal R}_0^{\reg}$ is a non-zero weight vector of $S[\tf^{\prime}] \hookrightarrow \daha$ of the weight (cf. (\ref{eq:hweight})):
\begin{equation}
\zeta_{r} = L c^{*} - \sum_{i=1}^n \left( c_i + m \right) \ep_i \quad \in (\tf^{\prime})^*, \label{eq:weight}
\end{equation}
where $c_i$ is the content of the square with the number $i$ in the skew diagram $D_r.$

\subsubsection{The isotropy subgroups} \label{sec:isotropy} Let $[\pi_{\kappa,\nu}(\ph_{x_r})\cdot v_0 ]$ be the image of $\pi_{\kappa,\nu}(\ph_{x_r})\cdot v_0$ in the quotient linear space $\Ve/\UU^{\prime}({\mathfrak b}^0)\Ve.$ This image is not zero if $r$ is a regular element of ${\mathcal R}_0,$ and is a weight vector of $S[\tf^{\prime}]$ of the same weight (\ref{eq:weight}) as $\pi_{\kappa,\nu}(\ph_{x_r})\cdot v_0 $ (cf.Proposition \ref{p:dahainv}). 
In this paragraph we fix an arbitrary regular $r \in {\mathcal R}_0$ and establish several facts about the isotropy subgroup of $\Cplx \ph$ where $\ph := [\pi_{\kappa,\nu}(\ph_{x_r})\cdot v_0 ]$ in $\bar{W}.$ Let $D:= D_r$ be the skew diagram that corresponds to $r.$ Let $D^+$ be the union of columns of height $>1$ in $D,$ and let $D^- = D\setminus D^+.$ We denote by $\bar{W}_{D^{+}}$ and  $\bar{W}_{D^{-}}$ the following subgroups of $\bar{W}:$
\begin{align*}
& \bar{W}_{D^{+}} = \langle s_i \: |\: \text{ $i$ and $i+1$ are inscribed in the same column of $D^+$} \rangle ,\\
& \bar{W}_{D^{-}} = \langle s_i \: |\: \text{ $i$ and $i+1$ are inscribed in the same row of $D^-$} \rangle. 
\end{align*}
Recall that $w \in \bar{W}$ such that $w(\ep_i) = \ep_{w(i)}$ acts on the set $\{1,\dots,n\}$ as $w : i \mapsto w(i).$ 
Then $\bar{W}_{D^{+}}$ is the subgroup of $\bar{W}$ which leaves invariant sets of numbers inscribed in  columns of $D^+,$ and which does not change any of the numbers inscribed in $D^{-}.$ On the other hand, $\bar{W}_{D^{-}}$ is the subgroup of $\bar{W}$ which leaves invariant sets of numbers inscribed in  rows of $D^-,$ 
and which does not change any of the numbers inscribed in $D^{+}.$ The subgroups $\bar{W}_{D^{+}}$ and  $\bar{W}_{D^{-}}$ are mutually commutative, let $\bar{W}_{D^{+}}\times \bar{W}_{D^{-}}$ be the subgroup of $\bar{W}$ generated by $\bar{W}_{D^{+}}$ and  $\bar{W}_{D^{-}}.$

\begin{propos} \label{p:pr1} 
Let $w \in \bar{W}_{D^{+}}.$ Then $w \cdot \ph = \ph.$ 
\end{propos}
\noindent\begin{proof} It is enough to show that $w \cdot \pi_{\kappa,\nu}(\ph_{x_r})\cdot v_0 = \pi_{\kappa,\nu}(\ph_{x_r})\cdot v_0 .$ A proof of this fact is essentially the same as the last part of the  proof of Theorem 2.4.4 in \cite{AST} and will be omitted here. \end{proof}   

\noindent Thus $\bar{W}_{D^{+}}$ is a subgroup of the isotropy group of $\Cplx \ph.$ In the next proposition we use restrictions imposed by the weight $\zeta := \zeta_r$ of $\ph$ to obtain an upper bound on the isotropy subgroup of $\Cplx \ph$ in $\bar{W}.$

\begin{propos}\label{p:pr2}
 Let $w \in \bar{W}$ and $w \cdot \ph \in \Cplx^* \ph.$ Then $w \in \bar{W}_{D^{+}}\times \bar{W}_{D^{-}}.$ 
\end{propos}
\noindent\begin{proof} Any  $w \in W$ and $\xi \in \tf^{\prime}$ satisfy the following relation in $\daha$ (cf. Proposition 1.3.4 in \cite{AST}):
\begin{equation}
w^{-1}\cdot \xi \cdot w = w^{-1}(\xi) + \sum_{\alpha \in S(w)} w(\alpha)(\xi) s_{\alpha}. \label{eq:brel}
\end{equation}
For each $w \in \bar{W}$ we define a complex number $\vep(w)$ as follows:
\begin{equation}
\vep(w) = \begin{cases} a & \text{ if $ w \cdot \ph = a \ph $ $(a \in \Cplx^*),$} \\  0  & \text{ if $ w \cdot \ph \not\in \Cplx^* \ph. $} \end{cases}
\end{equation}
Let $\zeta := \zeta_r$ so that the weight of $\ph$ is $\zeta.$ Suppose $w \in \bar{W}$ and $\vep(w) \neq 0.$ Relation (\ref{eq:brel}) gives 
\begin{equation}
\zeta - w(\zeta) = \sum_{\al \in S(w) } w(\al) \:\vep(s_{\al}). \label{eq:brel2}
\end{equation}
In the following lemma we prove the proposition in the case where $w$ is a reflection:
\begin{lemma} \label{l:reflect}
Let $\al \in \bar{R}_+$ and let $s_{\al}\cdot \ph = \pm \ph.$ Then $s_{\al} \in \bar{W}_{D^{\pm}}.$
\end{lemma}
\noindent\begin{proof} Let $\al = \al_{ij}$ $(i < j).$ Assuming $\vep(s_{\al_{ij}}) \neq 0$ the relation (\ref{eq:brel2}) gives 
\begin{align}
& \vep(s_{\al_{kj}}) = \vep(s_{\al_{ik}}) \quad ( i < k < j),\label{eq:loc1}\\ 
&c_j - c_i + \vep(s_{\al_{ij}}) = - \sum_{i < k < j} \vep(s_{\al_{kj}}), \label{eq:loc2}
\end{align}
where $c_i$ is the content of the square numbered $i$ in $D.$ Let $j=i+1,$ then $c_{i+1}-c_i + \vep(s_i) = 0.$ If $\vep(s_i) = 1$ then $i$ and $i+1$ necessarily belong to the same column of $D^+.$ If $\vep(s_i) = - 1$ then $i$ and $i+1$ necessarily belong to the same row of $D^-.$ Hence the lemma is proved for reflections of length 1. Now we continue by induction. Suppose the lemma is proved for all reflections of length less than $l(s_{\al_{ij}})>1.$

Let $\vep(s_{\al_{ij}}) = 1.$ Suppose $\exists k $ $(i < k < j)$ such that $\vep(s_{\al_{kj}}) = 1.$ Then by (\ref{eq:loc1}) we have $\vep(s_{\al_{ik}}) = 1,$ and by the induction assumption: $s_{\al_{kj}},s_{\al_{ik}} \in \bar{W}_{D^+},$ which implies that $i,k,j$ belong to the same column of $D^+.$ Hence $s_{\al_{ij}}   \in \bar{W}_{D^+}.$ Suppose $\exists k $ $(i < k < j)$ such that $\vep(s_{\al_{kj}}) = -1.$ 
Then by (\ref{eq:loc1}) we have $\vep(s_{\al_{ik}}) = -1,$ and by the induction assumption: $s_{\al_{kj}},s_{\al_{ik}} \in \bar{W}_{D^-},$ which implies that $i,k,j$ belong to the same row of $D^-,$ and hence $c_j - c_i + \vep(s_{\al_{ij}}) = j - i + 1.$ On the other hand, the right-hand side of (\ref{eq:loc2}) does not exceed $j-i-1.$ Therefore the assumption that $\exists k $ $(i < k < j)$ such that $\vep(s_{\al_{kj}}) = -1$ contradicts the  equation  (\ref{eq:loc2}). 
Now suppose $ \vep(s_{\al_{kj}}) = 0$ for all $k$ $(i < k <j).$ Then (\ref{eq:loc2}) gives $c_j - c_i + 1=0.$ Since $l(s_{\al_{ij}})>1,$  we have  $j>i+1$ and therefore $j$ can be located only in the area of $\Zint^2$ which is schematically depicted below as the set of squares inscribed with $j.$ 

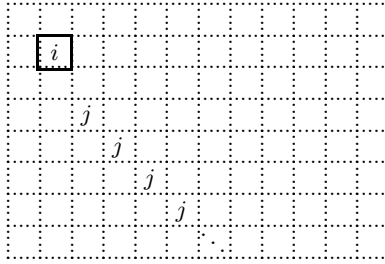
\begin{scriptsize}
\begin{center}
%%%%%%%%%%%%%%%%%%%%%%%%%%%%%%%%%  
\begin{picture}(100,80)(0,0)
\put(120,80){.}
\savebox{\hln}(120,0){\dln{2}{0}{60}{120}{0}}
\savebox{\vln}(0,80){\dln{0}{2}{40}{0}{80}}
\multiput(0,0)(0,10){9}{\usebox{\hln}}
\multiput(0,0)(10,0){13}{\usebox{\vln}}

\put(0,0){\begin{picture}(150,100)(-0.5,0)
{\thicklines
\put(10,60){\framebox(10,10){ $i$ }}
\multiput(20,40)(10,-10){4}{\makebox(10,10){ $j$ }} 
\put(60,2.5 ){\makebox(10,10){$\ddots$}}
}
\end{picture}}
\sbox{\hln}{}\sbox{\vln}{}
\end{picture}
\end{center}
\end{scriptsize}
Since $j$ is located in the diagram $D,$ the number $i+1$ must be inscribed directly below $i.$ Hence $s_i \in   \bar{W}_{D^+}$ and by Proposition \ref{p:pr1} $\vep(s_i) = 1,$ which contradicts the assumption that for all $k$ $(i<k<j)$ we have $\vep(s_{\al_{kj}}) = \vep(s_{\al_{ik}})  = 0.$

Now let $\vep(s_{\al_{ij}}) = -1.$ Suppose $\exists k $ $(i < k < j)$ such that $\vep(s_{\al_{kj}}) = -1.$ 
Then by (\ref{eq:loc1}) we have $\vep(s_{\al_{ik}}) = -1,$ and by the induction assumption: $s_{\al_{kj}},s_{\al_{ik}} \in \bar{W}_{D^-},$ which implies that $i,k,j$ belong to the same row of $D^-.$ Hence $s_{\al_{ij}}   \in \bar{W}_{D^-}.$ Suppose $\exists k $ $(i < k < j)$ such that $\vep(s_{\al_{kj}}) = 1.$ Then by (\ref{eq:loc1}) we have $\vep(s_{\al_{ik}}) = 1,$ and by the induction assumption: $s_{\al_{kj}},s_{\al_{ik}} \in \bar{W}_{D^+},$ which implies that $i,k,j$ belong to the same column of $D^+,$ and hence $c_j - c_i + \vep(s_{\al_{ij}}) = i-j - 1.$ 
On the other hand, the right-hand side of (\ref{eq:loc2}) is not less than $i - j + 1.$ Therefore the assumption that $\exists k $ $(i < k < j)$ such that $\vep(s_{\al_{kj}}) = 1$ contradicts the  equation  (\ref{eq:loc2}). Now suppose $ \vep(s_{\al_{kj}}) = 0$ for all $k$ $(i < k <j).$ Then (\ref{eq:loc2}) gives $c_j - c_i - 1=0.$ Therefore $j$ can be located only in the area of $\Zint^2$ which is schematically depicted below as the set of squares inscribed with $j.$ 
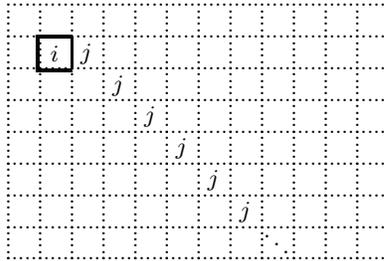
\begin{scriptsize}
\begin{center}
%%%%%%%%%%%%%%%%%%%%%%%%%%%%%%%%%  
\begin{picture}(100,80)(0,0)
\put(120,80){.}
\savebox{\hln}(120,0){\dln{2}{0}{60}{120}{0}}
\savebox{\vln}(0,80){\dln{0}{2}{40}{0}{80}}
\multiput(0,0)(0,10){9}{\usebox{\hln}}
\multiput(0,0)(10,0){13}{\usebox{\vln}}

\put(0,0){\begin{picture}(150,100)(-0.5,0)
{\thicklines
\put(10,60){\framebox(10,10){ $i$ }}
\multiput(20,60)(10,-10){6}{\makebox(10,10){ $j$ }} 
\put(80,2.5 ){\makebox(10,10){$\ddots$}}
}
\end{picture}}
\sbox{\hln}{}\sbox{\vln}{}
\end{picture}
\end{center}
\end{scriptsize}
If $j$ is located directly to the right of $i$ then either $s_i \in \bar{W}_{D^+}$ or $s_{j-1} \in \bar{W}_{D^+}$ because $j > i+1.$ If $j$ is located in any other square of the allowed area, then necessarily $s_i \in \bar{W}_{D^+}.$ Either case contradicts the assumption that  $\vep(s_{\al_{kj}})=\vep(s_{\al_{ik}}) = 0$ for all $k$ $(i < k <j).$ Thus the lemma is proved. \end{proof} 

We continue with the proof of the proposition. Define 
\begin{align*}
& \bar{R}_+^{\col} = \{ \al_{ij} \in   \bar{R}_+ \: | \: \text{ $i$ and $j$ are inscribed in the same  column of $D^+$}\} \\
& \bar{R}_+^{\row} = \{ \al_{ij} \in   \bar{R}_+ \: | \: \text{ $i$ and $j$ are inscribed in the same row of $D^-$}\}.
\end{align*}
Assuming $\vep(w) \neq 0$ the relation (\ref{eq:brel2}) and the above lemma give
\begin{equation}
w^{-1}(\zeta) - \zeta = \sum_{\{ \al \in \bar{R}_+^{\col} \: | \: w(\al) < 0 \}} \alpha \: \vep(s_{\alpha})+ \sum_{\{ \al \in \bar{R}_+^{\row} \: | \: w(\al) < 0 \} } \alpha \: \vep(s_{\alpha}). \label{eq:loc3}
\end{equation}
 Let $i$ be inscribed into a row ${\mathcal R}$ of  $D^-.$ Then the equation (\ref{eq:loc3}) gives 
\begin{equation}
c_i - c_{w(i)} = \sum_{\{ \:\fbox{{\tiny $k$}} \:\in \:{\mathcal R} \: | \: k >i,\: w(i)>w(k)\}}\vep(s_{\al_{ik}})-\sum_{\{ \:\fbox{{\tiny $k$}} \:\in \:{\mathcal R} \: | \: k <i,\: w(i)<w(k)\}}\vep(s_{\al_{ki}}) 
\end{equation}
By Lemma \ref{l:reflect} we have  $\vep(s_{\al_{ik}}), \vep(s_{\al_{ki}}) \in \{0,-1\}$ for any $k$ that appears in one of the sums in the above equation, and therefore the number $w(i)$ can be inscribed only within the infinite diagonal strip of $\Zint^2$ obtained from ${\mathcal R}$ by translations along the vector $(1,1).$ The picture below schematically shows this strip and the row ${\mathcal R}$ (in black).
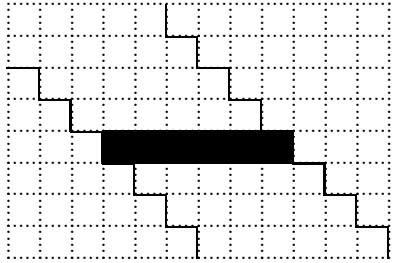
\begin{scriptsize}
\begin{center}
%%%%%%%%%%%%%%%%%%%%%%%%%%%%%%%%%  
\begin{picture}(100,80)(0,0)
\put(120,80){.}
\savebox{\hln}(120,0){\dln{2}{0}{60}{120}{0}}
\savebox{\vln}(0,80){\dln{0}{2}{40}{0}{80}}
\multiput(0,0)(0,10){9}{\usebox{\hln}}
\multiput(0,0)(10,0){13}{\usebox{\vln}}
\put(0,0){\begin{picture}(150,100)(-0.5,0)
{\thicklines
\put(30,30){\rule{72pt}{12.5pt}}}

\multiput(0,60)(10,-10){6}{\line(1,0){10}}
\multiput(50,70)(10,-10){7}{\line(1,0){10}}

\multiput(10,60)(10,-10){6}{\line(0,-1){10}}
\multiput(50,80)(10,-10){8}{\line(0,-1){10}}

\end{picture}}
\sbox{\hln}{}\sbox{\vln}{}
\end{picture}
\end{center}
\end{scriptsize}
The intersection of the strip with $D$ is exactly the row ${\mathcal R}.$ Hence $w(i)$ can be located only in ${\mathcal R}.$ Thus the set of numbers inscribed into ${\mathcal R}$   is preserved by $w.$ 
  
\noindent Let $i$ be inscribed into the leftmost column ${\mathcal C}$ of  $D^+.$ Then the equation (\ref{eq:loc3}) gives 
\begin{equation}
c_i - c_{w(i)} = \sum_{\{ \:\fbox{{\tiny $k$}} \:\in \:{\mathcal C} \: | \: k >i,\: w(i)>w(k)\}}\vep(s_{\al_{ik}})-\sum_{\{ \:\fbox{{\tiny $k$}} \:\in \:{\mathcal C} \: | \: k <i,\: w(i)<w(k)\}}\vep(s_{\al_{ki}}) \label{eq:loc31} 
\end{equation}
By Lemma \ref{l:reflect} and Proposition \ref{p:pr1} we have  $\vep(s_{\al_{ik}}), \vep(s_{\al_{ki}}) = 1$ for any $k$ that appears in one of the sums in the above equation, and therefore the number $w(i)$ can be inscribed only within the infinite diagonal strip of $\Zint^2$ obtained from ${\mathcal C}$ by translations along the vector $(1,1).$ 
The number $w(i),$ moreover, cannot be located to the left of ${\mathcal C}$ because $w$ leaves invariant sets of numbers inscribed into rows of $D^-$ (if any) positioned to the left of ${\mathcal C},$ nor can it be located below ${\mathcal C}$ because below ${\mathcal C}$ there are no squares of $D.$  
Therefore $w(i)$ can only be located in the right triangle of squares which has ${\mathcal C}$ as its left edge and has base of width equal to the height of ${\mathcal C}.$  The picture below schematically shows this triangle and, inside it, the column ${\mathcal C}$ (in black). 
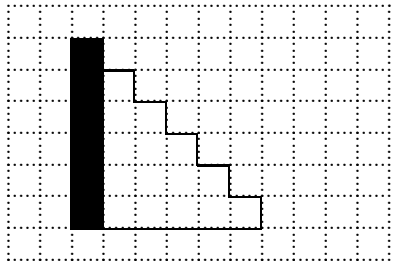
\begin{scriptsize}
\begin{center}
%%%%%%%%%%%%%%%%%%%%%%%%%%%%%%%%%  
\begin{picture}(100,80)(0,0)
\put(120,80){.}
\savebox{\hln}(120,0){\dln{2}{0}{60}{120}{0}}
\savebox{\vln}(0,80){\dln{0}{2}{40}{0}{80}}
\multiput(0,0)(0,10){9}{\usebox{\hln}}
\multiput(0,0)(10,0){13}{\usebox{\vln}}
\put(0,0){\begin{picture}(150,100)(-0.5,0)
{\thicklines
\put(20,10){\rule{12.5pt}{72pt}}}

\put(20,10){\line(1,0){60}}

\multiput(70,20)(-10,10){6}{\line(1,0){10}}
\multiput(80,10)(-10,10){6}{\line(0,1){10}}
\end{picture}}
\sbox{\hln}{}\sbox{\vln}{}
\end{picture}
\end{center}
\end{scriptsize}
Let $s \in \{1,\dots,n\}$ be such that the top square of ${\mathcal C}$ is inscribed with $s,$ without loss of generality we will assume in what follows that $s=1.$ Let $h$ be the height of ${\mathcal C},$ so that ${\mathcal C}$ is inscribed with $1,2,\dots,h$ from top to bottom. For each $i=1,\dots,h$ we have from (\ref{eq:loc31}):  
\begin{equation}
c_{i} - c_{w(i)} = \#\{h\geq k>i\:|\: w(i) > w(k)\} - \#\{1\leq k<i\:|\: w(i) < w(k)\}. \label{eq:loc32}
\end{equation}
Let now $i$ be such that $w(i)$ is maximal element of the set $\{ w(1),w(2),\dots,w(h)\}.$ 
Then the last equation gives $c_i - c_{w(i)} = h - i.$ In the triangle of allowed locations for $w(i)$ there is only one square with the content $c_i + i - h,$ and this is the bottom square of ${\mathcal C}.$ 
Thus $w(i)$ is located in this square, i.e. $w(i) = h.$ Since $w(i)$ is maximal among $w(1),w(2),\dots,w(h),$ all of the latter numbers are also located in ${\mathcal C}.$ Hence the set of numbers inscribed in ${\mathcal C}$ is left invariant by $w.$

We continue by induction. Assume that $w$ leaves invariant every  set of numbers inscribed into a column of $D^+$ located to the left of a given column ${\mathcal C}$ of  $D^+.$ Repeating the arguments given above for the case of the leftmost column of $D^{+}$ we find that $w$ leaves invariant the set of numbers inscribed in  ${\mathcal C}.$

Thus $w$ leaves invariant numbers inscribed into each column of $D^+$ and each row of $D^-.$ Hence $w \in \bar{W}_{D^+}\times \bar{W}_{D^-}.$ \end{proof}

\begin{propos} \label{p:pr3} Let $w \in \bar{W}_{D^-},$ and let $w \cdot \ph \in \Cplx^* \ph.$ Then $w \cdot \ph = (-1)^{l(w)} \ph.$ 
\end{propos}
\noindent\begin{proof} If $l(w) = 1$ the statement of the proposition follows from Lemma \ref{l:reflect}. We continue by induction. Let $l=l(w),$ $w = s_{i_1}s_{i_2}\cdots s_{i_l}$ $(s_{i_1},s_{i_2},\dots s_{i_l} \in \bar{W}_{D^-}).$ Suppose $w \cdot \ph = a \ph $ $(a \in \Cplx^*).$ The relation (\ref{eq:brel}) gives
\begin{equation}
a( w^{-1}(\zeta) - \zeta ) = \sum_{\al \in S(w)} \al \:\vep(w\cdot s_{\al}). 
\label{eq:loc4}\end{equation}
Each element of the set $\{ w\cdot s_{\al} \:|\: \al \in S(w) \}$ = $\{ s_{i_2}s_{i_3}\cdots s_{i_l},s_{i_1}s_{i_3}\cdots s_{i_l},\dots,\: s_{i_1}s_{i_2}\cdots s_{i_{l-1}} \}$ belongs to $\bar{W}_{D^-}$ and has length $( < l)$ equal to $l-1$ modulo 2. \mbox{} From (\ref{eq:loc4}) we have for any $i\in \{1,\dots,n\}:$
\begin{equation}
a(c_i - c_{w(i)}) = \sum_{\{k>i\:| w(k) < w(i)\}} \vep(w\cdot s_{\al_{ik}}) - \sum_{\{k<i\:| w(k) > w(i)\}} \vep(w\cdot s_{\al_{ki}}). \label{eq:loc5}
\end{equation}
Let $i$ be the minimal number from the set $\{ 1,\dots,n\}$ such that $w(i) \neq i.$ Then $i$ necessarily belongs to a row of $D^-$ of length $>1.$  Since $w \in \bar{W}_{D^-},$ the number $w(i)$ is located in the same row of $D^-$ as $i,$ and hence $c_i - c_{w(i)} = i - w(i).$ Then we have   
\begin{equation}
a(i - {w(i)}) = \sum_{\{k>i\:| w(k) < w(i)\}} \vep(w\cdot s_{\al_{ik}}). \label{eq:loc6}
\end{equation}
The set $\{\vep(w\cdot s_{\al_{ik}})\:| \: k > i,\:w(k) < w(i)\}$ can contain only zeroes, or, by the induction assumption, numbers $(-1)^{l-1}.$ Cardinality of this set is $w(i) - i.$ Let $m$ be the number of zeroes in this set, then (\ref{eq:loc6}) yields 
\bdm
a(i - w(i)) = (-1)^{l-1} ( w(i) - i - m).
\edm
Since $|a| = 1,$ this implies that $m=0$ and $a = (-1)^l.$ The proposition is proved. \end{proof}

\subsection{Yangian decompositions of vacuum modules of $\agl_N$}

\subsubsection{A family of Yangian highest weight vectors in the wedge product}
\label{sec:Yhwinwedge}
In Section \ref{sec:wedge} we defined the wedge product $\Vaff^{\wedge n}$ as  the quotient of  
\begin{equation}
\Vaff^{\otimes n}  \cong \Ve \otimes \left( \Cplx^N\right)^{\otimes n}
\end{equation}
by the linear subspace  $\Ker(T_1-1)+\cdots +\Ker(T_{n-1} - 1)$ where $T_i$ is the permutation of factors $i$ and $i+1$ in $\Vaff^{\otimes n}.$ 

\begin{defin} \label{def:skewdiagN}
We will say that a diagram $D \subset \Zint^2 $ of type ${\mathcal D}^m_L(N)$ if and only if $D$ is of type ${\mathcal D}_L^m$ {\em (cf. Definition \ref{def:skewdiag})} and  each  column of $D$ contains no more than $N$ squares.
\end{defin}
\noindent To every  $r \in {\mathcal R}_0^{\reg}$ such that $D_r$ is a skew diagram of type ${\mathcal D}^m_L(N)$ we associate the following vector: 
\begin{equation}
\vr_r = \uf_{p_1}\otimes \uf_{p_2} \otimes \cdots \otimes \uf_{p_n} \quad \in (\Cplx^N)^{\otimes n} \label{eq:vecr}
\end{equation}
where for every $i=1,\dots n$ we put $p_i = $ the height of the square \fra{i} in the column of $D_r$ which contains this square (thus $p_i =1$ if \fra{i} is the bottom square, etc.). Define  
\begin{equation}  \label{eq:psifin}
\psi_r^{(n)} = \wedge( \pi_{\kappa,\nu}(\ph_{x_r}) \cdot v_0\otimes \vr_r ) \quad \in \Vaff^{\wedge n}.  
\end{equation}
Recall that  we set $\kappa = L,$ $\nu(a) = a$ and assume that $n = Lm.$ 

In Section \ref{sec:wedge}, actions of $\asll_L$ were defined on $\Vaff^{\otimes n}$ and  $\Vaff^{\wedge n}.$ The subalgebra $\UU^{\prime}({\mathfrak b}^0)$ of $\UU^{\prime}(\asll_L)$ is defined in Section \ref{sec:DAHArep}, in the same section an action of $\asll_L$ is defined on $\Ve.$   

\begin{propos} \label{p:notinfin}
Let $r \in {\mathcal R}_0^{\reg}$ be such that $D_r$ is a skew diagram of type ${\mathcal D}^m_L(N).$ Then $\psi_r^{(n)} \not\in \UU^{\prime}({\mathfrak b}^0)\Vaff^{\wedge n}.$ 
\end{propos}
\noindent\begin{proof} $\bar{W}$ acts on $(\Cplx^N)^{\otimes n}$ by permutation of factors. Denote by $t_i$ the corresponding action of the $i$th simple reflection. By Proposition \ref{p:dahainv} $\bar{W}$ acts as well on the quotient space $\Ve/\UU^{\prime}({\mathfrak b}^0)\Ve.$ Denote the corresponding action of a simple reflection by $s_i.$ 
The tensor product of these two actions  gives  rise to a natural action of $\bar{W}$ on   $\Ve/\UU^{\prime}({\mathfrak b}^0)\Ve\otimes (\Cplx^N)^{\otimes n}.$ As in the previous section, we denote by  $[\pi_{\kappa,\nu}(\ph_{x_r}) \cdot v_0]$ the image of  $\pi_{\kappa,\nu}(\ph_{x_r}) \cdot v_0$ in $\Ve/\UU^{\prime}({\mathfrak b}^0)\Ve.$ The proof of the proposition is based on the following lemma.
\begin{lemma} We have:
\bdm 
\psi_r^{(n)} \in \UU^{\prime}({\mathfrak b}^0)\Vaff^{\wedge n} \qquad \Leftrightarrow \qquad [\pi_{\kappa,\nu}(\ph_{x_r}) \cdot v_0] \otimes \vr_r \in \sum_{i=1}^{n-1} \Ker(s_i\otimes t_i -1). 
\edm
\end{lemma}
\noindent\begin{proof} Observe that $\Ker(s_i\otimes t_i -1) = \Imm(s_i\otimes t_i +1)$ and $\Ker(T_i -1) = \Imm(T_i +1).$ The lemma is proved by taking into account the following elementary fact. Let $a_1,a_2,\dots$ and $b_1,b_2,\dots$ be two families of operators on a linear space ${\mathcal L},$ and let $\sum \Imm \: b_j$ be preserved by each $a_i$ and likewise  
 let $\sum \Imm \: a_j$ be preserved by each $b_i,$ so that operator $a_i,$ respectively  $b_i,$ is well defined on ${\mathcal L}/\sum \Imm \: b_j,$ respectively ${\mathcal L}/\sum \Imm \: a_j.$ For $x\in {\mathcal L}$ let $[[x]_a]_b,$ respectively $[[x]_b]_a,$ be the equivalence class of
 $x$ in the double quotient ${\mathcal L}/\sum \Imm \: a_j/\sum \Imm \: b_j,$ respectively ${\mathcal L}/\sum \Imm \: b_j/\sum \Imm \: a_j.$ Then $[[x]_a]_b = 0$ if and only if $[[x]_b]_a = 0.$ \end{proof}

\noindent Hence it is enough to prove that $A\cdot ([\pi_{\kappa,\nu}(\ph_{x_r}) \cdot v_0] \otimes \vr_r ) \neq 0$ where $
A = \sum_{w \in \bar{W}} \: (-1)^{l(w)}\: w $ 
is the total antisymmetrizer. By Proposition \ref{p:pr2} we have 
\bdm
A\cdot ([\pi_{\kappa,\nu}(\ph_{x_r}) \cdot v_0]\otimes \vr_r ) = \sum_{w \in \bar{W}_{D^+_r}\times \bar{W}_{D^-_r}} (-1)^{l(w)}\: w \:\cdot ( [\pi_{\kappa,\nu}(\ph_{x_r}) \cdot v_0]\otimes \vr_r ) \; + \;X
\edm
where $X \not\in \Cplx[\pi_{\kappa,\nu}(\ph_{x_r}) \cdot v_0] \otimes (\Cplx^N)^{\otimes n}$ or $X=0.$ By Proposition \ref{p:pr3} the first term in the right-hand side of the above equation is (cf.(\ref{eq:vecr})):
\bdm
 {\mathrm {crd}} \cdot   \sum_{w \in \bar{W}_{D^+_r}} (-1)^{l(w)}\: w \:\cdot ( [\pi_{\kappa,\nu}(\ph_{x_r}) \cdot v_0]\otimes \vr_r ) \; + \;X 
\edm
where ${\mathrm {crd}}$ is the cardinality of the isotropy subgroup of $ \Cplx [\pi_{\kappa,\nu}(\ph_{x_r}) \cdot v_0]$ in $\bar{W}_{D^-_r},$ and either $X \not\in \Cplx[\pi_{\kappa,\nu}(\ph_{x_r}) \cdot v_0] \otimes (\Cplx^N)^{\otimes n}$ or $X=0.$ Finally, by Proposition \ref{p:pr1}, the first term of the above expression is $ {\mathrm {crd}} \cdot  [\pi_{\kappa,\nu}(\ph_{x_r}) \cdot v_0]\otimes \vr_r^{\prime}$ where 
$
\vr_r^{\prime} = \sum_{w \in \bar{W}_{D^+_r}} (-1)^{l(w)}\: w \:\cdot \vr_r  \neq 0
$ 
as implied by  (\ref{eq:vecr}). \end{proof}

\noindent In Section \ref{sec:yangonwedge} we defined the Yangian action $\rho_{\kappa,\nu}(\YY(\gl_N))$ on 
the wedge product $\Vaff^{\wedge n}.$ The following proposition is shown by an elementary computation:
\begin{propos} \label{p:Yhwinwedge}
For each $r \in {\mathcal R}_0^{\reg}$ such that $D_r$ is a skew diagram of type ${\mathcal D}_L^m(N)$ the vector $\psi_r^{(n)} = \wedge(\pi_{\kappa,\nu}(\ph_{x_r}) \cdot v_0 \otimes \vr_r)$ is a highest weight vector of $\YY(\gl_N).$ We have 
\begin{align*}
&\rho_{\kappa,\nu}(T_{kl}(u))\cdot \psi_r^{(n)} = 0 \quad \text{if $ k > l,$}\\
&\rho_{\kappa,\nu}(T_{kk}(u))\cdot \psi_r^{(n)} = \prod_{i}\frac{u-m -c_i + 1}{u-m - c_i } \cdot \psi_r^{(n)}
\end{align*}
where the product is taken over $i \in \{1,\dots,n\}$ such that \fra{i} is a bottom square in a column of height $\geq k$ in $D_r.$ 
\end{propos}

\subsubsection{Yangian decompositions of   vacuum modules of $\agl_N$}\label{sec:yangdecofvac}

Let $\F_0$ be the zero charge component of the Fock space (cf. Section \ref{sec:wedge}). By Theorem \ref{t:Fockdec} we have the isomorphism of $\agl_N$-modules  $\F_0/\UU^{\prime}({\mathfrak b}^0)\F_0 \simeq  S_-\otimes V(L\Lambda_0).$ 

\noindent{\bf Semi-infinite skew diagrams.} Recall that in Section \ref{se:res} to any sequence  $\vec{h}=(h_i)_{i\in \Nat}$ such that $h_i \in \Zint_{\geq 0}$ we associated a semi-infinite skew diagram   
\begin{equation}
D(\vec{h}) = \quad 
\begin{picture}(210,20)(0,15) 
{\thicklines
\put(0,0){\line(1,0){70}}\put(0,10){\line(1,0){100}} 
\put(30,20){\line(1,0){70}}
\multiput(140,20)(0,10){2}{\line(1,0){70}}
%%%%%%%%%%
\multiput(0,0)(10,0){3}{\line(0,1){10}} \multiput(60,0)(10,0){2}{\line(0,1){10}} 
\multiput(30,10)(10,0){3}{\line(0,1){10}} \multiput(90,10)(10,0){2}{\line(0,1){10}} 
\multiput(140,20)(10,0){3}{\line(0,1){10}} \multiput(200,20)(10,0){2}{\line(0,1){10}} 
}
%%%%%%%%%%
\put(10,15){\vector(-1,0){10}}\put(20,15){\vector(1,0){10}}
\put(80,25){\vector(-1,0){50}}\put(90,25){\vector(1,0){50}}
\put(10,10){{\makebox(10,10){{\scriptsize $h_1$}}}}\put(80,20){{\makebox(10,10){ {\scriptsize $h_2$} }}}
\put(20,0){\makebox(40,10){$\dots$}}\put(50,10){\makebox(40,10){$\dots$}}
\put(160,20){\makebox(40,10){$\dots$}}
\multiput(200,35)(6,4){3}{$\cdot$}
\end{picture}
\end{equation}
\hfill \\ 
where each row has $L$ squares. Note that there are no vertical gaps between any two consecutive rows, and that we identify two skew diagrams which can be obtained one from another by translations in $\Zint^2.$ The {\em vacuum skew diagram} is defined as the skew diagram that corresponds to the sequence 
\bdm
\vec{h}^{\vac} = ( \underbrace{0,\dots,0,L}_N,\underbrace{0,\dots,0,L}_N,\dots \:).
\edm
Recall that in Definition \ref{def:siskewD} we introduced the notion of a skew Young diagram of type ${\mathcal D}_L(N).$ In what follows we will understand that ``a diagram'' means a skew diagram of type ${\mathcal D}_L(N)$ unless  stated otherwise. 

\noindent{\bf Finite parts of a semi-infinite skew diagram.} For any  diagram $D(\vec{h})$ define its {\em degree} as 
\bdm
\deg D(\vec{h})   = \sum_{i=1}^{\infty}i( h_i - \hv_i).
\edm 
The degree of any diagram is non-negative, and the only diagram of zero degree is the vacuum diagram $D(\vec{h}^{\vac}).$ Observe  that $h_i = \hv_i$ if $i > Nd$ where $d = \deg D(\vec{h}).$ 

Let $l$ be any integer such that  $l \geq d.$ Let $\ov{\!D}(\vec{h})$ be the (finite) skew diagram of type  ${\mathcal D}_L^{Nl}(N)$ (cf. Definitions \ref{def:skewdiag} and \ref{def:skewdiagN}) obtained from the semi-infinite skew diagram $D(\vec{h})$ as follows: rows of $\ov{\!D}(\vec{h})$ are the first $Nl$ rows of $D(\vec{h}).$ 
These $Nl$ rows are positioned in the plane $\Rea^2$ so that the vertical coordinate of the squares in the lowest row of $\ov{\!D}(\vec{h})$ is $Nl,$ and the horizontal coordinate of the squares in the leftmost column of  $\ov{\!D}(\vec{h})$ is $1 + \sum_{i=1}^{\infty} \hv_i - h_i.$ 
We will call $\ov{\!D}(\vec{h})$  {\em a finite part} of $D(\vec{h}).$ Note that a finite part is not uniquely defined but depends on the choice of $l \geq d.$

Let $n= lNL.$ Let $r=r_{\ov{\!D}(\vec{h})}$ be the regular sequence associated with the finite skew diagram $\ov{\!D}(\vec{h})$ (cf. Proposition \ref{p:regdiag}). That is $r=(r_1\leq r_2 \leq \cdots \leq r_n)$ where the integers $r_i$ are defined as follows: 
for each $s \in \Zint:$  $\#\{ r_i = s \}=$ the number of squares of $\ov{\!D}(\vec{h})$ in the $s$th column of $\Zint^2.$       
The vector $\psi_r^{(n)}$ (cf. Equation (\ref{eq:psifin})) belongs to the subspace $V_{n}^d$ of the wedge product $\Vaff^{\wedge n}$ (cf. Section \ref{sec:inter-and-stab}).  

\noindent{\bf A vector of the Fock space associated with a finite part of a semi-infinite skew diagram.} For each diagram $D(\vec{h})$ we fix an arbitrary  finite part $\ov{\!D}(\vec{h}),$ and  define        
\begin{equation}
\Psi_{\ov{\!D}(\vec{h})} = \vro_{l}^{d}(\psi_r^{(n)}) = \psi_r^{(n)}\wedge | - n \rangle  \quad \text{ where  $\;r = r_{\ov{\!D}(\vec{h})}.$}
\end{equation}
The vector $\Psi_{\ov{\!D}(\vec{h})} $ is a homogeneous vector of degree $d$ in the charge component $\F_0$ of the Fock space (cf. Proposition \ref{p:VFockisom}). This vector {\em a priori} depends on not only on the diagram $D(\vec{h})$ but also on the choice of the finite part $\ov{\!D}(\vec{h}).$
\begin{propos} \label{p:notin} Let $D(\vec{h})$ be a diagram and let $\ov{\!D}(\vec{h})$ be an arbitrary finite part of this diagram. We have 
\bdm  \Psi_{\ov{\!D}(\vec{h})} \not \in \UU^{\prime}({\mathfrak b}^0) \F_0. \edm 
\end{propos}
\noindent\begin{proof} Suppose  the contrary: $\Psi_{\ov{\!D}(\vec{h})} \in \UU^{\prime}({\mathfrak b}^0) \F_0.$ Then there exist vectors $v_1,\dots,v_{L-1};$ $\tilde{v}_1,\dots,\tilde{v}_{L-1}$ $\in \F_0^{d};$ $v_0 \in   \F_0^{d-1},$ where $d = \deg D(\vec{h}),$ such that  
\bdm \Psi_{\ov{\!D}(\vec{h})} = \sum_{a=1}^{L-1} h_a \cdot \tilde{v}_a + \sum_{a=1}^{L-1} f_a \cdot {v}_a + f_0 \cdot v_0.
\edm 
Let $n= lNL,$ where $l$ is the integer associated with $\ov{\!D}(\vec{h})$ as in the preceding paragraph $(l\geq d).$ By Proposition \ref{p:VFockisom} we have\begin{align}
& v_a = (\vro_l^d)^{-1}(v_a) \wedge | - n \rangle, \quad \tilde{v}_a = (\vro_l^d)^{-1}(\tilde{v}_a) \wedge | - n \rangle \quad (a=1,\dots,L-1)   \\
& v_0 = (\vro_{l-1}^{d-1})^{-1}(v_0) \wedge | - n + NL\rangle = (\vro_l^{d-1})^{-1}({v}_0) \wedge | - n \rangle. 
\end{align}
Here $(\vro_l^d)^{-1}(v_a)$ and $(\vro_l^d)^{-1}(\tilde{v}_a)$ belong to $V_n^{d},$ whereas  $(\vro_{l-1}^{d-1})^{-1}(v_0)$ belongs to $V_{n-NL}^{d-1}$ and $(\vro_{l}^{d-1})^{-1}(v_0)$ belongs to $V_{n}^{d-1}.$ Observe that $f_a\cdot |- n \rangle = h_a \cdot |- n \rangle  = 0 $ for all $a=1,\dots,L-1;$ and that  
\bdm 
u_{-n+NL}\wedge u_{-n+NL - 1} \wedge \cdots \wedge u_{-n+1} \wedge \left(f_0\cdot |-n \rangle \right) = 0. 
\edm 
This gives $ \Psi_{\ov{\!D}(\vec{h})} = \psi^{(n)}\wedge | - n\rangle $ where  
\bdm
\psi^{(n)} = \sum_{a=1}^{L-1} h_a \cdot(\vro_l^d)^{-1}(\tilde{v}_a ) + \sum_{a=1}^{L-1} f_a \cdot (\vro_l^d)^{-1}({v}_a) + f_0 \cdot (\vro_{l}^{d-1})^{-1}(v_0 ) .
\edm
The last vector belongs to $V_n^{d},$ because $\UU^{\prime}({\mathfrak b}^0)$ acting on $\Vaff^{\wedge n}$ leaves the subspace $V_n$ (cf. Equation \ref{eq:Vs}) invariant. On the other hand we have by definition: $\Psi_{\ov{\!D}(\vec{h})} = \psi_r^{(n)}\wedge | - n\rangle$ $(r = r_{\ov{\!D}(\vec{h})})$ where $\psi_r^{(n)}$ belongs to $V_n^{d}$ as well. 
Hence by Proposition \ref{p:VFockisom} we have $\psi_r^{(n)} = \psi^{(n)}.$ However, $\psi^{(n)}$ belongs to $\UU^{\prime}({\mathfrak b}^0)\Vaff^{\wedge n}$ and therefore the equality $\psi_r^{(n)} = \psi^{(n)}$ is impossible by Proposition \ref{p:notinfin}. \end{proof}

\noindent{\bf Yangian highest weight vectors in the vacuum module of $\agl_N.$} To every square in a diagram $D(\vec{h})$ we assign the {\em content}  as follows. Let $s_1$ be the leftmost  square in the first row of $D(\vec{h}).$ 
Its content is defined as $c(s_1) = 1 + \sum_{i=1}^{\infty} \hv_i - h_i.$ 
Set the origin of the coordinate system in $\Zint^2$ (cf. Section \ref{sec:regdiag}) at $s_1.$ Then the  content $c(s)$ of any other square $s=(i,j)$ is defined as $c(s) = c(s_1) + j - i.$ With  every diagram $D(\vec{h})$ we associate $N-1$ polynomials $P_1^{\vec{h}}(u),\dots, P_{N-1}^{\vec{h}}(u)$ in the formal variable $u$ as follows: 
\begin{equation}
P_k^{\vec{h}}(u) = \prod_c ( u - k - c + 1) \label{eq:DP}
\end{equation}
where the product is taken over contents of bottom squares of columns of height $k$ in the diagram $D(\vec{h}).$ Each $P_k^{\vec{h}}(u)$ is indeed a polynomial because $D(\vec{h})$ has only finite number of columns of any height less than $N.$ With  every diagram $D(\vec{h})$ we also associate the rational function 
\begin{equation}
\Omega_{\vec{h}}(u) = \prod_{k=1}^{\infty} \frac{ u + 2 - r_1  - k - \sum_{i=1}^{k-1} h_i  }{ u + 1 - k - \sum_{i=1}^{k-1} \hv_i  }
\label{eq:QD}
\end{equation}
where  $r_1 = 1 + \sum_{i=1}^{\infty} \hv_i - h_i.$ Note that $\Omega_{\vec{h}}(u)$ is indeed a rational function, because $r_1 + \sum_{i=1}^{k-1} h_i  = 1+ \sum_{i=1}^{k-1} \hv_i$ for all large enough $k.$  

Let $\ov{\!D}(\vec{h})$ be an arbitrary finite part of $D(\vec{h}),$ and let   $\left[\Psi_{\ov{\!D}(\vec{h})}\right]$ be the image of $\Psi_{\ov{\!D}(\vec{h})}$ in the quotient $\F_0/\UU^{\prime}({\mathfrak b}^0)\F_0 \simeq S_-\otimes V(L\Lambda_0).$ By Proposition  \ref{p:notin} we have $\left[\Psi_{\ov{\!D}(\vec{h})}\right] \neq 0.$

\begin{propos} \label{p:hwww} \mbox{} \\
For each semi-infinite diagram $D(\vec{h})$ let $\ov{\!D}(\vec{h})$ be an arbitrary finite part of $D(\vec{h}).$ Then $\left[\Psi_{\ov{\!D}(\vec{h})}\right]$ is a highest weight vector of the Yangian action $\rho(\YY(\gl_N))$ on $\F_0/\UU^{\prime}({\mathfrak b}^0)\F_0 $ $\simeq  $ $  S_-\otimes V(L\Lambda_0).$ We have
\begin{align*}
&\rho(\Delta(u))\cdot \left[\Psi_{\ov{\!D}(\vec{h})}\right] = \frac{\Omega_{\vec{h}}(u)}{\Omega_{\vec{h}}(u-L)}\cdot \left[\Psi_{\ov{\!D}(\vec{h})}\right],\\
&\frac{\rho(T_{k+1,k+1}(u-k))}{\rho(T_{k,k}(u-k))}\cdot \left[\Psi_{\ov{\!D}(\vec{h})}\right] = \frac{P_k^{\vec{h}}(u-1)}{P_k^{\vec{h}}(u)} \cdot \left[\Psi_{\ov{\!D}(\vec{h})}\right]\qquad (k=1,\dots,N-1).
\end{align*}
\end{propos}
\noindent\begin{proof} Definition and Proposition \ref{dp:YFock} and Proposition \ref{p:Yhwinwedge} show that $\Psi_{\ov{\!D}(\vec{h})}$ is a highest weight vector of the Yangian action $\rho_{\kappa,\nu}(\YY(\gl_N))$ on $\F_0.$ The definition of  $\rho(\YY(\gl_N))$ (cf. Section \ref{sec:Yonirr}) shows that $\left[\Psi_{\ov{\!D}(\vec{h})}\right] $ is a highest weight vector of this Yangian action. The eigenvalue of the quantum determinant $\rho(\Delta(u))$ and eigenvalues of $\rho(T_{k,k}(u))$ follow from Proposition \ref{p:Yhwinwedge}. \end{proof}

\noindent Notice that neither $\Omega_{\vec{h}}(u)$ nor polynomials $P_k^{\vec{h}}(u)$ depend on the choice of a finite part of $D(\vec{h}),$ but only on $D(\vec{h})$ itself. This is a consequence of Definition and Proposition \ref{dp:YFock}. 

\begin{propos} \label{p:lindep}
For each semi-infinite diagram $D(\vec{h})$ fix an arbitrary finite part $\ov{\!D}(\vec{h}).$ Then vectors $\left[\Psi_{\ov{\!D}(\vec{h})}\right] $ where $D(\vec{h})$ runs through the set of all semi-infinite diagrams are linearly independent.
\end{propos}
\noindent\begin{proof} Observe that for each diagram $D(\vec{h})$ the sequence $(r_1 + k + h_1 + h_2 + \cdots + h_{k-1})_{k \in \Nat}$ is strictly increasing. Therefore this sequence, and hence $D(\vec{h}),$ is uniquely restored from  the rational function $\Omega_{\vec{h}}(u).$ Hence the spectrum of the quantum determinant $\rho(\Delta(u))$ separates the vectors  $\left[\Psi_{\ov{\!D}(\vec{h})}\right] .$ \end{proof} 

\noindent{\bf Elementary modules of $\YY(\gl_N).$} 
A non-zero vector $v$ in an arbitrary $\YY(\gl_N)$-module is called a highest weight vector if it is annihilated by all the coefficients of the series $T_{kl}(u)$ with $k>l,$ and is an eigenvector of the coefficients of the series $T_{kk}(u)$ for all $k=1,\dots,N.$ Let $V$ be any irreducible finite-dimensional $\YY(\gl_N)$-module. The module $V$ contains a highest weight vector $v$ which is unique up to a scalar multiplier. Moreover,     
\begin{equation}
\frac{T_{k+1,k+1}(u-k)}{T_{k,k}(u-k)} \cdot v = \frac{P_k^V(u-1)}{P_k^V(u)}\cdot v  \qquad (k=1,\dots,N-1) \label{eq:Yhwdef}
\end{equation}
for certain monic polynomials $P_1^V(u),\dots.P^V_{N-1}(u)$ with complex coefficients. These $N-1$ polynomials are called the {\em Drinfeld polynomials} of the module $V.$ 
Two irreducible modules with the same Drinfeld polynomials may differ only by the automorphism $\omega_f$ of the algebra $\YY(\gl_N)$ which sends $T(u)$ into $f(u) T(u)$ for a formal $\Cplx$-valued series $f(u)$ of the form $f(u) = 1 + u^{-1} f_1 + u^{-2}f_2 + \cdots\:. $ 
Note that the dimension of any $\YY(\gl_N)$-module which contains a highest weight vector which  satisfies  (\ref{eq:Yhwdef}) is not less than the dimension of $V.$

Let $\lambda/\mu$ be a skew Young diagram such that every column of $\lambda/\mu$ contains no more than $N$ squares. To any such skew Young diagram and a complex number $a$ one associates an irreducible $\YY(\gl_N)$-module $V_{\lambda/\mu}(a).$ This module is defined up to an automorphism of the form $\omega_f$ by the Drinfeld polynomials
\begin{equation}
P^{V_{\lambda/\mu}(a)}_k(u) = \prod_{c}(u + c + a)\qquad (k=1,\dots,N-1)
\end{equation}
where the product is taken over contents of bottom squares of columns of height $k$ in $\lambda/\mu.$ Explicit construction for modules of the form $V_{\lambda/\mu}(a)$ is given in \cite{NT1,NT2} where they are called the {\em elementary} modules. 
For our purposes it will be sufficient to quote the formula for the $\sll_N$-character of an elementary module. Let $V$ be any finite-dimensional $\sll_N$-module. The character of $V$ is the Laurent polynomial in variables $x_1,\dots,x_N$ such that $x_1 x_2 \cdots x_N = 1$ defined as   
\bdm
\Ch_V(x) = \operatorname{Tr}\left( x_1^{e_1}x_2^{e_2}\cdots x_N^{e_N}\right)
\edm
where $e_s \in {\mathfrak h}_N$  are defined from  $H_s = e_s - e_{s+1}$ $(s=1,\dots,N)$ and $e_1 + \cdots + e_N = 0.$ %If $V$ and $V^{\prime}$ are two $\sll_N$-modules, we will write $\Ch_V(x) \geq \Ch_{V^{\prime}}(x)$ if and only if $\dim V \geq \dim V^{\prime}.$

A tableau of a skew Young diagram on numbers $1,\dots,N$ is called {\em semi-standard} if the numbers inscribed into this tableau increase strictly downward along every column and increase weakly from left to right along every row. The skew Schur function associated with a skew Young diagram $\lambda/\mu$ is defined as  
\bdm
s_{\lambda/\mu}(x) = \sum_{T} x_1^{m_1} x_2^{m_2}\cdots x_N^{m_N}
\edm 
where the sum is taken over all semi-standard tableaux of shape $\lambda/\mu$ on numbers $1,\dots,N;$ and $m_s$ is the number of squares of  $\lambda/\mu$ inscribed with $s.$ The following fact is proved in \cite{KKN1}.
\begin{propos} \label{p:smallchar}
$\qquad \qquad \Ch_{V_{\lambda/\mu}(a)}(x) = s_{\lambda/\mu}(x).$
\end{propos}
The algebra $\YY(\gl_N)$ has the following two automorphisms: ${\tau}_a$ defined for $a\in \Cplx$ by the assignment $\tau_a : T(u) \mapsto T(u+a);$ and the involution $\sigma$ defined by $\sigma : T(u) \mapsto T(-u)^{-1}.$ 
Both of these automorphisms are identical on  the subalgebra $\UU(\sll_N).$ Let $V^{\tau_a}$ and $V^{\sigma}$ be  modules obtained from an irreducible $\YY(\gl_N)$-module $V$ by a pull-back through the corresponding automorphisms. These modules are irreducible, and their Drinfeld polynomials are expressed in terms of the Drinfeld polynomials of $V$ as  
\begin{eqnarray}
& P_k^{V^{\tau_a}}(u) = P_k^V(u+a) \qquad & (k=1,\dots,N-1), \label{eq:DPofshift}\\
& P_k^{V^{\sigma}}(u) = (-1)^{\deg P_k^{V}} P_k^V(-u+k-1)\qquad & (k=1,\dots,N-1). \label{eq:DPofsigma}
\end{eqnarray}

\noindent{\bf Yangian decomposition of $S_-\otimes V(L\Lambda_0).$}For each  semi-infinite diagram $D(\vec{h})$ fix an arbitrary finite part $\ov{\!D}(\vec{h}).$ Let $V_{\ov{\!D}(\vec{h})}$ be the Yangian module generated from the highest weight vector $\left[\Psi_{\ov{\!D}(\vec{h})}\right]$ by the action $\rho(\YY(\gl_N)).$ Observe that by (\ref{eq:DPofshift}) and (\ref{eq:DPofsigma}) the  polynomial $P_k^{\vec{h}}(u)$ (cf.(\ref{eq:DP}))
 coincides with the $k$th Drinfeld polynomial of the irreducible module $V_{\ov{\!D}(\vec{h})}(a)^{\sigma}$ for a certain choice of $a$ dependent on $D(\vec{h})$ only. Therefore 
\begin{equation}
 \dim V_{\ov{\!D}(\vec{h})}  \geq \dim {V_{\ov{\!D}(\vec{h})}(a)^{\sigma}} = \dim {V_{\ov{\!D}(\vec{h})}(a)}. \label{eq:ineq}
\end{equation}

The standard grading is introduced on the Fock module $S_-  = \Cplx[B(-1),B(-2),\dots]$ of the Heisenberg algebra $H$ by setting the degree of a monomial $$ (B(-1))^{n_1}(B(-2))^{n_2} (B(-3))^{n_3} \cdots $$ to be $ n_1 + 2 n_2 + 3 n_3 + \cdots.$ And the standard grading on $V(L\Lambda_0)$ is defined by setting the degree of a monomial $ F_{s_1} F_{s_2} \cdots F_{s_n} \cdot | L \Lambda_0 \rangle $ to be equal $\# \{ i \: | \: s_i = 0 \}.$ 
A grading on the vacuum module $S_-\otimes V(L\Lambda_0)$ of $\agl_N = H \oplus \asll_N$ is defined by setting $\deg( v_1 \otimes v_2) = \deg(v_1) + \deg(v_2).$ 
In this grading each homogeneous component of  degree $d$ of $S_-\otimes V(L\Lambda_0)$ is a finite-dimensional module of $\sll_N,$ we denote this module by $(S_-\otimes V(L\Lambda_0))_d.$ The character of  $S_-\otimes V(L\Lambda_0)$ then is defined as      
\bdm
\Ch_{S_-\otimes V(L\Lambda_0)}(q,x) = \sum_{d=0}^{\infty} \: q^d \: \Ch_{(S_-\otimes V(L\Lambda_0))_d}(x),
\edm 
and the specialized character -- as $\Ch_{S_-\otimes V(L\Lambda_0)}(q,1^{N}).$ \\ \mbox{} \hfill

\noindent The eigenvalue of the quantum determinant separates between modules $V_{\ov{\!D}(\vec{h})}$ which correspond to different diagrams $D(\vec{h})$ (cf. proof of Proposition \ref{p:lindep}). 
Therefore, by (\ref{eq:ineq}) and Proposition  \ref{p:smallchar}, the specialized character of the vacuum module $S_-\otimes V(L\Lambda_0)$ dominates the expression 
\bdm
\sum_{D(\vec{h})} \: q^{\sum_{i=1}^{\infty} i(h_i - \hv_i)} \: s_{\ov{\!D}(\vec{h})}(1^N)
\edm
where the sum is taken over all skew diagrams of type ${\mathcal D}_L(N).$ On the other hand a straightforward modification of the character formula for $V(L\Lambda_0)$ proven in \cite{KKN2} gives the following character formula for $ S_-\otimes V(L\Lambda_0):$    

\begin{propos} \label{p:char} We have:
\bdm 
 \Ch_{S_-\otimes V(L\Lambda_0)}(q,x) = \sum_{D(\vec{h})} \: q^{\sum_{i=1}^{\infty} i(h_i - \hv_i)} \: s_{\ov{\!D}(\vec{h})}(x)
\edm
where the sum is taken over all skew diagrams of type ${\mathcal D}_L(N).$
\end{propos}
\noindent Hence  $\dim V_{\ov{\!D}(\vec{h})} = \dim V_{\ov{\!D}(\vec{h})}(a)^{\sigma}$  for every semi-infinite skew diagram $D(\vec{h}).$ And  $S_-\otimes V(\Lambda_0)$ is a direct sum of modules $V_{\ov{\!D}(\vec{h})}.$ The following lemma is straightforward: 
\begin{lemma} Let $V$ be an irreducible $\YY(\gl_N)$-module with Drinfeld polynomials \\ $P_1(u),\dots,P_{N-1}(u).$ Let $V^{\prime}$ be another $\YY(\gl_N)$-module containing a highest weight vector $v$ such that 
\bdm 
\frac{T_{k+1,k+1}(u-k)}{T_{k,k}(u-k)} \cdot v = \frac{P_k(u-1)}{P_k(u)}\cdot v  \qquad (k=1,\dots,N-1), 
\edm
and $\dim V^{\prime} = \dim V.$ Then $V^{\prime}$ is isomorphic to $V$ up to a $Y(\gl_N)$-automorphism of the form $\omega_f : T(u) \mapsto f(u) T(u)$ where $f(u) \in 1 + u^{-1}\Cplx[[u^{-1}]].$  
\end{lemma}
\noindent Hence, for any semi-infinite diagram $D(\vec{h})$ and any finite part $\ov{\!D}(\vec{h})$ of this diagram, the $Y(\gl_N)$-modules  $V_{\ov{\!D}(\vec{h})}$ and $V_{\ov{\!D}(\vec{h})}(a)^{\sigma}$ are isomorphic up to an automorphism of the form $\omega_f.$ 
In particular,  $V_{\ov{\!D}(\vec{h})}$  is irreducible, and its $\sll_N$-character is the skew Schur function $s_{\ov{\!D}(\vec{h})}(x).$ Note that modules  $V_{\ov{\!D}(\vec{h})}$ with different finite parts $\ov{\!D}(\vec{h})$ 
associated with the same semi-infinite diagram $D(\vec{h})$ are isomorphic, because their Drinfeld polynomials and the eigenvalues of the quantum determinant (cf.(\ref{eq:DP}),(\ref{eq:QD}) and Proposition \ref{p:hwww}) depend only on $D(\vec{h}).$ We denote an arbitrary representative of the equivalence class of $V_{\ov{\!D}(\vec{h})}$ by $V_{{D}(\vec{h})}.$  
Summarizing the working after  Proposition \ref{p:char} we obtain the irreducible decomposition of $S_-\otimes V(L\Lambda_0)$  as the Yangian module. The result is Theorem \ref{t:Yangdec}.
%------------------------------------------------------------------------------
%References
%------------------------------------------------------------------------------
\newcommand{\BOOK}[6]{\bibitem[{#6}]{#1}{#2}, {\it #3} (#4)#5.}
\newcommand{\JPAPER}[8]{\bibitem[{#8}]{#1} {#2}, {\it #3},
{#4} {\bf #5} (#6) #7.}
\newcommand{\JPAPERS}[9]{\bibitem[{#9}]{#1} {\sc #2}, `#3', {\it #4} #5 #6,
#7 #8.}

\end{document}